\documentclass[a4paper,fleqn]{cas-sc}

\usepackage[numbers]{natbib}

\def\tsc#1{\csdef{#1}{\textsc{\lowercase{#1}}\xspace}}
\tsc{WGM}
\tsc{QE}

\usepackage{amssymb}
\usepackage{latexsym}
\usepackage{amsmath}
\usepackage{mathrsfs}

\usepackage{subcaption}
\usepackage{float}

\usepackage{caption}
\usepackage{booktabs}
\usepackage{tabularx}

\usepackage{siunitx}

\DeclareMathOperator{\grad}{\nabla}
\DeclareMathOperator{\dive}{\nabla\cdot}
\DeclareMathOperator{\gradx}{\nabla_{\mathbf{x}}}
\DeclareMathOperator{\gradxi}{\nabla_{\boldsymbol{\xi}}}
\DeclareMathOperator{\divex}{\nabla_{\mathbf{x}}\cdot}
\DeclareMathOperator{\divexi}{\nabla_{\boldsymbol{\xi}}\cdot}
\DeclareMathOperator{\sgn}{\text{sgn}}

\usepackage{amsthm}
\newtheorem{theorem}{Theorem}[section]
\newtheorem{assumption}[theorem]{Assumption}

\newtheorem{definition}[theorem]{Definition}

\usepackage{placeins}
\usepackage{ulem}

\definecolor{dartmouthgreen}{rgb}{0.05, 0.5, 0.06}
\definecolor{mustard}{rgb}{0.976, 0.651, 0.008}

\begin{document}

\let\WriteBookmarks\relax
\def\floatpagepagefraction{1}
\def\textpagefraction{.001}
	
\shorttitle{AP IMEX schemes for Euler equations: non-ideal gases}    
	
\shortauthors{G. Orlando et al.}  
	
\title[mode = title]{Asymptotic-preserving IMEX schemes for the Euler equations of non-ideal gases}  
	
\author[1]{Giuseppe Orlando}[orcid=0000-0002-7119-4231]
\cormark[1]
\ead{giuseppe.orlando@polytechnique.edu}
	
\author[2]{Luca Bonaventura}[orcid=0000-0002-1994-0217]
\ead{luca.bonaventura@polimi.it}
	
\affiliation[1]{organization={CMAP, CNRS, \'{E}cole polytechnique, Institut Polytechnique de Paris},
addressline={Route de Saclay}, 
city={Palaiseau},
postcode={91120}, 
country={France}}	
	
\affiliation[2]{organization={Dipartimento di Matematica, Politecnico di Milano},
addressline={Piazza Leonardo da Vinci 32}, 
city={Milano},
postcode={20133}, 
country={Italy}}
	
\cortext[1]{Corresponding author}
	
\begin{abstract}
We analyze schemes based on a general Implicit-Explicit (IMEX) time discretization for the compressible Euler equations of gas dynamics, showing that they are asymptotic-preserving (AP) in the low Mach number limit. The analysis is carried out for a general equation of state (EOS). We consider both a single asymptotic length scale and two length scales. We then show that, when coupling these time discretizations with a Discontinuous Galerkin (DG) space discretization with appropriate fluxes, a numerical method effective for a wide range of Mach numbers is obtained. A number of benchmarks for ideal gases and their non-trivial extension to non-ideal EOS validate the performed analysis.
\end{abstract}
	
	
\begin{highlights}
    \item Analysis of the asymptotic-preserving (AP) property of a general class of IMEX schemes for a general EOS
    \item Non-trivial extension to a general EOS of the asymptotic analysis of two length scale models 
    \item Development of a high-order numerical method in combination with a DG space discretization effective for a wide range of Mach numbers
    \item Development of an AP scheme without operator splitting, flux splitting or relaxation techniques
    \item Non-trivial extension of classical benchmarks for low Mach flows to the SG-EOS and to the general cubic EOS 
\end{highlights}
	
	
\maketitle

\section{Introduction}
\label{sec:intro}

The compressible Euler equations of gas dynamics are the standard mathematical model in several applications such as atmosphere dynamics \cite{steppeler:2003}, combustion or astrophysics. For these equations, one can consider two opposite regimes. In the first one, the   acoustic waves are much faster than the local fluid velocity, while in the second one the fluid moves at high speed and compressibility plays a key role. The relevant non-dimensional number which identifies the regime is the local Mach number $M_{loc}$, defined as $M_{loc} = \frac{\left|\mathbf{u}\right|}{c}$, where $\left|\mathbf{u}\right|$ is the magnitude of the local fluid velocity and $c$ is the speed of sound. When the Mach number tends to zero, under suitable conditions, the compressible Euler equations converge to the incompressible Euler equations, see \cite{feireisl:2016, klainerman:1981}, and the references therein for the analysis of singular limits of compressible flows. Weakly compressible flows are an example of problem with multiple length and time scales. The design of efficient and stable numerical schemes for such models is a challenging task and typically requires a specific numerical treatment of the terms related to compressibility effects. 

The concept of asymptotic-preserving (AP) schemes has been introduced for this purpose, see, e.g., \cite{haack:2012}. Consider a continuous physical model $\mathcal{M}^{\varepsilon}$ which involves a small perturbation parameter $\varepsilon \ll 1$. Denote by $\mathcal{M}^{0}$ the limit of $\mathcal{M}^{\varepsilon}$ when $\varepsilon \to 0$, e.g. the incompressible Euler equations in our framework. Let now $\mathcal{M}^{\varepsilon}_{\Delta t}$ be a time discretization method which provides a consistent discretization of $\mathcal{M}^{\epsilon}$. The scheme $\mathcal{M}^{\varepsilon}_{\Delta t}$ is said to be asymptotic-preserving (AP) if its stability condition is independent of $\varepsilon$ and if its limit $\mathcal{M}^{\varepsilon}_{\Delta t}$ for $\varepsilon \to 0$ provides a consistent discretization of the continuous limit model $\mathcal{M}^{0}$. We analyze here the Euler equations of gas dynamics and the parameter $\varepsilon$ is represented by the Mach number $M$, as defined in Section \ref{sec:ap_model}. Since the seminal contribution \cite{klein:1995}, several AP schemes for Euler equations have been proposed in the literature, see among many others \cite{abbate:2019, boscheri:2020, chalons:2013, cordier:2012, dellacherie:2010a, klein:2001, kucera:2022, noelle:2014, thomann:2019} and the references therein. Methods that work at all values of the Mach number (including $M \geq 1$) are also available, see for example the seminal paper \cite{park:2005}. While a complete review of all the different approaches for low Mach flows is out of the scope of the present work, we briefly outline some of the strategies proposed in the literature to deal with low Mach flows, in order to highlight the main differences with the numerical method considered here. Following the discussion in \cite{klein:1995}, a class of AP methods \cite{chalons:2013, chalons:2016} proposes to decouple acoustic and transport phenomenon, leading to the so-called Lagrange-Projection  schemes. In these approaches, an operator splitting is applied, solving first the transport subsystem and dealing with acoustic effects afterwards. Following again \cite{klein:1995}, another class of AP schemes \cite{cordier:2012, noelle:2014} considers a splitting of the fluxes into non-stiff and stiff parts. More specifically, effects of global compression or long-wave acoustics are considered explicitly and then an implicit pressure correction is applied. Another class of popular methods are the so-called pressure correction schemes. They extend the projection techniques widely used for incompressible flows \cite{chorin:1967, orlando:2022a, temam:1969} and, starting from \cite{harlow:1968, harlow:1971}, several approaches have been proposed \cite{hennink:2021, herbin:2014, therme:2014}. Finally, a Suliciu type relaxation scheme \cite{suliciu:1990}, splitting the pressure in a slow and a fast acoustic part, was proposed in \cite{thomann:2019}, whereas a Jin-Xin type relaxation method, building a linear hyperbolic relaxation system with a small dissipative correction to approximate the Euler equations, was presented in \cite{abbate:2019}. 

We analyze here the AP properties of a general class of Implicit-Explicit (IMEX) time discretization schemes. The key observation is that, as first proposed in \cite{casulli:1984}, it suffices to adopt an implicit treatment of the pressure gradient term within the momentum equation and of the pressure work term in the energy equation to remove the acoustic CFL restriction and to decouple acoustic and transport effects, see also Appendix \ref{app:eigenvalues}. Similar approaches have been proposed, e.g., in \cite{boscheri:2020, busto:2021}. Here, we consider a general equation of state (EOS), to which only a small number of studies have been devoted \cite{abbate:2019, cordier:2012, dellacherie:2010a}. In particular, the single spatial scale analysis performed in \cite{klein:1995} was first extended to the general EOS case in \cite{dellacherie:2010a}. Here, the corresponding extension to a general EOS is introduced also for the case of two length scales. Notice that several low Mach schemes have been proposed for a barotropic equation of state \cite{bruel:2019, grenier:2013, herbin:2021}. As discussed in \cite{klein:1995}, the assumption of a barotropic fluid, for which a direct relation between the pressure and the density exists, restricts the analysis to constant-entropy data and the limit case is an incompressible flow with constant density. However, large amplitude density fluctuations are crucial for an accurate description of reacting flows \cite{klein:1995}, for atmospheric applications, and for the analysis of relevant fluid dynamics instabilities, as we will see in Section \ref{sec:num}. Finally, we show that a high-order numerical method effective for a wide range of Mach number values can be obtained coupling these time discretizations with a Discontinuous Galerkin (DG) space discretization \cite{giraldo:2020} with appropriate fluxes. In a recent work, Jung and Perrier \cite{jung:2024} analyzed the behaviour of the DG method for low Mach regimes, showing under which conditions a low Mach number accurate method is obtained. We discuss the practical implications of these results for our method, which however is shown to provide accurate results for Mach number values corresponding to fluids typically modelled as incompressible. The numerical verification is based on the higher order extension of the IMEX-DG method proposed in \cite{orlando:2023c, orlando:2022b, orlando:2023b}.

The paper is structured as follows. In Section \ref{sec:ap_model}, we present the formal limits of the continuous model considering both a single length scale and two length scales. In Section \ref{sec:ap_num}, we show the AP property of a general class of IMEX-RK methods, whereas in Section \ref{sec:space} we discuss some details of the DG formulation that allows us to obtain a numerical scheme effective for a wide range of Mach numbers. In Section \ref{sec:num}, some numerical results are presented to verify the robustness of the proposed approach with $M < 1$ and $M \ll 1$, using the higher order extension of the numerical method developed in \cite{orlando:2023c, orlando:2022b, orlando:2023b}. Finally, some conclusions and perspectives for future work are discussed in Section \ref{sec:conclu}.

\section{Asymptotic analysis for the continuous model}
\label{sec:ap_model}

Our goal is to discuss here the limit of the fully compressible Euler equations of gas dynamics as the Mach number goes to zero. For this purpose, we introduce the Euler equations and recall their non-dimensional formulation. Let $\Omega \subset \mathbb{R}^{d}, 1 \le d \le 3$ be a connected open bounded set with a sufficiently smooth boundary $\partial\Omega$ and denote by $\mathbf{x}$ the spatial coordinates and by $t$ the temporal coordinate. The mathematical model reads as follows:
\begin{eqnarray}\label{eq:euler_energy}
    \frac{\partial\rho}{\partial t} + \dive\left(\rho\mathbf{u}\right) &=& 0 \nonumber \\
    \frac{\partial\rho\mathbf{u}}{\partial t} + \dive\left(\rho\mathbf{u} \otimes \mathbf{u}\right) + \grad p &=& \mathbf{0} \\
    \frac{\partial\rho E}{\partial t} + \dive\left[\left(\rho E + p\right)\mathbf{u}\right] &=& 0. \nonumber 
\end{eqnarray}
Here $\rho$ is the density, $\mathbf{u}$ is the fluid velocity, $p$ is the pressure, and $E$ is the total energy per unit of mass. The previous set of equations has to be completed by an equation od state (EOS). Notice that no external source terms, such as gravity terms, are considered in \eqref{eq:euler_energy}. The total energy $\rho E$ can be rewritten as $\rho E = \rho e + \rho k$, where $e$ is the internal energy and $k = \frac{1}{2}\left|\mathbf{u}\right|^{2}$ is the kinetic energy per unit of mass. We also introduce the specific enthalpy $h = e + \frac{p}{\rho}$ and we notice that one can rewrite the energy flux as
\begin{equation}
    \left(\rho E + p\right)\mathbf{u} = \left(e + k + \frac{p}{\rho}\right)\rho\mathbf{u} = \left(h + k\right)\rho\mathbf{u}.
\end{equation}
Hence, \eqref{eq:euler_energy} can be rewritten as
\begin{eqnarray}\label{eq:euler}
    \frac{\partial\rho}{\partial t} + \dive\left(\rho\mathbf{u}\right) &=& 0 \nonumber \\
    \frac{\partial\rho\mathbf{u}}{\partial t} + \dive\left(\rho\mathbf{u} \otimes \mathbf{u}\right) + \grad p &=& \mathbf{0} \\
    \frac{\partial\rho E}{\partial t} + \dive\left[\left(h + k\right)\rho\mathbf{u}\right] &=& 0. \nonumber 
\end{eqnarray}
We now proceed to recall the non-dimensional version of system \eqref{eq:euler}, along the lines of the analysis presented, e.g., in \cite{klein:2001}, to which we refer for a more extensive discussion of the  underlying hypotheses. We introduce reference scaling values $\mathcal{T}, \mathcal{L}$, and $\mathcal{U}$ for time, length, and velocity, respectively. We also introduce reference values $\mathcal{P}$ for the pressure and $\mathcal{R}$ for the density. The Buckingham $\pi$ theorem \cite{buckingham:1914} states that there are $n - m$ relevant non-dimensional parameters that characterize the model, where $n$ is the number of independent physical variables and $m$ is the rank of the matrix which associates to each physical variable its unit of measure. Here, $n=5$ and $m=3$, as it can be easily verified and discussed in detail in \cite{klein:2001}. Hence, there are 2 non-dimensional parameters associated to \eqref{eq:euler}. We assume that the internal energy scales as $\mathcal{I} \approx \frac{\mathcal{P}}{\mathcal{R}}$ and that the total energy scales as $\mathcal{E} \approx \mathcal{I} + \mathcal{U}^{2}$. Finally, we assume that the specific enthalpy scales as $\mathcal{H} \approx \mathcal{I} + \frac{\mathcal{P}}{\mathcal{R}}$. We then introduce the following non-dimensional parameters
\begin{equation}
    St = \frac{\mathcal{L}}{\mathcal{T}\mathcal{U}} \qquad M^{2} = \frac{\mathcal{R}\mathcal{U}^{2}}{\mathcal{P}}
\end{equation}
and notice that
\begin{equation}
    \frac{\mathcal{I} + \frac{\mathcal{P}}{\mathcal{R}}}{\mathcal{E}} \approx  \frac{2\frac{\mathcal{P}}{\mathcal{R}}}{\frac{\mathcal{P}}{\mathcal{R}} + \mathcal{U}^{2}} = \frac{\frac{2}{M^{2}}}{\frac{1}{M^{2}} + 1} = \mathcal{O}\left(1\right) \qquad 
    \frac{U^{2}}{\mathcal{E}} \approx \frac{\mathcal{U}^{2}}{\frac{\mathcal{P}}{\mathcal{R}} + \mathcal{U}^{2}} = \frac{1}{1 + \frac{1}{M^{2}}} = \mathcal{O}\left(M^{2}\right).
\end{equation}
As a consequence, the non-dimensional version of \eqref{eq:euler} reads as follows:
\begin{eqnarray}
    St\frac{\partial\rho}{\partial t} + \dive\left(\rho\mathbf{u}\right) &=& 0 \nonumber \\
    St\frac{\partial\rho\mathbf{u}}{\partial t} + \dive\left(\rho\mathbf{u} \otimes \mathbf{u}\right) + \frac{1}{M^{2}}\grad p &=& \mathbf{0} \\
    St\frac{\partial\rho E}{\partial t} + \dive\left[\left(h + kM^{2}\right)\rho\mathbf{u}\right] &=& 0, \nonumber 
\end{eqnarray}
where, with a slight abuse of notation, the non-dimensional variables are denoted with the same symbols of the dimensional ones. Finally, as customary in the literature, see, e.g., \cite{klein:2001, munz:2003}, we assume that $St \approx 1$, so as to obtain
\begin{eqnarray}\label{eq:euler_adim}
    \frac{\partial\rho}{\partial t} + \dive\left(\rho\mathbf{u}\right) &=& 0 \nonumber \\
    \frac{\partial\rho\mathbf{u}}{\partial t} + \dive\left(\rho\mathbf{u} \otimes \mathbf{u}\right) + \frac{1}{M^{2}}\grad p &=& \mathbf{0} \\
    \frac{\partial\rho E}{\partial t} + \dive\left[\left(h + kM^{2}\right)\rho\mathbf{u}\right] &=& 0. \nonumber 
\end{eqnarray}
Our goal is to present the formal limit of the continuous model both in the case of single length scale and two length scales. Notice that the asymptotic limit for single length scale for a general EOS was already present in \cite{dellacherie:2010a}. Before achieving the proposed goal, we present the EOS that will be employed for the numerical simulations in Section \ref{sec:num}.

\subsection{The equation of state}
\label{ssec:eos}	

System \eqref{eq:euler_adim} has to be completed with an equation of state (EOS). In this work, we will consider the ideal gas law, the stiffened gas EOS (SG-EOS) \cite{metayer:2016} and the general cubic EOS \cite[p.~221]{sandler:2017}, \cite[p.~119]{vidal:2001}, even though we point out that the analyses which will be carried out in Sections \ref{ssec:ap_single_length}, \ref{ssec:ap_two_scale_length}, and \ref{sec:ap_num} are valid for a general EOS.

For an ideal gas, the equation that links together pressure, density, and internal energy is given by
\begin{equation}\label{eq:ideal_gas}
    p = \left(\gamma - 1\right)\rho e = \left(\gamma - 1\right)\left(\rho E - \frac{1}{2}M^{2}\rho\mathbf{u} \cdot \mathbf{u}\right). 
\end{equation}
Notice that \eqref{eq:ideal_gas} is valid only for a constant value of the ratio $\gamma $ between the specific heat at constant pressure and the specific heat at constant volume \cite{vidal:2001}. The analogous relation for the SG-EOS reads as follows:
\begin{equation}\label{eq:sg_eos}
    p = \left(\gamma - 1\right)\left(\rho e - \rho q_{\infty}\right) - \gamma\pi_{\infty} = \left(\gamma - 1\right)\left(\rho E - \frac{1}{2}M^{2}\rho\mathbf{u} \cdot \mathbf{u} - \rho q_{\infty}\right) - \gamma\pi_{\infty},
\end{equation} 
with $q_{\infty}$ and $\pi_{\infty}$ representing constant parameters which determine the characteristics of the fluid. Notice that for $q_{\infty} = \pi_{\infty} = 0$ in \eqref{eq:sg_eos}, we recover \eqref{eq:ideal_gas}. Finally, for the general cubic EOS the equation linking together internal energy, density and temperature is given by \cite{orlando:2022b}, \cite[p.~118]{vidal:2001}
\begin{equation}\label{eq:int_energy_general_cubic_eos}
    e = e^{\sharp}(T) + \frac{a(T) + T\frac{da}{dT}}{b}U\left(\rho, b, r_{1}, r_{2}\right).
\end{equation}
Here, $e^{\sharp}(T)$ denotes the internal energy of an ideal gas at temperature $T$, $r_{1}$ and $r_{2}$ are suitable constants, whereas the parameters $a(T), b$ determine fluid characteristics \cite{vidal:2001}. More specifically, $a(T)$ is related to intermolecular forces, while $b$, the so called co-volume, takes into account the volume occupied by the molecules. The expression of $U$ is:
\begin{equation}
    U\left(\rho, b, r_{1}, r_{2}\right) = \frac{1}{r_{1} - r_{2}}\log\left(\frac{1 - \rho b r_{1}}{1 - \rho b r_{2}}\right).
\end{equation}
Notice that, for $r_{1} \to 0$ and $r_{2} \to 0$, then $U \to -b\rho$, which corresponds to the van der Waals EOS. For $r_{1} = -1 -\sqrt{2}, r_{2} = -1 + \sqrt{2}$, we get the Peng-Robinson EOS \cite{orlando:2022b}, \cite[p.~231,p.~482]{sandler:2017}, \cite[p.~118]{vidal:2001}. For the sake of simplicity, we will assume in our numerical experiments that the coefficient $a(T)$ and the quantity $\frac{de^{\sharp}}{dT}$ are constants. We refer to \cite{orlando:2022b, orlando:2025} for the specific numerical treatment of the general cubic EOS in the more general case
$$\frac{da}{dT} \neq 0 \qquad \frac{d^2e^{\sharp}}{dT^2}\neq 0.$$
Nevertheless, we recall once more that the analyses which will be carried out in Sections \ref{ssec:ap_single_length}, \ref{ssec:ap_two_scale_length}, and \ref{sec:ap_num} are valid for a general EOS, without requiring these simplifying assumptions. Finally, the equation linking pressure, density and temperature for the general cubic EOS can be expressed as follows:
\begin{equation}\label{eq:pres_general_cubic_eos}
    p = \frac{\rho R T}{1 - \rho b} - \frac{a(T)\rho^{2}}{\left(1 - \rho b r_{1}\right)\left(1 - \rho b r_{2}\right)},
\end{equation}
with $R$ denoting the specific gas constant. We refer to \cite{cowperthwaite:1969} for a detailed discussion of the relationship between \eqref{eq:pres_general_cubic_eos} and \eqref{eq:int_energy_general_cubic_eos}. Notice that for $a = b = 0$, the equation for an ideal gas equation  $p = \rho RT$ is obtained. If $a(T)$ is constant, \eqref{eq:pres_general_cubic_eos} can be easily inverted so as to provide $T(\rho,p)$, i.e.
\begin{equation}\label{eq:T_general_cubic_eos}
    T = \frac{1 - \rho b}{R}\left(\frac{p}{\rho} + \frac{a\rho}{\left(1 - \rho br_{1}\right)\left(1 - \rho br_{2}\right)}\right).
\end{equation}
Hence, substituting \eqref{eq:T_general_cubic_eos} into \eqref{eq:int_energy_general_cubic_eos}, the equation that links internal energy, pressure, and density that we consider for our numerical simulations is the following
\begin{equation}\label{eq:general_cubic_eos}
    e = \frac{1 - \rho b}{\gamma - 1}\left(\frac{p}{\rho} + \frac{a\rho}{\left(1 - \rho br_{1}\right)\left(1 - \rho br_{2}\right)}\right) + \frac{a}{b}U\left(\rho, b, r_{1}, r_{2}\right),
\end{equation}
with $\gamma$ denoting the specific heats ratio associated to $e^{\sharp}$. We also recall here the expression of the speed of sound, which will be employed to compute the acoustic Courant number (see Section \ref{sec:num}). The speed of sound is defined for a generic equation of state as \cite{orlando:2022b, vidal:2001}:
\begin{equation}\label{eq:speed_sound}
    c^{2} = \frac{\partial p}{\partial \rho}\bigg\rvert_{s} = \frac{\frac{p}{\rho^{2}} - \frac{\partial e}{\partial\rho}}{\frac{\partial e}{\partial p}} =-\frac{\frac{\partial h}{\partial\rho}}{\frac{\partial e}{\partial p}},
\end{equation}
with $s$ denoting the specific entropy. Hence, for the ideal gas law \eqref{eq:ideal_gas} we obtain
\begin{equation}
    c^{2} = \gamma\frac{p}{\rho}.
\end{equation}
For the SG-EOS, one has instead
\begin{equation}\label{eq:c2_SG_EOS}
    c^{2} = \gamma\frac{p + \pi_{\infty}}{\rho}.
\end{equation}
Finally, the speed of sound for the general cubic EOS reads as follows:
\begin{eqnarray}\label{eq:c2_general_cubic}
    c^{2} &=& \gamma\frac{p}{\rho}\frac{1}{1 - \rho b} - \frac{a\rho}{1 - \rho b}\left(\frac{\frac{\partial U}{\partial\rho}}{b}\left(\gamma - 1\right) + \frac{1 - 2\rho b}{\left(1 - \rho br_{1}\right)\left(1 - \rho br_{2}\right)}\right) \nonumber \\
    &-& ab\rho^{2}\frac{r_{1}\left(1 - \rho br_{2}\right) + r_{2}\left(1 - \rho b r_{1}\right)}{\left(1 - \rho br_{1}\right)^{2}\left(1 - \rho br_{2}\right)^{2}},
\end{eqnarray}
with
\begin{equation}
    \frac{\partial U}{\partial\rho} = -\frac{b}{\left(1 - \rho br_{1}\right)\left(1 - \rho br_{2}\right)}.
\end{equation}
Notice once more that, \eqref{eq:c2_general_cubic} is valid only if $\frac{da}{dT} = 0$ and $\frac{de^{\sharp}}{dT}$ is constant.

\subsection{Asymptotic expansion for single length scale}
\label{ssec:ap_single_length}

In this Section, we analyze the formal limit of \eqref{eq:euler_adim} as $M \to 0$ assuming that the solution depends on a single length scale. We consider the following expansion for density, velocity, and pressure, respectively:
\begin{eqnarray}
    \rho(\mathbf{x},t) &=& \bar{\rho}(\mathbf{x},t) + M\rho^{'}(\mathbf{x},t) + M^{2}\rho^{''}(\mathbf{x},t) + \mathcal{O}(M^{3}) \label{eq:rho_expansion} \\
    \mathbf{u}(\mathbf{x},t) &=& \bar{\mathbf{u}}(\mathbf{x}, t) + M\mathbf{u}^{'}(\mathbf{x},t) + M^{2}\mathbf{u}^{''}(\mathbf{x},t) + \mathcal{O}(M^{3}) \label{eq:u_expansion} \\
    p(\mathbf{x},t) &=& \bar{p}(\mathbf{x},t) + Mp^{'}(\mathbf{x},t) + M^{2}p^{''}(\mathbf{x},t) + \mathcal{O}(M^{3}). \label{eq:p_expansion}
\end{eqnarray}
From now on, for the sake of simplicity in the notation, we omit the explicit dependence on space and time for all the variables. Substituting \eqref{eq:rho_expansion} and \eqref{eq:u_expansion} into the continuity equation in \eqref{eq:euler_adim}, the leading order term relation is
\begin{equation}\label{eq:continuity_limit}
    \frac{\partial\bar{\rho}}{\partial t} + \dive\left(\bar{\rho}\bar{\mathbf{u}}\right) = 0.
\end{equation}
For what concerns the momentum balance, the first two terms in the expansion reduce to
\begin{equation}\label{eq:momentum_limit}
    \grad\bar{p} = \mathbf{0}, \qquad \grad p^{'} = \mathbf{0},
\end{equation}
which implies that  $\bar{p},  p^{'}$ do not depend on space. Moreover, the second order term reads as follows:
\begin{equation}\label{eq:momentum_limit_second_order}
    \frac{\partial\bar{\rho}\bar{\mathbf{u}}}{\partial t} + \dive\left(\bar{\rho}\bar{\mathbf{u}} \otimes \bar{\mathbf{u}}\right) + \grad p^{''} = \mathbf{0},
\end{equation}
where $p^{''}$ represents a dynamical pressure \cite{cordier:2012, thomann:2019}, namely the standard pressure variable for incompressible flows \cite{klein:1995}. Finally, the leading order term for the energy equation is
\begin{equation}\label{eq:energy_limit}
    \frac{\partial\bar{\rho}e\left(\bar{p}, \bar{\rho}\right)}{\partial t} + \dive\left(\bar{\rho}h\left(\bar{p}, \bar{\rho}\right)\bar{\mathbf{u}}\right) = 0.
\end{equation}
Notice that here we do not assume a Hilbert expansion for the internal energy $e$, and that $e\left(\bar{p}, \bar{\rho}\right)$ and $h\left(\bar{p}, \bar{\rho}\right)$ denote the expressions obtained from the equation of state evaluated at $\bar{p}, \bar{\rho}$. Other contributions in the literature, such as \cite{kucera:2022}, assume a Hilbert expansion also for the energy. The limit model obtained is the same in the case of single length scale, provided that $\frac{\overline{\rho e}}{\bar{\rho}} = e\left(\bar{p}, \bar{\rho}\right)$, whereas some differences can arise in the case of the two length scale model. Since $\bar{\rho}e\left(\bar{p}, \bar{\rho}\right) = \bar{\rho}h\left(\bar{p}, \bar{\rho}\right) - \bar{p}$, we obtain
\begin{equation}
    \frac{\partial\bar{\rho}h\left(\bar{p}, \bar{\rho}\right)}{\partial t} - \frac{\partial\bar{p}}{\partial t} + \dive\left(\bar{\rho}h\left(\bar{p}, \bar{\rho}\right)\bar{\mathbf{u}}\right) = 0,
\end{equation}
or, equivalently, thanks to \eqref{eq:continuity_limit}
\begin{equation}\label{eq:limit_model_tmp}
    \bar{\rho}\left(\frac{\partial h\left(\bar{p}, \bar{\rho}\right)}{\partial t} + \bar{\mathbf{u}} \cdot \grad h\left(\bar{p}, \bar{\rho}\right)\right) - \frac{\partial\bar{p}}{\partial t} = 0.
\end{equation}
From \eqref{eq:limit_model_tmp}, we get
\begin{equation}
    \bar{\rho}\frac{\partial h\left(\bar{p}, \bar{\rho}\right)}{\partial\bar{\rho}}\left(\frac{\partial\bar{\rho}}{\partial t} + \bar{\mathbf{u}} \cdot \grad\bar{\rho}\right) + \bar{\rho}\frac{\partial h\left(\bar{p}, \bar{\rho}\right)}{\partial\bar{p}}\left(\frac{\partial\bar{p}}{\partial t} + \bar{\mathbf{u}} \cdot \grad\bar{p}\right) - \frac{\partial\bar{p}}{\partial t} = 0.
\end{equation}
Thanks to \eqref{eq:continuity_limit} and \eqref{eq:momentum_limit}, we obtain
\begin{equation}\label{eq:energy_limit_incomp}
    -\bar{\rho}^{2}\frac{\partial h\left(\bar{p}, \bar{\rho}\right)}{\partial\bar{\rho}}\left(\dive\bar{\mathbf{u}}\right) + \bar{\rho}\frac{\partial e\left(\bar{p}, \bar{\rho}\right)}{\partial\bar{p}}\frac{d\bar{p}}{dt} = 0. 
\end{equation}
If $\bar{\rho} \neq 0$ and $\frac{\partial h\left(\bar{p}, \bar{\rho}\right)}{\partial\bar{\rho}} \neq 0$, as it holds away from vacuum, thanks to \eqref{eq:speed_sound}, relation \eqref{eq:energy_limit_incomp} can be rewritten as
\begin{equation}\label{eq:dive_u}
    \dive\bar{\mathbf{u}} = -\frac{1}{\bar{\rho}c^{2}\left(\bar{p}, \bar{\rho}\right)}\frac{d\bar{p}}{dt}
\end{equation}
Summing up, the asymptotic limit of \eqref{eq:euler_adim} is
\begin{eqnarray}\label{eq:euler_adim_ap}
    \frac{\partial\bar{\rho}}{\partial t} + \dive\left(\bar{\rho}\bar{\mathbf{u}}\right) &=& 0 \nonumber \\
    \grad\bar{p} &=& \mathbf{0} \nonumber \\
    \grad p^{'} &=& \mathbf{0} \\
    \frac{\partial\bar{\rho}\bar{\mathbf{u}}}{\partial t} + \dive\left(\bar{\rho}\bar{\mathbf{u}} \otimes \bar{\mathbf{u}}\right) + \grad p^{''} &=& \mathbf{0} \nonumber \\
    \frac{\partial\bar{\rho}e\left(\bar{p}, \bar{\rho}\right)}{\partial t} + \dive\left(\bar{\rho}h\left(\bar{p}, \bar{\rho}\right)\bar{\mathbf{u}}\right) &=& 0, \nonumber
\end{eqnarray}
or, equivalently,
\begin{eqnarray}\label{eq:euler_adim_ap_incomp}
    \frac{\partial\bar{\rho}}{\partial t} + \dive\left(\bar{\rho}\bar{\mathbf{u}}\right) &=& 0 \nonumber \\
    \grad\bar{p} &=& \mathbf{0} \nonumber \\
    \grad p^{'} &=& \mathbf{0} \\
    \frac{\partial\bar{\rho}\bar{\mathbf{u}}}{\partial t} + \dive\left(\bar{\rho}\bar{\mathbf{u}} \otimes \bar{\mathbf{u}}\right) + \grad p^{''} &=& \mathbf{0} \nonumber \\
    \dive\bar{\mathbf{u}} &=& -\frac{1}{\bar{\rho}c^{2}\left(\bar{p}, \bar{\rho}\right)}\frac{d\bar{p}}{dt} 0. \nonumber
\end{eqnarray}
The asymptotic limit \eqref{eq:euler_adim_ap_incomp} was already present in \cite{dellacherie:2010a} and represents the extension to non-ideal gases of the system of equations derived in \cite{klein:1995}. Analogous relations have been derived in \cite{cordier:2012} for the case $\frac{\partial\bar{p}}{\partial t} = 0$. Under periodic or free-slip boundary conditions, thanks to the divergence theorem, we have
$$\int_{\Omega}\dive\bar{\mathbf{u}}d\Omega = 0,$$
so that, by integrating \eqref{eq:dive_u} on $\Omega$, we find $\frac{d\bar{p}}{dt} = 0$. However, as is evident from the last relation in \eqref{eq:dive_u}, a time dependent pressure with large amplitude variations imposed by Dirichlet boundary conditions leads to a non-incompressible flow, i.e. $\dive\bar{\mathbf{u}} \neq 0$, as we will verify numerically in Section \ref{ssec:open_tube}. Hence, under periodic or free-slip boundary conditions or if $\frac{d\bar{p}}{dt} = 0$, all the equations of state lead to the same limit, namely the incompressible Euler equations. On the other hand, if $\frac{d\bar{p}}{dt} \neq 0$, then $\dive\bar{\mathbf{u}}$ depends on the specific EOS and on its parameters. For the ideal gas law \eqref{eq:ideal_gas}, we obtain
\begin{equation}\label{eq:dive_u_IG}
    \dive\bar{\mathbf{u}} = -\frac{1}{\gamma}\frac{d\log\bar{p}}{dt}.
\end{equation}
Hence, the compressibility of a fluid described by the ideal gas law is uniform in space and changes only in time. This is no longer valid for a general EOS, as we will also show in Section \ref{sec:num}.

\subsection{Asymptotic expansion for two length scales}
\label{ssec:ap_two_scale_length}

In this Section, following \cite{barsukow:2021, klein:1995}, 
we try to account for the fact that, for sufficiently small values of the Mach number, two decoupled spatial scales can be identified. More specifically, since the speed of sound $c$ is much larger than the typical flow velocity $\left|\mathbf{u}\right|$ and if a unique time scale is considered, the typical length scale associated to acoustic phenomena is much larger than that associated with the material flow. In order to properly highlight this fact, we assume that the solution depends on the material scale variable $\mathbf{x}$ and also on the acoustic scale variable $\boldsymbol{\xi} = M\mathbf{x}$. Separate equations will then be derived for the material information, which moves at speed $\left|\mathbf{u}\right|$, and for the acoustic information, which moves approximately at the speed of sound $c$ \cite{bruel:2019, galie:2024}. Relevant applications which show the interaction between the two scales arise in reacting flows \cite{klein:2002}, in the interaction of shocks with large density gradients \cite{brouillette:2002} and in atmospheric models, as we will show in Section \ref{ssec:baroclinic}. An analogous analysis can be performed considering a single length scale and two time scales, as done, e.g., in \cite{bruel:2019, galie:2024}. In an asymptotic analysis with two spatial scales, we consider the following expansion for any dependent variable:
\begin{equation}
    f\left(\mathbf{x}, \boldsymbol{\xi}, t\right) = \bar{f}\left(\mathbf{x}, \boldsymbol{\xi}, t\right) + Mf^{'}\left(\mathbf{x}, \boldsymbol{\xi}, t\right) + M^{2}f^{''}\left(\mathbf{x}, \boldsymbol{\xi}, t\right) + \mathcal{O}(M^{3}),
\end{equation} 
so that a large scale spatial derivative operator appears in the asymptotic expansion. More specifically, we get
\begin{equation}\label{eq:gradient_two_scale_length}
    \grad f = \gradx f + M\gradxi f.
\end{equation}
One can easily notice from \eqref{eq:gradient_two_scale_length} that the leading order relations are not modified introducing $\boldsymbol{\xi}$, provided that we reinterpret $\grad\square$ and $\dive\square$ as $\gradx\square$ and $\divex\square$, respectively. Equations \eqref{eq:momentum_limit} and \eqref{eq:momentum_limit_second_order} change because of $\boldsymbol{\xi}$. Indeed, since
\begin{equation}
    \grad p = \gradx\bar{p} + M\left(\gradx p^{'} + \gradxi\bar{p}\right) + M^{2}\left(\gradx p^{''} + \gradxi p^{'}\right) + \mathcal{O}(M^{3}),
\end{equation}
we obtain
\begin{eqnarray}
    \gradx p^{'} + \gradxi\bar{p} &=& \mathbf{0} \label{eq:momentum_limit_first_order_two_scale} \\
    \frac{\partial\bar{\rho}\bar{\mathbf{u}}}{\partial t} + \dive\left(\bar{\rho}\bar{\mathbf{u}} \otimes \bar{\mathbf{u}}\right) + \gradx p^{''} + \gradxi p^{'} &=& \mathbf{0}. \label{eq:momentum_limit_second_order_two_scale}
\end{eqnarray}
We also consider the first order term of the continuity equation, which reduces to
\begin{equation}
    \frac{\partial\rho^{'}}{\partial t} + \divex\left(\rho^{'}\bar{\mathbf{u}}\right) + \divex\left(\bar{\rho}\mathbf{u}^{'}\right) + \divexi \left(\bar{\rho}\bar{\mathbf{u}}\right) = 0.
\end{equation}
Finally, we consider the first order term of the energy equation, which reads as follows:
\begin{equation}\label{eq:energy_limit_first_order_two_scale}
    \frac{\partial\rho^{'}e\left(p^{'},\rho^{'}\right)}{\partial t} + \divex\left(\rho^{'}h\left(p^{'},\rho^{'}\right)\bar{\mathbf{u}}\right) + \divex\left(\bar{\rho}h\left(\bar{p},\bar{\rho}\right)\mathbf{u}^{'}\right) + \divexi\left(\bar{\rho}\bar{\mathbf{u}}h\left(\bar{p},\bar{\rho}\right)\right) = 0.
\end{equation}
Notice that, relation \eqref{eq:energy_limit_first_order_two_scale} implicitly assumes that a Hilbert expansion holds for $\rho e$, so that the first order contribution for $\rho e$ reduces to $\rho^{'}{e}\left(p^{'},\rho^{'}\right)$. However, other options are possible; as an example, assuming a Hilbert expansion for $e$ would lead to
\begin{equation}
    \frac{\partial\left(\bar{\rho}e^{'} + \rho^{'}\bar{e}\right)}{\partial t} + \divex\left(\rho^{'}\bar{e}\bar{\mathbf{u}}\right) + \divex\left(\bar{\rho}e^{'}\bar{\mathbf{u}}\right) + \divex\left(\bar{\rho}\bar{e}\mathbf{u}^{'}\right) + \divex\left(\bar{p}\mathbf{u}^{'}\right) + \divex\left(p^{'}\bar{\mathbf{u}}\right) + \divexi\left(\bar{\rho}\bar{\mathbf{u}}\bar{h}\right) = 0.
\end{equation}
Summing up, the asymptotic limit of \eqref{eq:euler_adim} for a two-scale analysis is
\begin{eqnarray}\label{eq:euler_adim_ap_two_scale}
    \frac{\partial\bar{\rho}}{\partial t} + \divex\left(\bar{\rho}\bar{\mathbf{u}}\right) &=& 0 \nonumber \\
    \gradx\bar{p} &=& \mathbf{0} \nonumber \\
    \gradx p^{'} + \gradxi\bar{p} &=& \mathbf{0} \\
    \frac{\partial\bar{\rho}\bar{\mathbf{u}}}{\partial t} + \divex\left(\bar{\rho}\bar{\mathbf{u}} \otimes \bar{\mathbf{u}}\right) + \gradx p^{''} + \gradxi p^{'} &=& \mathbf{0} \nonumber \\
    \frac{\partial\bar{\rho}e\left(\bar{p},\bar{\rho}\right)}{\partial t} + \divex\left(\bar{\rho}h\left(\bar{p},\bar{\rho}\right)\bar{\mathbf{u}}\right) &=& 0 \nonumber \\
    \frac{\partial\rho^{'}}{\partial t} + \divex\left(\rho^{'}\bar{\mathbf{u}}\right) + \divex\left(\bar{\rho}\mathbf{u}^{'}\right) + \divexi \left(\bar{\rho}\bar{\mathbf{u}}\right) &=& 0 \nonumber \\
    \frac{\partial\rho^{'}e\left(p^{'},\rho^{'}\right)}{\partial t} + \divex\left(\rho^{'}h\left(p^{'},\rho^{'}\right)\bar{\mathbf{u}}\right) + \divex\left(\bar{\rho}h\left(\bar{p},\bar{\rho}\right)\mathbf{u}^{'}\right) + \divexi\left(\bar{\rho}\bar{\mathbf{u}}h\left(\bar{p},\bar{\rho}\right)\right)&=& 0. \nonumber
\end{eqnarray}

Following the discussion in \cite{klein:1995}, we then focus on the regime in which variations on the material scale are negligible and only variations on the large acoustic scale are relevant. Starting from \eqref{eq:euler_adim_ap_two_scale}, these assumptions imply that
\begin{eqnarray}\label{eq:euler_adim_ap_two_scale_averaged}
    \frac{\partial\bar{\rho}}{\partial t} &=& 0 \nonumber \\
    \frac{\partial\bar{\rho}\bar{\mathbf{u}}}{\partial t} + \gradxi p^{'} &=& \mathbf{0} \nonumber \\
    \gradxi\bar{p} &=& \mathbf{0} \\
    \frac{\partial\bar{p}}{\partial t} &=& 0 \nonumber \\
    \frac{\partial\rho^{'}}{\partial t} + \divexi\left(\bar{\rho}\bar{\mathbf{u}}\right) &=& 0 \nonumber \\
    \frac{\partial\rho^{'}e\left(p^{'},\rho^{'}\right)}{\partial t} + \divexi\left(\bar{\rho}\bar{\mathbf{u}}h\left(\bar{p},\bar{\rho}\right)\right) &=& 0. \nonumber
\end{eqnarray}
The relation $\frac{\partial\bar{p}}{\partial t} = 0$ is a direct consequence of the fact that $\divex\bar{\mathbf{u}} = 0$, since we neglect variations on the material scale. Moreover, we notice that $\bar{p}$ reduces to a constant. In the particular case of the ideal gas law \eqref{eq:ideal_gas}, system \eqref{eq:euler_adim_ap_two_scale_averaged} reduces to
\begin{eqnarray}\label{eq:euler_adim_ap_two_scale_averaged_IG}
    \frac{\partial\bar{\rho}}{\partial t} &=& 0 \nonumber \\
    \frac{\partial\bar{\mathbf{u}}}{\partial t} + \frac{1}{\bar{\rho}\left(\boldsymbol{\xi}\right)}\gradxi p^{'} &=& \mathbf{0} \nonumber \\
    \gradxi\bar{p} &=& \mathbf{0} \\
    \frac{\partial\rho^{'}}{\partial t} + \divexi\left(\bar{\rho}\bar{\mathbf{u}}\right) &=& 0 \nonumber \\
    \frac{\partial\bar{p}}{\partial t} &=& 0 \nonumber \\
    \frac{\partial p^{'}}{\partial t} + \gamma\bar{p}\divexi\bar{\mathbf{u}} &=& 0. \nonumber
\end{eqnarray}
Taking the time derivative of the last equation, we obtain 
\begin{equation}\label{eq:wave_acoustic}
    \frac{\partial^{2}p^{'}}{\partial t^{2}} = \divexi\left(c\left(\bar{p}, \bar{\rho}\right)^{2}\gradxi p^{'}\right),
\end{equation}
with $c\left(\bar{p}, \bar{\rho}\right)^{2} = \gamma\frac{\bar{p}}{\bar{\rho}}$. Equation \eqref{eq:wave_acoustic} is the wave equation for $p^{'}$ already derived in \cite{klein:1995}. The time derivative of the first order term of the energy equation in \eqref{eq:euler_adim_ap_two_scale_averaged} reduces to
\begin{equation}\label{eq:energy_limit_first_order_two_scale_time_derivative}
    \frac{\partial^{2}\rho^{'}e\left(\rho^{'},p^{'}\right)}{\partial t^{2}} = \divexi\left(h\left(\bar{p},\bar{\rho}\right)\gradxi p^{'}\right).
\end{equation}
Starting from \eqref{eq:energy_limit_first_order_two_scale_time_derivative}, one can verify that \eqref{eq:wave_acoustic} is valid also for the SG-EOS \eqref{eq:sg_eos}. Indeed, since
$$\frac{\partial^{2}\rho^{'}e\left(\rho^{'},p^{'}\right)}{\partial t^{2}} = \frac{1}{\gamma - 1}\frac{\partial^{2} p^{'}}{\partial t^{2}} + q_{\infty}\frac{\partial^{2}\rho^{'}}{\partial t^{2}}$$
and
$$\frac{\partial^{2}\rho^{'}}{\partial t^{2}} = \divexi\gradxi p^{'},$$
we obtain
\begin{equation}
    \frac{\partial^{2}p^{'}}{\partial t^{2}} = \divexi\left[\left(\gamma - 1\right)\left(h\left(\bar{p},\bar{\rho}\right) - q_{\infty}\right)\gradxi p^{'}\right].
\end{equation}
Since
$$h\left(\bar{p},\bar{\rho}\right) = \frac{\gamma\left(\bar{p} + \pi_{\infty}\right)}{\bar{\rho}\left(\gamma - 1\right)} + q_{\infty},$$
we recover relation \eqref{eq:wave_acoustic} thanks to \eqref{eq:c2_SG_EOS}. Relation \eqref{eq:wave_acoustic} is instead in general no longer valid for a general EOS and supplementary terms arise for the general cubic EOS \eqref{eq:general_cubic_eos}.

\section{Asymptotic analysis for a class of IMEX-RK schemes}
\label{sec:ap_num}

We analyze now the AP property of a general class of Implicit-Explicit Runge-Kutta (IMEX-RK) methods for the time discretization of system \eqref{eq:euler_adim}. Following \cite{casulli:1984, dumbser:2016b}, we couple implicitly the energy equation to the momentum one, while the continuity equation is discretized in a fully explicit fashion. As a result, at each stage of the IMEX-RK method, we will obtain a nonlinear Helmholtz equation for the pressure, which is solved through a fixed point procedure \cite{dumbser:2016b, orlando:2022b}. The time discretization is based on an IMEX-RK scheme \cite{kennedy:2003}, as done in \cite{orlando:2023c, orlando:2022b, orlando:2023b}. IMEX-RK schemes are represented compactly by the companion Butcher tableaux \cite{butcher:2008}:
\begin{center}
    \begin{tabular}{c|c}
	$\mathbf{c}$ & $\mathbf{A}$ \\
	\hline
	& $\mathbf{b}^{T}$
    \end{tabular}
    \qquad
    \begin{tabular}{c|c}
	$\tilde{\mathbf{c}}$ & $\tilde{\mathbf{A}}$  \\
	\hline
	& $\tilde{\mathbf{b}}^{T}$
    \end{tabular}
\end{center}
with $\mathbf{A} = \left\{a_{lm}\right\}, \mathbf{b} = \left\{b_{l}\right\}, \mathbf{c} = \left\{c_{l}\right\}, \tilde{\mathbf{A}} = \left\{\tilde{a}_{lm}\right\}, \tilde{\mathbf{b}} = \left\{\tilde{b}_{l}\right\}$, and $\tilde{\mathbf{c}} = \left\{\tilde{c}_{l}\right\}$, $l,m = 1\dots s$, where $s$ denotes the number of stages of the method. Notice that matrix $\mathbf{A}$ corresponds to the explicit part of the scheme, i.e. $a_{i,j} = 0 $ for $j\geq i$, while $\tilde{\mathbf{A}}$ corresponds to the implicit part of the scheme. Coefficients $a_{lm}, \tilde{a}_{lm}, c_{l}, \tilde{c}_{l}, b_{l}$, and $\tilde{b}_{l}$ are determined so that the method is consistent of a given order. In particular, the following relation has to be satisfied \cite{kennedy:2003}:
\begin{equation}
    \sum_{l=1}^{s}b_{l} =  \sum_{l=1}^{s}\tilde{b}_{l} = 1.
\end{equation}
We then introduce the following Definition, which characterizes different IMEX-RK schemes according to the structure of the implicit method:
\begin{definition}
    An IMEX-RK method is said to be of \textbf{type I} \cite{boscarino:2024, pareschi:2005} if the matrix $\tilde{\mathbf{A}}$ is invertible. It is said to be of \textbf{type II} \cite{boscarino:2024, kennedy:2003} if the matrix $\tilde{\mathbf{A}}$ can be written in the form
    $$\tilde{\mathbf{A}} = 
       \begin{pmatrix}
	   0 & 0 \\
	   \tilde{\mathbf{a}} & \tilde{\boldsymbol{\mathcal{A}}}
        \end{pmatrix},
    $$
    with $\tilde{\mathbf{a}} = \left(\tilde{a}_{21}, \dots, \tilde{a}_{s1}\right)^{\top} \in \mathbb{R}^{s-1}$ and the matrix $\tilde{\boldsymbol{\mathcal{A}}} \in \mathbb{R}^{(s-1) \times (s-1)}$ is invertible. In the special case $\tilde{\mathbf{a}} = 0$, $b_1 = 0$, the method is said of \textbf{type ARS} (Ascher, Ruuth and Spiteri) \cite{ascher:1997} and the implicit method is reducible to a method using $s-1$ stages.
\end{definition}
We assume that the implicit scheme is a Diagonally Implicit Runge-Kutta (DIRK) method, namely $\tilde{a}_{lm} = 0$ for $l > m$. Following the Butcher tableaux introduced above, for a time dependent problem
\begin{equation}
    \frac{d\mathbf{y}}{dt} = \mathbf{f}_{E}\left(\mathbf{y}, t\right) + \mathbf{f}_{I}\left(\mathbf{y}, t\right),
\end{equation}
the generic $l$-stage of an IMEX-RK method can be defined as
\begin{eqnarray}\label{eq:stage_imex}
    \mathbf{v}^{(n,l)} &=& \mathbf{v}^{n} + \Delta t\sum_{m=1}^{l - 1}a_{lm}\mathbf{f}_{E}\left(\mathbf{v}^{(n,m)}, t^{n} + c_{m}\Delta t\right) + \Delta t\sum_{m=1}^{l}\tilde{a}_{lm}\mathbf{f}_{I}\left(\mathbf{v}^{(n,m)}, t^{n} + \tilde{c}_{m}\Delta t\right),
\end{eqnarray}
where $l = 1, \dots, s$, $\Delta t$ is the time discretization step, $\mathbf{v}^{n} \approx \mathbf{y}\left(t^{n}\right)$, $\mathbf{f}_{E}$ is the term treated explicitly, and $\mathbf{f}_{I}$ is the term treated implicitly. After computation of the intermediate stages, the updated solution is computed as follows:
\begin{equation}\label{eq:update_imex}
    \mathbf{v}^{n+1} = \mathbf{v}^{n} + \Delta t\sum_{l = 1}^{s}b_{l}\mathbf{f}_{E}\left(\mathbf{v}^{(n,l)}, t^{n} + c_{l}\Delta t\right) + \Delta t\sum_{l = 1}^{s}\tilde{b}_{l}\mathbf{f}_{I}\left(\mathbf{v}^{(n,l)}, t^{n} + \tilde{c}_{l}\Delta t\right).
\end{equation}
The formulation \eqref{eq:stage_imex}-\eqref{eq:update_imex} is valid for an IMEX scheme of arbitrary order. We recall that implicit methods of order higher than one cannot be unconditionally total variation diminishing (TVD) for hyperbolic problems \cite{gottlieb:2001}. This also holds for IMEX methods \cite{boscheri:2020, dimarco:2018}. In this work, as done, e.g., in \cite{boscheri:2020}, we do not focus on this limit imposed by high order schemes and we consider therefore numerical methods which, in principle, may not guarantee $L^{\infty}$-stability. Notice also that the existence of the Hilbert expansion \eqref{eq:discrete_expansion} can be justified only for smooth functions \cite{kucera:2022}. The development of a numerical treatment to avoid this issue goes beyond the scope of the present work and will be carried out as future development. For our analysis, we assume
\begin{equation}\label{eq:imex_compatibility}
    \sum_{m=1}^{s}a_{lm} = c_{l} \qquad \qquad \sum_{m=1}^{s}\tilde{a}_{lm} = \tilde{c}_{l}.
\end{equation}
Relation \eqref{eq:imex_compatibility} is an usual assumption for Runge-Kutta schemes \cite{boscarino:2016, kennedy:2003}, which simplifies the order conditions and, moreover, guarantees that a method of at least first order is employed at each stage. Notice that, for IMEX-RK methods of type I, $c_{1} \neq \tilde{c}_{1}$ because of \eqref{eq:imex_compatibility}. For the following analyses, we consider methods of type I for which $c_{l} = \tilde{c}_{l}$ for $l > 1$ and methods of type II with $\mathbf{c} = \tilde{\mathbf{c}}$. The assumption $\mathbf{c} = \tilde{\mathbf{c}}$ also allows to simplify the order conditions and has been employed, e.g., in \cite{ascher:1997, boscheri:2020}. A generic stage of the Euler equations reads as follows:
\begin{eqnarray}\label{eq:stage_euler}
    \rho^{(n,l)} &=& \rho^{n} - \sum_{m=1}^{l-1}a_{lm}\Delta t\dive\left(\rho^{(n,m)}\mathbf{u}^{(n,m)}\right) \nonumber \\
    \rho^{(n,l)}\mathbf{u}^{(n,l)} + \frac{1}{M^{2}}\tilde{a}_{ll}\Delta t\grad p^{(n,l)} &=& \rho^{n}\mathbf{u}^{n} - \sum_{m=1}^{l-1}a_{lm}\Delta t\dive\left(\rho^{(n,m)}\mathbf{u}^{(n,m)} \otimes \mathbf{u}^{(n,m)}\right) \nonumber \\
    &-& \frac{1}{M^{2}}\sum_{m=1}^{l-1}\tilde{a}_{lm}\Delta t\grad p^{(n,m)} \\
    \rho^{(n,l)}E^{(n,l)} + \tilde{a}_{ll}\Delta t\dive\left(h^{(n,l)}\rho^{(n,l)}\mathbf{u}^{(n,l)}\right) &=& \rho^{n}E^{n} - \sum_{m=1}^{l-1}\tilde{a}_{lm}\Delta t\dive\left(h^{(n,m)}\rho^{(n,m)}\mathbf{u}^{(n,m)}\right) \nonumber \\
    &-& \sum_{m=1}^{l-1}a_{lm}\Delta t M^{2}\dive\left(k^{(n,m)}\rho^{(n,m)}\mathbf{u}^{(n,m)}\right). \nonumber
\end{eqnarray}
We analyze now the behaviour of the time semi-discretization as $M \to 0$, so as to verify that it provides a consistent semi-discretization for the two limit models identified in Section \ref{ssec:ap_single_length} and Section \ref{ssec:ap_two_scale_length}, respectively.

\subsection{Asymptotic analysis in the single length scale case}
\label{ssec:ap_num_single_length}

In this Section, we consider the limit model \eqref{eq:euler_adim_ap}-\eqref{eq:euler_adim_ap_incomp}. Following, e.g., \cite{kucera:2022}, we make the assumption that, at each stage, the discrete quantities admit a formal expansion analogous to the continuous case.

\begin{assumption}\label{ass:discrete_expansion}
    The physical variables $\rho, \mathbf{u},$ and $p$ admit at each stage a formal Hilbert expansion of the form (written, e.g., for $\rho^{n}$)
    \begin{equation}\label{eq:discrete_expansion}
        \rho^{n}\left(\mathbf{x}\right) = \bar{\rho}^{n}\left(\mathbf{x}\right) + M\rho^{',n}\left(\mathbf{x}\right) + M^{2}\rho^{'',n}\left(\mathbf{x}\right) + \mathcal{O}(M^{3}).
    \end{equation}
\end{assumption}

\noindent
We also make the following assumption:

\begin{assumption}\label{ass:gradp0}
    In the case of schemes of type II that are not of type ARS, the initial datum $p^{0}$ is \textit{well-prepared}, namely $\grad\bar{p}^{0} = \grad p^{',0} = \mathbf{0}$.
\end{assumption}

\noindent
Then, the following result holds:

\begin{theorem}
    Under Assumption \ref{ass:discrete_expansion} and Assumption \ref{ass:gradp0}, \eqref{eq:stage_euler} provides a consistent discretization of \eqref{eq:euler_adim_ap}-\eqref{eq:euler_adim_ap_incomp} in the limit $M\rightarrow 0.$
\end{theorem}

\begin{proof}
    We plug asymptotic expansions of the form \eqref{eq:discrete_expansion} into \eqref{eq:stage_euler}. The discrete limit system reads therefore as follows:
    \begin{eqnarray}\label{eq:stage_euler_ap_single_length}
        \bar{\rho}^{(n,l)} &=& \bar{\rho}^{n} - \sum_{m=1}^{l-1}a_{lm}\Delta t\dive\left(\bar{\rho}^{(n,m)}\bar{\mathbf{u}}^{(n,m)}\right) \nonumber \\
        \tilde{a}_{ll}\grad\bar{p}^{(n,l)} &=& -\sum_{m=1}^{l-1}\tilde{a}_{lm}\grad\bar{p}^{(n,m)} \nonumber \\
        \tilde{a}_{ll}\grad p^{',(n,l)} &=& -\sum_{m=1}^{l-1}\tilde{a}_{lm}\grad p^{',(n,m)} \\
        \bar{\rho}^{(n,l)}\bar{\mathbf{u}}^{(n,l)} + \tilde{a}_{ll}\Delta t\grad p^{'',(n,l)} &=& \bar{\rho}^{n}\bar{\mathbf{u}}^{n} - \sum_{m=1}^{l-1}a_{lm}\Delta t\dive\left(\bar{\rho}^{(n,m)}\bar{\mathbf{u}}^{(n,m)} \otimes \bar{\mathbf{u}}^{(n,m)}\right) \nonumber \\
        &-& 
        \sum_{m=1}^{l-1}\tilde{a}_{lm}\Delta t\grad p^{'',(n,m)} \nonumber \\
        \bar{\rho}^{(n,l)}e\left(\bar{\rho}^{(n,l)}, \bar{p}^{(n,l)}\right) + \tilde{a}_{ll}\Delta t\dive\left(h\left(\bar{\rho}^{(n,l)}, \bar{p}^{(n,l)}\right)\bar{\rho}^{(n,l)}\bar{\mathbf{u}}^{(n,l)}\right) &=& \bar{\rho}^{n}e\left(\bar{\rho}^{n}, \bar{p}^{n}\right) - \sum_{m=1}^{l-1}\tilde{a}_{lm}\Delta t\dive\left(h\left(\bar{\rho}^{n}, \bar{p}^{n}\right)\bar{\rho}^{(n,m)}\bar{\mathbf{u}}^{(n,m)}\right). \nonumber
    \end{eqnarray} 
    First, we focus on the leading order terms of the momentum equation. For schemes of type I, since $\tilde{a}_{11} \neq 0$, we obtain $\grad\bar{p}^{(n,1)} = \mathbf{0}$ and therefore $\grad\bar{p}^{(n,l)} = \mathbf{0}$ for $l \ge 1$. For schemes of type ARS, since $\tilde{a}_{l1} = 0$, we obtain $\grad\bar{p}^{(n,l)} = \mathbf{0}$ for $l > 1$. For the other schemes, we need Assumption \eqref{ass:gradp0} to obtain a consistent discretization. Analogous considerations hold for the discretization of $\grad p^{'}$. The consistency of the remaining relations is a direct consequence of the consistency of the IMEX method. Nevertheless, we want to show that the last relation yields a consistent discretization for \eqref{eq:energy_limit_incomp}, so as to prove that \eqref{eq:stage_euler_ap_single_length} is a consistent discretization of \eqref{eq:euler_adim_ap_incomp}. After a few manipulations, taking into account that $\grad\bar{p}^{(n,l)} = \mathbf{0}$, we get
    \begin{eqnarray}\label{eq:stage_energy_incomp_single_length}
        &&\bar{\rho}^{(n,l)}h\left(\bar{\rho}^{(n,l)}, \bar{p}^{(n,l)}\right) - \bar{p}^{(n,l)} + \tilde{a}_{ll}\Delta t\left[\bar{\rho}^{(n,l)}h\left(\bar{\rho}^{(n,l)}, \bar{p}^{(n,l)}\right)\left(\dive\bar{\mathbf{u}}^{(n,l)}\right) + \mathbf{\bar{u}}^{(n,l)} \cdot \grad\bar{\rho}^{(n,l)}\frac{\partial\bar{\rho}^{(n,l)}h\left(\bar{\rho}^{(n,l)}, \bar{p}^{(n,l)}\right)}{\partial\bar{\rho}^{(n,l)}}\right] = \nonumber \\
        &&\bar{\rho}^{n}h\left(\bar{\rho}^{n}, \bar{p}^{n}\right) - \bar{p}^{n} - \sum_{m=1}^{l-1}\tilde{a}_{lm}\Delta t\left[\bar{\rho}^{(n,m)}h\left(\bar{\rho}^{n}, \bar{p}^{n}\right)\left(\dive\bar{\mathbf{u}}^{(n,m)}\right) + \mathbf{\bar{u}^{(n,m)}} \cdot \grad\bar{\rho}^{(n,m)}\frac{\partial\bar{\rho}^{(n,m)}h\left(\bar{\rho}^{n}, \bar{p}^{n}\right)}{\partial\bar{\rho}^{(n,m)}}\right].
    \end{eqnarray}
    From now, for the sake of simplicity in the notation, we denote $h\left(\bar{\rho}^{(n,l)}, \bar{p}^{(n,l)}\right)$ by $\bar{h}^{(n,l)}$ and $h\left(\bar{\rho}^{n}, \bar{p}^{n}\right)$ by $\bar{h}^{n}$. The error obtained applying \eqref{eq:stage_energy_incomp_single_length} to the exact solution reads therefore as follows:
    \newpage
    \begin{eqnarray}\label{eq:lte_energy_incomp_single_length}
        \hat{\tau}^{(n,l)} &=& \bar{\rho}\left(\mathbf{x}, t^{n} + c_{l}\Delta t\right)\bar{h}\left(\mathbf{x}, t^{n} + c_{l}\Delta t\right) - \bar{\rho}\left(\mathbf{x}, t^{n}\right)\bar{h}\left(\mathbf{x}, t^{n}\right) - \left[\bar{p}\left(\mathbf{x}, t^{n} + c_{l}\Delta t\right) - \bar{p}\left(\mathbf{x}, t^{n}\right)\right] \nonumber \\
        &+& \tilde{a}_{ll}\Delta t\left[\bar{\rho}\left(\mathbf{x}, t^{n} + c_{l}\Delta t\right)\bar{h}\left(\mathbf{x}, t^{n} + c_{l}\Delta t\right)\left(\dive\bar{\mathbf{u}}\left(\mathbf{x}, t^{n} + c_{l}\Delta t\right)\right)\right] \nonumber \\
        &+& \tilde{a}_{ll}\Delta t\left[\bar{\mathbf{u}}\left(\mathbf{x}, t^{n} + c_{l}\Delta t\right) \cdot \grad\bar{\rho}\left(\mathbf{x}, t^{n} + c_{l}\Delta t\right)\frac{\partial\bar{\rho}\bar{h}}{\partial\bar{\rho}}\left(\mathbf{x}, t^{n} + c_{l}\Delta t\right)\right] \\
        &+& \sum_{m=1}^{l-1}\tilde{a}_{lm}\Delta t\left[\bar{\rho}\left(\mathbf{x}, t^{n} + c_{m}\Delta t\right)\bar{h}\left(\mathbf{x}, t^{n} + c_{m}\Delta t\right)\left(\dive\bar{\mathbf{u}}\left(\mathbf{x}, t^{n} + c_{m}\Delta t\right)\right)\right] \nonumber \\
        &+& \sum_{m=1}^{l-1}\tilde{a}_{lm}\Delta t\left[\bar{\mathbf{u}}\left(\mathbf{x}, t^{n} + c_{m}\Delta t\right) \cdot \grad\bar{\rho}\left(\mathbf{x}, t^{n} + c_{m}\Delta t\right)\frac{\partial\bar{\rho}\bar{h}}{\partial\bar{\rho}}\left(\mathbf{x}, t^{n} + c_{m}\Delta t\right)\right] \nonumber .
    \end{eqnarray}
    Since $\bar{\rho}\bar{h} = \bar{\rho}\bar{h}\left(\bar{\rho},\bar{p}\right)$, using a Taylor expansion, we get
    \begin{eqnarray}
        \bar{\rho}\left(\mathbf{x}, t^{n} + c_{l}\Delta t\right)\bar{h}\left(\mathbf{x}, t^{n} + c_{l}\Delta t\right) &=& \bar{\rho}\left(\mathbf{x}, t^{n}\right)\bar{h}\left(\mathbf{x}, t^{n}\right) + \frac{\partial\bar{\rho}\bar{h}}{\partial\bar{\rho}}\left(\mathbf{x}, t^{n}\right)\left[\bar{\rho}\left(\mathbf{x}, t^{n} + c_{l}\Delta t\right) - \bar{\rho}\left(\mathbf{x}, t^{n}\right)\right] \nonumber \\
        &+& \frac{\partial\bar{\rho}\bar{h}}{\partial\bar{p}}\left(\mathbf{x}, t^{n}\right)\left[\bar{p}\left(\mathbf{x}, t^{n} + c_{l}\Delta t\right) - \bar{p}\left(\mathbf{x}, t^{n}\right)\right] \\
        &+& o\left(\bar{\rho}\left(\mathbf{x}, t^{n} + c_{l}\Delta t\right) - \bar{\rho}\left(\mathbf{x}, t^{n}\right)\right)+ o\left(\bar{p}\left(\mathbf{x}, t^{n} + c_{l}\Delta t\right) - \bar{p}\left(\mathbf{x}, t^{n}\right)\right). \nonumber
    \end{eqnarray}
    Employing now the discretization of the continuity equation in \eqref{eq:stage_euler_ap_single_length}, we obtain for $l > 1$
    \begin{eqnarray}
        \bar{\rho}\left(\mathbf{x}, t^{n} + c_{l}\Delta t\right)\bar{h}\left(\mathbf{x}, t^{n} + c_{l}\Delta t\right) &=& \bar{\rho}\left(\mathbf{x}, t^{n}\right)\bar{h}\left(\mathbf{x}, t^{n}\right) - \sum_{m=1}^{l-1}a_{lm}\Delta t \frac{\partial\bar{\rho}\bar{h}}{\partial\bar{\rho}}\left(\mathbf{x}, t^{n}\right)\dive\left(\bar{\rho}\left(\mathbf{x}, t^{n} + c_{m}\Delta t\right)\bar{\mathbf{u}}\left(\mathbf{x}, t^{n} + c_{m}\Delta t\right)\right) \nonumber \\
        &+& \frac{\partial\bar{\rho}\bar{h}}{\partial\bar{p}}\left(\mathbf{x}, t^{n}\right)\left[\bar{p}\left(\mathbf{x}, t^{n} + c_{l}\Delta t\right) - \bar{p}\left(\mathbf{x}, t^{n}\right)\right] \\
        &+& o\left(\bar{\rho}\left(\mathbf{x}, t^{n} + c_{l}\Delta t\right) - \bar{\rho}\left(\mathbf{x}, t^{n}\right)\right) + o\left(\bar{p}\left(\mathbf{x}, t^{n} + c_{l}\Delta t\right) - \bar{p}\left(\mathbf{x}, t^{n}\right)\right), \nonumber
    \end{eqnarray}
    or, equivalently,
    \begin{eqnarray}\label{eq:rhoh_taylor_expansion}
        \bar{\rho}\left(\mathbf{x}, t^{n} + c_{l}\Delta t\right)\bar{h}\left(\mathbf{x}, t^{n} + c_{l}\Delta t\right) &=& \bar{\rho}\left(\mathbf{x}, t^{n}\right)\bar{h}\left(\mathbf{x}, t^{n}\right) - \sum_{m=1}^{l-1}a_{lm}\Delta t \frac{\partial\bar{\rho}\bar{h}}{\partial\bar{\rho}}\left(\mathbf{x}, t^{n}\right)\bar{\rho}\left(\mathbf{x}, t^{n} + c_{m}\Delta t\right)\dive\bar{\mathbf{u}}\left(\mathbf{x}, t^{n} + c_{m}\Delta t\right) \nonumber \\
        &-& \sum_{m=1}^{l-1}a_{lm}\Delta t \frac{\partial\bar{\rho}\bar{h}}{\partial\bar{\rho}}\left(\mathbf{x}, t^{n}\right)\bar{\mathbf{u}}\left(\mathbf{x}, t^{n} + c_{m}\Delta t\right) \cdot \grad\bar{\rho}\left(\mathbf{x}, t^{n} + c_{m}\Delta t\right) \nonumber \\
        &+& \frac{\partial\bar{\rho}\bar{h}}{\partial\bar{p}}\left(\mathbf{x}, t^{n}\right)\left[\bar{p}\left(\mathbf{x}, t^{n} + c_{l}\Delta t\right) - \bar{p}\left(\mathbf{x}, t^{n}\right)\right] \\
        &+& o\left(\bar{\rho}\left(\mathbf{x}, t^{n} + c_{l}\Delta t\right) - \bar{\rho}\left(\mathbf{x}, t^{n}\right)\right) + o\left(\bar{p}\left(\mathbf{x}, t^{n} + c_{l}\Delta t\right) - \bar{p}\left(\mathbf{x}, t^{n}\right)\right). \nonumber
    \end{eqnarray}
    Using again a Taylor expansion, we get 
    \begin{eqnarray}\label{eq:rho_u_p_taylor_expansion}
        \bar{\rho}\left(\mathbf{x}, t^{n} + c_{l}\Delta t\right) = \bar{\rho}\left(\mathbf{x}, t^{n}\right) + c_{l}\Delta t\frac{\partial\bar{\rho}}{\partial t}\left(\mathbf{x}, t^{n}\right) + \mathcal{O}(\Delta t^2) \nonumber \\
        \bar{\mathbf{u}}\left(\mathbf{x}, t^{n} + c_{l}\Delta t\right) = \bar{\mathbf{u}}\left(\mathbf{x}, t^{n}\right) + c_{l}\Delta t\frac{\partial\bar{\mathbf{u}}}{\partial t}\left(\mathbf{x}, t^{n}\right) + \mathcal{O}(\Delta t^2) \\
        \bar{p}\left(\mathbf{x}, t^{n} + c_{l}\Delta t\right) = \bar{p}\left(\mathbf{x}, t^{n}\right) + c_{l}\Delta t\frac{\partial\bar{p}}{\partial t}\left(\mathbf{x}, t^{n}\right) + \mathcal{O}(\Delta t^2). \nonumber
    \end{eqnarray}
    Substituting \eqref{eq:rhoh_taylor_expansion} and \eqref{eq:rho_u_p_taylor_expansion} into \eqref{eq:lte_energy_incomp_single_length}, we obtain
    \newpage
    \begin{eqnarray}
        \hat{\tau}^{(n,l)} &=& -\sum_{m=1}^{l-1}a_{lm}\Delta t \frac{\partial\bar{\rho}\bar{h}}{\partial\bar{\rho}}\left(\mathbf{x}, t^{n}\right)\bar{\rho}\left(\mathbf{x}, t^{n}\right)\dive\bar{\mathbf{u}}\left(\mathbf{x}, t^{n}\right) - \sum_{m=1}^{l-1}a_{lm}\Delta t \frac{\partial\bar{\rho}\bar{h}}{\partial\bar{\rho}}\left(\mathbf{x}, t^{n}\right)\bar{\mathbf{u}}\left(\mathbf{x}, t^{n}\right) \cdot \grad\bar{\rho}\left(\mathbf{x}, t^{n}\right) \nonumber \\
        &+& \left(\frac{\partial\bar{\rho}\bar{h}}{\partial\bar{p}}\left(\mathbf{x}, t^{n}\right) - 1\right)c_{l}\Delta t\frac{\partial\bar{p}}{\partial t}\left(\mathbf{x}, t^{n}\right) + \sum_{m=1}^{l}\tilde{a}_{lm}\Delta t\left[\bar{\rho}\left(\mathbf{x}, t^{n}\right)\bar{h}\left(\mathbf{x}, t^{n}\right)\dive\bar{\mathbf{u}}\left(\mathbf{x}, t^{n}\right)\right] \\
        &+& \sum_{m=1}^{l}\tilde{a}_{lm}\Delta t\left[\bar{\mathbf{u}}\left(\mathbf{x}, t^{n}\right) \cdot \grad\bar{\rho}\left(\mathbf{x}, t^{n}\right)\frac{\partial\bar{\rho}\bar{h}}{\partial\bar{\rho}}\left(\mathbf{x}, t^{n}\right)\right] + \mathcal{O}\left(\Delta t^2\right) \nonumber.
    \end{eqnarray}
    Since $\sum\limits_{m=1}^{l}\tilde{a}_{lm} =$ $\sum\limits_{m=1}^{l-1}a_{lm} =$ $c_{l}$ \eqref{eq:imex_compatibility} and $\bar{\rho}\left(\mathbf{x}, t^{n}\right)\bar{h}\left(\mathbf{x}, t^{n}\right) - \frac{\partial\bar{\rho}\bar{h}}{\partial\bar{\rho}}\left(\mathbf{x}, t^{n}\right)\bar{\rho}\left(\mathbf{x}, t^{n}\right) = -\bar{\rho}^{2}\left(\mathbf{x}, t^{n}\right)\frac{\partial\bar{h}}{\partial\bar{\rho}}\left(\mathbf{x}, t^{n}\right)$, we obtain
    \begin{equation}
        \hat{\tau}^{(n,l)} = -c_{l}\Delta t\bar{\rho}^{2}\left(\mathbf{x}, t^{n}\right) \frac{\partial\bar{h}}{\partial\bar{\rho}}\left(\mathbf{x}, t^{n}\right)\dive\bar{\mathbf{u}}\left(\mathbf{x}, t^{n}\right) + \left(\frac{\partial\bar{\rho}\bar{h}}{\partial\bar{p}}\left(\mathbf{x}, t^{n}\right) - 1\right)c_{l}\Delta t\frac{\partial\bar{p}}{\partial t}\left(\mathbf{x}, t^{n}\right) + \mathcal{O}\left(\Delta t^2\right) = \mathcal{O}\left(\Delta t^2\right),
    \end{equation}
    thanks to \eqref{eq:energy_limit_incomp}. 
    Finally, the update stage for the energy equation reads as follows:
    \begin{equation}\label{eq:energy_update}
        \bar{\rho}^{n+1}e\left(\bar{\rho}^{n+1},\bar{p}^{n+1}\right) = \bar{\rho}^{n}e\left(\bar{\rho}^{n},\bar{p}^{n}\right) - \sum_{l=1}^{s}\tilde{b}_{l}\dive\left(\bar{\rho}^{(n,m)}h\left(\bar{\rho}^{(n,m)},\bar{p}^{(n,m)}\right)\bar{\mathbf{u}}^{(n,m)}\right)
    \end{equation}
    The error obtained applying \eqref{eq:energy_update} to the exact solution is
    \begin{eqnarray}
        \hat{\tau}^{n+1} &=& \bar{\rho}\left(\mathbf{x}, t^{n} + \Delta t\right)e\left(\bar{\rho}\left(\mathbf{x}, t^{n} + \Delta t\right),\bar{p}\left(\mathbf{x}, t^{n} + \Delta t\right)\right) -\bar{\rho}\left(\mathbf{x}, t^{n}\right)e\left(\bar{\rho}\left(\mathbf{x}, t^{n}\right),\bar{p}\left(\mathbf{x}, t^{n}\right)\right) \nonumber \\
        &-& \sum_{l=1}^{s}\Delta t\tilde{b}_{l}\dive\left(\bar{\rho}\left(\mathbf{x}, t^{n} + \tilde{c}_{l}\Delta t\right)h\left(\bar{\rho}\left(\mathbf{x}, t^{n} + \tilde{c}_{l}\Delta t\right),\bar{p}\left(\mathbf{x}, t^{n} + \tilde{c}_{l}\Delta t\right)\right)\bar{\mathbf{u}}\left(\mathbf{x}, t^{n} + \tilde{c}_{l}\Delta t\right)\right).
    \end{eqnarray}
    Thanks to a Taylor expansion, we get
    \begin{equation}
        \hat{\tau}^{n+1} = \Delta t\left[\frac{\partial\bar{\rho}e\left(\bar{\rho}, \bar{p}\right)}{\partial t}\left(\mathbf{x}, t^{n}\right) + \sum_{l=1}^{s}\tilde{b}_{l}\dive\left(\bar{\rho}\left(\mathbf{x}, t^{n}\right)h\left(\bar{\rho}\left(\mathbf{x}, t^{n}\right),\bar{p}\left(\mathbf{x}, t^{n} \right)\right)\bar{\mathbf{u}}\left(\mathbf{x}, t^{n}\right)\right)\right] + \mathcal{O}\left(\Delta t^2\right).
    \end{equation}
    Since $\sum\limits_{l=1}^{s} \tilde{b}_{l} = 1$ and thanks to \eqref{eq:energy_limit_incomp}, we obtain
    \begin{equation}
        \hat{\tau}^{n+1} = \Delta t\left[\frac{\partial\bar{\rho}e\left(\bar{\rho}, \bar{p}\right)}{\partial t}\left(\mathbf{x}, t^{n}\right) + \dive\left(\bar{\rho}\left(\mathbf{x}, t^{n}\right)h\left(\bar{\rho}\left(\mathbf{x}, t^{n}\right),\bar{p}\left(\mathbf{x}, t^{n} \right)\right)\bar{\mathbf{u}}\left(\mathbf{x}, t^{n}\right)\right)\right] + \mathcal{O}\left(\Delta t^2\right) = \mathcal{O}\left(\Delta t^2\right).
    \end{equation}
    The consistency for the remaining relations can be shown in an analogous manner.
\end{proof}

Since we are considering an implicit coupling between the momentum and the energy balance, the stability condition of the numerical method does not depend on $M$ or on the acoustic speed of sound (see, e.g., \cite{casulli:1984, dumbser:2016b, tavelli:2017}), meaning that \eqref{eq:stage_euler} provides an AP scheme for \eqref{eq:euler_adim_ap}-\eqref{eq:euler_adim_ap_incomp}.   Only a mild CFL-type restriction based on the flow velocity is necessary for these schemes \cite{dumbser:2016b}. Schemes of type I and of type ARS are also strongly asymptotic-preserving, namely they are asymptotic-preserving for general initial data. On the contrary, schemes of type II which are not of type ARS turn out to be weakly asymptotic-preserving. Indeed, they require a well-prepared initial datum for the pressure (see Assumption \ref{ass:gradp0}). However, when the limit model reduces to the incompressible Euler equations thanks to suitable boundary conditions, we do not need a divergence-free initial velocity field to recover the incompressible limit, as we will also verify numerically in Section \ref{ssec:density_layering}. This is not valid for all the AP methods presented in the literature, see, e.g., \cite{noelle:2014}, which has been subsequently corrected in \cite{bispen:2017}, or \cite{thomann:2019}.

The AP property guarantees the consistency of the discretization as $M \to 0$, but it does not imply that a scheme preserves its order of accuracy as $M \to 0$. In this latter case, the scheme is said to be asymptotically-accurate (AA). Since the seminal paper \cite{pareschi:2005}, it is quite established that \texttt{L}-stability is necessary to guarantee asymptotic accuracy. A Runge-Kutta scheme is said to be $\texttt{L}$-stable \cite{wanner:1996} if it is \texttt{A}-stable and $R(z) \to 0$ as $z \to \infty$, where $R(z)$ is the stability function. Following the result in \cite{wanner:1996}, a \texttt{L}-stable scheme results from the combination of a \texttt{A}-stable scheme with a stiffly-accurate (SA) scheme, i.e. a scheme for which the update stage is identical to the last internal stage. However, for methods of type II, this combination does not necessarily lead to a \texttt{L}-stable scheme, because the matrix $\tilde{\mathbf{A}}$ is not invertible \cite{boscarino:2009}. In the case of methods of type II which are stiffly-accurate, a supplementary condition is required to obtain the \texttt{L}-stability, i.e. \cite{boscarino:2009}
\begin{equation}\label{eq:Lstab_type_II}
    R_{\infty} = \sum_{m=2}^{s}\hat{w}_{sm}\tilde{a}_{m1} = 0,
\end{equation}
where $\hat{w}_{lm}$ denotes the elements of the inverse of $\tilde{\boldsymbol{\mathcal{A}}}$. Hence, SA schemes of type ARS are also \texttt{L}-stable.

\subsection{Asymptotic analysis for two length scales}
\label{ssec:ap_num_two_scale_length}

In this Section, we consider the limit model \eqref{eq:euler_adim_ap_two_scale}. We replace Assumption \ref{ass:discrete_expansion} with the following one:

\begin{assumption}\label{ass:discrete_expansion_two_scale}
    The physical variables $\rho, \mathbf{u}$, and $p$ admit at each stage a formal Hilbert expansion of the form (written, e.g., for $\rho^{n}$)
    \begin{equation}\label{eq:discrete_expansion_two_scale}
        \rho^{n}\left(\mathbf{x}, \boldsymbol{\xi}\right) = \bar{\rho}^{n}\left(\mathbf{x}, \boldsymbol{\xi}\right) + M\rho^{',n}\left(\mathbf{x}, \boldsymbol{\xi}\right) + M^{2}\rho^{'',n}\left(\mathbf{x}, \boldsymbol{\xi}\right) + \mathcal{O}(M^{3}),
    \end{equation}
    with $\boldsymbol{\xi} = M\mathbf{x}$.
\end{assumption}

\noindent
Moreover, we replace Assumption \ref{ass:gradp0} with the 

\begin{assumption}\label{ass:gradp0_two_scale}
    In the case of IMEX-RK schemes of type II that are not of type ARS, the initial datum $p^{0}$ is \textit{well-prepared},  namely $\gradx\bar{p}^{0} = \gradx\bar{p}^{'} = \gradxi p^{'} = \mathbf{0}$.
\end{assumption}

\noindent
Then, the following result holds:

\begin{theorem}
    Under Assumption \ref{ass:discrete_expansion_two_scale} and Assumption \ref{ass:gradp0_two_scale}, \eqref{eq:stage_euler} provides an AP scheme for \eqref{eq:euler_adim_ap_two_scale}.
\end{theorem}

\begin{proof}
    As pointed out for the continuous model, the leading order term relations do not change when also introducing the acoustic length scale variable $\boldsymbol{\xi} = M\mathbf{x}$. We plug asymptotic expansion of the form \eqref{eq:discrete_expansion_two_scale} into the semi-discretized momentum equation, so as to obtain for the first order term
    \begin{equation}
        \tilde{a}_{ll}\gradx p^{',(n,l)} + \tilde{a}_{ll}\gradxi\bar{p}^{(n,l)} = -\sum_{m=1}^{l-1}\tilde{a}_{lm}\gradx p^{',(n,m)} -\sum_{m=1}^{l-1}\tilde{a}_{lm}\gradxi\bar{p}^{(n,m)}.
    \end{equation}
    Analogous considerations to those reported in Section \ref{ssec:ap_num_single_length} hold for the above relation. More specifically, for schemes of type I, since $\tilde{a}_{11} \neq 0$, we obtain
    $$\gradx p^{',(n,l)} +\gradxi\bar{p}^{(n,l)} = \mathbf{0} \qquad \text{for } l \ge 1.$$
    For what concerns schemes of type ARS, since $\tilde{a}_{l1} = 0, l=1 \dots s,$ we get
    $$\gradx p^{',(n,l)} +\gradxi\bar{p}^{(n,l)} = \mathbf{0} \qquad \text{for } l > 1.$$
    Assumption \ref{ass:gradp0_two_scale} is instead necessary to obtain a consistent discretization of \eqref{eq:momentum_limit_first_order_two_scale} for the other schemes. For what concerns the second order term, we get
    \begin{eqnarray}\label{eq:momemtum_two_scale_IMEX}
        \bar{\rho}^{(n,l)}\bar{\mathbf{u}}^{(n,l)} + \tilde{a}_{ll}\Delta t\left(\gradx p^{'',(n,l)} + \gradxi p^{',(n,l)}\right) &=& \bar{\rho}^{n}\bar{\mathbf{u}}^{n} - \sum_{m=1}^{l-1}a_{lm}\Delta t\divex\left(\bar{\rho}^{(n,m)}\bar{\mathbf{u}}^{(n,m)} \otimes \bar{\mathbf{u}}^{(n,m)}\right) \\
        &-& \sum_{m=1}^{l-1}\tilde{a}_{lm}\Delta t\left(\gradx p^{'',(n,m)} + \gradxi p^{',(n,m)}\right). \nonumber
    \end{eqnarray}
    One can easily verify that \eqref{eq:momemtum_two_scale_IMEX} is a consistent discretization of \eqref{eq:momentum_limit_second_order_two_scale}. Indeed, the error obtained applying \eqref{eq:momemtum_two_scale_IMEX} to the exact solution reads as follows:
    \newpage
    \begin{eqnarray}
        \hat{\tau}_{\xi}^{(n,l)} &=& \bar{\rho}\left(\mathbf{x}, t^{n} + c_{l}\Delta t\right)\bar{\mathbf{u}}\left(\mathbf{x}, t^{n} + c_{l}\Delta t\right) - \bar{\rho}\left(\mathbf{x}, t^{n}\right)\bar{\mathbf{u}}\left(\mathbf{x}, \Delta t\right) \nonumber \\
        &+& \sum_{m=1}^{l-1}a_{lm}\Delta t\divex\left(\bar{\rho}\left(\mathbf{x}, t^{n} + c_{m}\Delta t\right)\bar{\mathbf{u}}\left(\mathbf{x}, t^{n} + c_{m}\Delta t\right) \otimes \bar{\mathbf{u}}\left(\mathbf{x}, t^{n} + c_{m}\Delta t\right)\right) \nonumber \\
        &+& \sum_{m=1}^{l-1}\tilde{a}_{lm}\Delta t\left(\gradx p^{''}\left(\mathbf{x}, t^{n} + c_{m}\Delta t\right) + \gradxi p^{'}\left(\mathbf{x}, t^{n} + c_{m}\Delta t\right)\right)
    \end{eqnarray}
    Thanks to a Taylor expansion, we get
    \begin{eqnarray}
        \hat{\tau}_{\xi}^{(n,l)} &=& \Delta t\left(c_{l}\frac{\partial\bar{\rho}\bar{\mathbf{u}}}{\partial t}\left(\mathbf{x}, t^{n}\right) + \sum_{m=1}^{l-1}a_{lm}\divex\left(\bar{\rho}\left(\mathbf{x}, t^{n}\right)\bar{\mathbf{u}}\left(\mathbf{x}, t^{n}\right) \otimes \bar{\mathbf{u}}\left(\mathbf{x}, t^{n}\right)\right) + \sum_{m=1}^{l-1}\tilde{a}_{lm}\left(\gradx p^{''}\left(\mathbf{x}, t^{n}\right) + \gradxi p^{'}\left(\mathbf{x}, t^{n}\right)\right)\right) \nonumber \\
        &+& \mathcal{O}(\Delta t^{2}).
    \end{eqnarray}
    Since $\sum\limits_{m=1}^{l}\tilde{a}_{lm} = \sum\limits_{m=1}^{l-1}a_{lm} = c_{l}$ \eqref{eq:imex_compatibility}, we obtain
    \begin{eqnarray}
        \hat{\tau}_{\xi}^{(n,l)} &=& c_{l}\Delta t\left(\frac{\partial\bar{\rho}\bar{\mathbf{u}}}{\partial t}\left(\mathbf{x}, t^{n}\right) + \divex\left(\bar{\rho}\left(\mathbf{x}, t^{n}\right)\bar{\mathbf{u}}\left(\mathbf{x}, t^{n}\right) \otimes \bar{\mathbf{u}}\left(\mathbf{x}, t^{n}\right)\right) + \left(\gradx p^{''}\left(\mathbf{x}, t^{n}\right) + \gradxi p^{'}\left(\mathbf{x}, t^{n}\right)\right)\right) \nonumber \\
        &+& \mathcal{O}(\Delta t^{2}) = \mathcal{O}(\Delta t^{2}).
    \end{eqnarray}
    Finally, for the first order term in the energy equation, we obtain
    \begin{eqnarray}\label{eq:energy_two_scale_IMEX}
        &&\rho^{',(n,l)}e^{',(n,l)} + \tilde{a}_{ll}\Delta t\divex\left(h^{',(n,l)}\rho^{',(n,l)}\mathbf{u}^{',(n,l)}\right) + \tilde{a}_{ll}\Delta t\divexi\left(\bar{h}^{(n,l)}\bar{\rho}^{(n,l)}\bar{\mathbf{u}}^{(n,l)}\right) = \nonumber \\
        &&\bar{\rho}^{n}\bar{e}^{n} - \sum_{m=1}^{l-1}\tilde{a}_{lm}\Delta t\divex\left(h^{',(n,m)}\rho^{',(n,m)}\mathbf{u}^{',(n,m)}\right) - \sum_{m=1}^{l-1}\tilde{a}_{lm}\Delta t\divexi\left(\bar{h}^{(n,m)}\bar{\rho}^{(n,m)}\bar{\mathbf{u}}^{(n,m)}\right).
    \end{eqnarray}
    Analogously, one can show that \eqref{eq:energy_two_scale_IMEX} is a consistent discretization of \eqref{eq:energy_limit_first_order_two_scale}. Similar computations show the consistency for the final update stage. Hence, \eqref{eq:stage_euler} provides an AP scheme for \eqref{eq:euler_adim_ap_two_scale}.
\end{proof}

\section{Spatial discretization}
\label{sec:space}

In this Section, we briefly outline the spatial discretization for \eqref{eq:stage_euler}, which is based on the Discontinuous Galerkin (DG) method \cite{giraldo:2020} as implemented in the \texttt{deal.II} library \cite{arndt:2023, bangerth:2007}. We use quadrilateral elements and the corresponding polynomial spaces $Q_{r}$ of degree $r$ \cite{quarteroni:2008}. More specifically, the shape functions correspond to the products of Lagrange polynomials for the support points of $\left(r + 1\right)$-order Gauss-Lobatto quadrature rule in each coordinate direction, where $r$ is the polynomial degree. However, the proposed approach can also be applied to tetrahedral meshes and $P$-spaces. We consider a decomposition of the domain $\Omega$ into a family of quadrilaterals $\mathcal{T}_{\mathcal{H}}$ and denote each element by $K$. We denote by $\mathcal{E}$ the set of all the element faces, so that $\mathcal{E} = \mathcal{E}^{I} \cup \mathcal{E}^{B}$, with $\mathcal{E}^{I}$ and $\mathcal{E}^{B}$ denoting the subset of interior and boundary faces, respectively. A face $\Gamma \in \mathcal{E}_{I}$ shares two elements, $K^{+}$ with outward unit normal $\mathbf{n}^{+}$ and $K^{-}$ with outward unit normal $\mathbf{n}^{-}$. Finally, we denote by $\mathbf{n}$ the outward unit normal for a face $\Gamma \in \mathcal{E}^{B}$. Hence, following, e.g., \cite{arnold:2002}, for a scalar function $\varphi$, we define the jump as
\begin{equation}
    \left[\left[\varphi\right]\right] = \varphi^{+}\mathbf{n}^{+} + \varphi^{-}\mathbf{n}^{-} \text{ if } \Gamma \in \mathcal{E}^{I} \qquad \left[\left[\varphi\right]\right] = \varphi\mathbf{n} \text{ if } \Gamma \in \mathcal{E}^{B},
\end{equation}
where we define the average as
\begin{equation}
    \left\{\left\{\varphi\right\}\right\} = \frac{1}{2}\left(\varphi^{+} + \varphi^{-}\right) \text{ if } \Gamma \in \mathcal{E}^{I} \qquad \left\{\left\{\varphi\right\}\right\} = \varphi \text{ if } \Gamma \in \mathcal{E}^{B}.
\end{equation}
Analogous definitions apply for a vector function $\boldsymbol{\varphi}$. More specifically, we define
\begin{alignat}{2}
    \left[\left[\boldsymbol{\varphi}\right]\right] &= \boldsymbol{\varphi}^{+} \cdot \mathbf{n}^{+} + \boldsymbol{\varphi}^{-} \cdot \mathbf{n}^{-} \text{ if } \Gamma \in \mathcal{E}^{I} \qquad &&\left[\left[\boldsymbol{\varphi}\right]\right] = \boldsymbol{\varphi} \cdot \mathbf{n} \text{ if } \Gamma \in \mathcal{E}^{B} \\
    \left\{\left\{\boldsymbol{\varphi}\right\}\right\} &= \frac{1}{2}\left(\boldsymbol{\varphi}^{+} + \boldsymbol{\varphi}^{-}\right) \text{ if } \Gamma \in \mathcal{E}^{I} \qquad &&\left\{\left\{\boldsymbol{\varphi}\right\}\right\} = \boldsymbol{\varphi} \text{ if } \Gamma \in \mathcal{E}^{B}.
\end{alignat}
For vector functions, it is also useful to define a tensor jump as follows:
\begin{equation}
    \left<\left<\boldsymbol{\varphi}\right>\right> = \boldsymbol{\varphi}^{+} \otimes \mathbf{n}^{+} + \boldsymbol{\varphi}^{-} \otimes \mathbf{n}^{-} \text{ if } \Gamma \in \mathcal{E}^{I} \qquad \left<\left<\boldsymbol{\varphi}\right>\right> = \boldsymbol{\varphi} \otimes \mathbf{n} \text{ if } \Gamma \in \mathcal{E}^{B}.
\end{equation}
Given these definitions, the weak formulation for the momentum equation at each stage \eqref{eq:stage_euler} reads as follows \cite{orlando:2022b, orlando:2023b}:
\begin{equation}
    \mathbf{A}^{(n,l)}\mathbf{U}^{(n,l)} + \mathbf{B}^{(n,l)}\mathbf{P}^{(n,l)} = \mathbf{F}^{(n,l)},
\end{equation}
with $\mathbf{U}^{(n,l)}$ denoting the vector of the degrees of freedom associated to the velocity field and $\mathbf{P}^{(n,l)}$ denoting the vector of the degrees of freedom associated to the pressure. Here we have set
\begin{eqnarray}
    A_{ij}^{(n,l)} &=& \sum_{K \in \mathcal{T}_{\mathcal{H}}} \int_{K} \rho^{(n,l)}\boldsymbol{\varphi}_{j} \cdot \boldsymbol{\varphi}_{i}d\Omega \\
    B_{ij}^{(n,l)} &=& \sum_{K \in \mathcal{T}_{\mathcal{H}}}\int_{K} -\tilde{a}_{ll}\frac{\Delta t}{M^{2}}\dive\boldsymbol{\varphi}_{i}\Psi_{j}d\Omega + \sum_{\Gamma \in \mathcal{E}}\int_{\Gamma}\tilde{a}_{ll}\frac{\Delta t}{M^{2}}\left\{\left\{\Psi_{j}\right\}\right\}\left[\left[\boldsymbol{\varphi}_{i}\right]\right]d\Sigma \\
    F_{i}^{(n,l)} &=& \sum_{K \in \mathcal{T}_{\mathcal{H}}}\int_{K}\rho^{n}\mathbf{u}^{n} \cdot \boldsymbol{\varphi}_{i}d\Omega + \sum_{m = 1}^{l-1}\sum_{K \in \mathcal{T}_{\mathcal{H}}}\int_{K}a_{lm}\Delta t\left(\rho^{(n,m)}\mathbf{u}^{(n,m)} \otimes \mathbf{u}^{(n,m)}\right) : \grad\boldsymbol{\varphi}_{i}d\Omega \nonumber \nonumber \\
    &+& \sum_{m=1}^{l-1}\sum_{K \in \mathcal{T}_{\mathcal{H}}}\int_{K}\tilde{a}_{lm}\frac{\Delta t}{M^{2}}p^{(n,m)}\dive\boldsymbol{\varphi}_{i}d\Omega \nonumber \\
    &-& \sum_{m=1}^{l-1}\sum_{\Gamma \in \mathcal{E}}\int_{\Gamma}a_{lm}\Delta t\left\{\left\{\rho^{(n,m)}\mathbf{u}^{(n,m)} \otimes \mathbf{u}^{(n,m)}\right\}\right\} : \left<\left<\boldsymbol{\varphi}_{i}\right>\right>d\Sigma \\
    &-& \sum_{m=1}^{l-1}\sum_{\Gamma \in \mathcal{E}}\int_{\Gamma}a_{lm}\Delta t\frac{\lambda^{(n,m)}}{2}\left<\left<\rho^{(n,m)}\mathbf{u}^{(n,m)}\right>\right> : \left<\left<\boldsymbol{\varphi}_{i}\right>\right>d\Sigma - \sum_{m=1}^{l-1}\sum_{\Gamma \in \mathcal{E}}\int_{\Gamma}\tilde{a}_{lm}\frac{\Delta t}{M^{2}}\left\{\left\{p^{(n,m)}\right\}\right\}\left[\left[\boldsymbol{\varphi}_{i}\right]\right]d\Sigma, \nonumber
\end{eqnarray}
with $\boldsymbol{\varphi}_{i}$ and $\Psi_{i}$ denoting the basis function of the space of polynomial functions employed to discretize the velocity and the pressure, respectively. Following the discussion in \cite{orlando:2023c, orlando:2022b}, one can notice that we employ a centered flux for the quantities defined implicitly and upwind-biased flux for the quantities computed explicitly. The choice of the upwind-biased flux influences the numerical dissipation. Ideally, a flux appropriate  for all Mach numbers should be used, as done e.g. in \cite{park:2005}. In order to obtain a numerical method effective for a wide range of Mach numbers, we take
\begin{equation}
    \lambda^{(n,m)} = \max\left[f\left(M_{loc}^{+,(n,m)}\right) \left(\left|\mathbf{u}^{+,(n,m)}\right| + \frac{1}{M}c^{+,(n,m)}\right), f\left(M_{loc}^{-,(n,m)}\right)\left(\left|\mathbf{u}^{-,(n,m)}\right| + \frac{1}{M}c^{-,(n,m)}\right)\right],
\end{equation}
with $M_{loc}^{\pm,(n,m)} = \frac{M\left|\mathbf{u}\right|^{\pm,(n,m)}}{c^{\pm,(n,m)}}$ and $f\left(M_{loc}\right) = \min\left(1, M_{loc}\right)$. This choice corresponds to the convex combination between a centered flux and a Rusanov flux \cite{rusanov:1962} discussed in \cite{abbate:2019}. More specifically, for a generic flux $\hat{\mathbf{F}}$, we employ
\begin{equation}
    \hat{\mathbf{F}} = \left(1 - f\left(M_{loc}\right)\right)\mathbf{F}_{c} + f\left(M_{loc}\right)\mathbf{F}_{R},
\end{equation}
with $\mathbf{F}_{c}$ and $\mathbf{F}_{R}$ denoting the centered flux and the Rusanov flux, respectively. Hence, for $M_{loc} \approx 1$, we obtain the Rusanov flux, whereas for $M_{loc} \ll 1$, we obtain a local Lax-Friedrichs flux. Analogously, the energy equation in \eqref{eq:stage_euler} can be expressed as
\begin{equation}
    \mathbf{C}^{(n,l)}\mathbf{U}^{(n,l)} + \mathbf{D}^{(n,l)}\mathbf{P}^{(n,l)} = \mathbf{G}^{(n,l)},
\end{equation}
with
\begin{eqnarray}
    C_{ij}^{(n,l)} &=& \sum_{K \in \mathcal{T}_{\mathcal{H}}}\int_{K} -\tilde{a}_{ll}\Delta t h^{(n,l)}\rho^{(n,l)}\boldsymbol{\varphi_{j}} \cdot \grad\Psi_{i}d\Omega + \sum_{\Gamma \in \mathcal{E}}\int_{\Gamma}\tilde{a}_{ll}\Delta t\left\{\left\{h^{(n,l)}\rho^{(n,l)}\boldsymbol{\varphi_{j}}\right\}\right\} \cdot \left[\left[\Psi_{i}\right]\right]d\Sigma \\
    D_{ij}^{(n,l)} &=& \sum_{K \in \mathcal{T}_{\mathcal{H}}}\int_{K}\rho^{(n,l)}e^{(n,l)}(\rho^{(n,l)}, \Psi_{j})\Psi_{i}d\Omega \\
    G_{i}^{(n,l)} &=& \sum_{m=1}^{l-1}\sum_{K \in \mathcal{T}_{\mathcal{H}}}\int_{K} \rho^{(n,l)}E^{(n,l)}\psi_{i}d\Omega \nonumber \\
    &+& \sum_{m=1}^{l-1}\sum_{K \in \mathcal{T}_{\mathcal{H}}}\int_{K}a_{lm}\Delta tM^{2}\left(k^{(n,m)}\rho^{(n,m)}\mathbf{u}^{(n,m)}\right) \cdot \grad\Psi_{i}d\Omega \nonumber \\
    &+& \sum_{m=1}^{l-1}\sum_{K \in \mathcal{T}_{\mathcal{H}}}\int_{K} \tilde{a}_{lm}\Delta t\left(h^{(n,m)}\rho^{(n,m)}\mathbf{u}^{(n,m)}\right) \cdot \grad\Psi_{i}d\Omega \nonumber \\
    &-& \sum_{m=1}^{l-1}\sum_{\Gamma \in \mathcal{E}}\int_{\Gamma} a_{lm}\Delta tM^{2}\left\{\left\{k^{(n,m)}\rho^{(n,m})\mathbf{u}^{(n,m)}\right\}\right\} \cdot \left[\left[\Psi_{i}\right]\right]d\Sigma \nonumber \\
    &-& \sum_{m=1}^{l-1}\sum_{\Gamma \in \mathcal{E}} \int_{\Gamma} \tilde{a}_{lm}\Delta t\left\{\left\{h^{(n,m)}\rho^{(n,m})\mathbf{u}^{(n,m)}\right\}\right\} \cdot \left[\left[\Psi_{i}\right]\right]d\Sigma \nonumber \\
    &-& \sum_{m=1}^{l-1}\sum_{\Gamma \in \mathcal{E}}\int_{\Gamma} a_{lm}\Delta tM^{2}\frac{\lambda^{(n,m)}}{2}\left[\left[\rho^{(n,m)}k^{(n,m)}\right]\right] \cdot \left[\left[\Psi_{i}\right]\right]d\Sigma 
    \nonumber \\
    &-& \sum_{m=1}^{l-1}\sum_{\Gamma \in \mathcal{E}}\int_{\Gamma} \tilde{a}_{lm}\Delta t\frac{\lambda^{(n,m)}}{2}\left[\left[\rho^{(n,m)}e^{(n,m)}\right]\right] \cdot \left[\left[\Psi_{i}\right]\right]d\Sigma \nonumber \\
    &-& \sum_{K \in \mathcal{T}_{h}}\int_{K}M^{2}\rho^{(n,l)}k^{(n,l)}\Psi_{i}d\Omega
    -\sum_{\Gamma \in \mathcal{E}}\int_{\Gamma}\tilde{a}_{ll}\Delta t\frac{\lambda^{(n,l)}}{2}\left[\left[\rho^{(n,l)}e^{(n,l)}\right]\right] \cdot \left[\left[\Psi_{i}\right]\right]d\Sigma.
\end{eqnarray}
Notice that, the upwind flux has been slightly modified with respect to the one employed in \cite{orlando:2022b}, so as to guarantee the preservation of uniform velocity and pressure fields (see the discussion in \cite{orlando:2023c}). Formally, one can derive
\begin{equation}
    \mathbf{U}^{(n,l)} = \left(\mathbf{A}^{(n,l)}\right)^{-1}\left(\mathbf{F}^{(n,l)} - \mathbf{B}^{(n,l)}\mathbf{P}^{(n,l)}\right), 
\end{equation}
so as to obtain
\begin{equation}
    \mathbf{D}^{(n,l)}\mathbf{P}^{(n,l)} + \mathbf{C}^{(n,l)}\left(\mathbf{A}^{(n,l)}\right)^{-1}\left(\mathbf{F}^{(n,l)} - \mathbf{B}^{(n,l)}\mathbf{P}^{(n,l)}\right) = \mathbf{G}^{(n,l)}.
\end{equation}
The above system can be solved following the fixed point procedure described in \cite{dumbser:2016b, orlando:2022b}. More specifically, setting $\mathbf{P}^{(n,l,0)} = \mathbf{P}^{(n,l-1)}$, one solves for $L = 0, \dots, \tilde{L}$
\begin{equation}
    \left(\mathbf{D}^{(n,l,L)} - \mathbf{C}^{(n,l,L)}\left(\mathbf{A}^{(n,l)}\right)^{-1}\mathbf{B}^{(n,l)}\right)\mathbf{P}^{(n,l,L+1)} = \mathbf{G}^{(n,l,L)} - \mathbf{C}^{(n,l,L)}\left(\mathbf{A}^{(n,l)}\right)^{-1}\mathbf{F}^{(n,l,L)}
\end{equation}
and then updates the velocity solving
\begin{equation}
    \mathbf{A}^{(n,l)}\mathbf{U}^{(n,l,L)} = \mathbf{F}^{(n,l,L)} - \mathbf{B}^{(n,l)}\mathbf{P}^{(n,l,L+1)}.
\end{equation}
Notice that, as discussed for the time discretization in Section \ref{sec:ap_num}, the employed spatial discretization is not TVD for $r > 0$. A discussion of possible approaches to overcome this issue is out of the scope of the work. However, a number of approaches have been proposed to obtain essentially monotone schemes using high order DG methods, see e.g. \cite{dumbser:2016a, orlando:2023a}.

The DG method naturally allows for high-order accuracy. However, as discussed in \cite{jung:2024}, its accuracy in the very low Mach regime depends on the numerical flux and on the shape of the elements. More specifically, a mesh of triangular/tetrahedral elements is needed to guarantee accuracy at all Mach numbers. A low Mach number fix for the Euler equations resolved employing the finite volume method on Cartesian grids was presented in \cite{barsukow:2021}. We will further discuss this point in Section \ref{ssec:isentropic_vortex}.

\section{Numerical results}
\label{sec:num}

The analysis outlined in Sections \ref{sec:ap_num} and \ref{sec:space} is now validated in a number of  benchmarks covering the $M < 1$ and $M \ll 1$ regimes. The implementation is carried out in the framework of the \texttt{deal.II} library \cite{arndt:2023, bangerth:2007}.
We use a time discretization based on the third order IMEX scheme presented in \cite{kennedy:2003}, for which the coefficients of both the explicit and implicit methods are reported in the Butcher tableaux Table \ref{tab:rk3_butch}. Moreover, in order to assess the convergence properties of the method and to exploit the high-order accuracy provided by the DG method, we also consider in Section \ref{ssec:isentropic_vortex} the fourth order time discretization scheme proposed in \cite{calvo:2001}, for which the coefficients of the explicit and of the implicit companion method are reported in Table \ref{tab:rk4_butch_ARS}. Notice that the implicit method of the fourth order time discretization scheme is of type ARS. One can easily check that the implicit companion methods of both schemes are stiffly-accurate. Hence, the implicit scheme in Table \ref{tab:rk4_butch_ARS} is \texttt{L}-stable. For what concerns the implicit scheme in Table \ref{tab:rk3_butch}, relation \eqref{eq:Lstab_type_II} leads to $R_{\infty} = 0$ and therefore it is also \texttt{L}-stable.

\begin{table}[pos=H]
    \begin{minipage}{0.9\textwidth}
	\begin{center}
            \begin{tabular}{c|cccc}
                $0$ & $0$ & $0$ & $0$ & $0$ \\ [2mm]
                $\frac{1767732205903}{2027836641118}$ & $\frac{1767732205903}{2027836641118}$ & $0$ & $0$ & $0$ \\ [2mm]
                $\frac{3}{5}$ & $\frac{5535828885825}{10492691773637}$ & $\frac{788022342437}{10882634858940}$ & $0$ & $0$ \\ [2mm]
                $1$ & $\frac{6485989280629}{16251701735622}$ & $\frac{-4246266847089}{9704473918619}$ & $	\frac{10755448449292}{10357097424841}$ & $0$ \\ [2mm]
                \hline \\ [-0.35cm]
	        & $\frac{1471266399579}{7840856788654}$ & 
                $\frac{-4482444167858}{7529755066697}$ & $\frac{11266239266428}{11593286722821}$ & $\frac{1767732205903}{4055673282236}$
	    \end{tabular}
        \end{center}
    \end{minipage}
    \vskip 0.3cm
    \begin{minipage}{0.9\textwidth}
	\begin{center}
            \begin{tabular}{c|cccc}
				$0$ & $0$ & $0$ & $0$ & $0$ \\ [2mm]
                $\frac{1767732205903}{2027836641118}$ & $\frac{1767732205903}{4055673282236}$ & $\frac{1767732205903}{4055673282236}$ & $0$ & $0$ \\ [2mm]
                $\frac{3}{5}$ & $\frac{2746238789719}{10658868560708}$ &
                $\frac{-640167445237}{6845629431997}$ & $\frac{1767732205903}{4055673282236}$ & $0$ \\ [2mm]
                $1$ & $\frac{1471266399579}{7840856788654}$ & 
                $\frac{-4482444167858}{7529755066697}$ & $\frac{11266239266428}{11593286722821}$ & $\frac{1767732205903}{4055673282236}$ \\
                \hline \\ [-0.35cm]
                & $\frac{1471266399579}{7840856788654}$ & 
                $\frac{-4482444167858}{7529755066697}$ & $\frac{11266239266428}{11593286722821}$ & $\frac{1767732205903}{4055673282236}$
		\end{tabular}
	\end{center}
    \end{minipage}
    \caption{\it Butcher tableaux of the third order time discretization scheme. Top: explicit method. Bottom: implicit method.}
    \label{tab:rk3_butch}
\end{table}

\begin{table}[pos=H]
    \begin{minipage}{0.9\textwidth}
	\begin{center}
            \begin{tabular}{c|cccccc}
		      $0$ & $0$ & $0$ & $0$ & $0$ & $0$ & $0$\\[1.5mm]
                $\frac{1}{4}$ & $\frac{1}{4}$ & $0$ & $0$ & $0$ & $0$ & $0$ \\[1.5mm]
                $\frac{3}{4}$ & $-\frac{1}{4}$ & $1$ & $0$ & $0$ & $0$ & $0$ \\[1.5mm]
                $\frac{11}{20}$ & $-\frac{13}{100}$ & $\frac{43}{75}$ & $\frac{8}{75}$ & $0$ & $0$ & $0$ \\[1.5mm]
                $\frac{1}{2}$ & $-\frac{6}{85}$ & $\frac{42}{85}$ & $\frac{179}{1360}$ & $-\frac{15}{272}$ & $0$ & $0$ \\[1.5mm]
                $1$ & $0$ & $\frac{79}{24}$ & $-\frac{5}{8}$ & $\frac{25}{2}$ & $-\frac{85}{6}$ & $0$ \\
                \hline \\ [-0.35cm]
                & $0$ &	$\frac{25}{24}$ & $-\frac{49}{48}$ & $\frac{125}{16}$ & $-\frac{85}{12}$ & $\frac{1}{4}$
		\end{tabular}
	\end{center}
    \end{minipage}
    \vskip 0.3cm
    \begin{minipage}{0.9\textwidth}
	\begin{center}
            \begin{tabular}{c|cccccc}
				$0$ & $0$ & $0$ & $0$ & $0$ & $0$ & $0$\\[1.5mm]
                $\frac{1}{4}$ & $0$ & $\frac{1}{4}$ & $0$ & $0$ & $0$ & $0$ \\[1.5mm]
                $\frac{3}{4}$ & $0$ & $\frac{1}{2}$ & $\frac{1}{4}$ & $0$ & $0$ & $0$ \\[1.5mm]
                $\frac{11}{20}$ & $0$ & $\frac{17}{50}$ & $-\frac{1}{25}$ & $\frac{1}{4}$ & $0$ & $0$ \\[1.5mm]
                $\frac{1}{2}$ & $0$ & $\frac{371}{1360}$ & $-\frac{137}{2720}$ & $\frac{15}{544}$ & $\frac{1}{4}$ & $0$ \\[1.5mm]
                $1$ & $0$ &	$\frac{25}{24}$ & $-\frac{49}{48}$ & $\frac{125}{16}$ & $-\frac{85}{12}$ & $\frac{1}{4}$ \\[1.5mm]
                \hline\\[-0.35cm]
                & $0$ &	$\frac{25}{24}$ & $-\frac{49}{48}$ & $\frac{125}{16}$ & $-\frac{85}{12}$ & $\frac{1}{4}$
		\end{tabular}
	\end{center}
    \end{minipage}
    \caption{\it Butcher tableaux of the fourth order time discretization scheme. Top: explicit method. Bottom: implicit method.}
    \label{tab:rk4_butch_ARS}
\end{table}

Notice that both the schemes are of type II. Results employing numerical schemes of type I can be found in \cite{orlando:2025}. We set $\mathcal{H} = \min\left\{\text{diam}\left(K\right) | K \in \mathcal{T}_{\mathcal{H}}\right\}$ and we define two Courant numbers, one based on the speed of sound (acoustic Courant number), denoted by $C$, and one based on the local velocity of the flow (advective Courant number), denoted by $C_{u}$:
\begin{equation}
    C = \frac{1}{M}rc\frac{\Delta t}{\mathcal{H}}\sqrt{d} \qquad C_{u} = ru\frac{\Delta t}{\mathcal{H}}\sqrt{d},
\end{equation}
where $c$ is the speed of sound and $u$ is the magnitude of the flow velocity. Recall that $r$ denotes the polynomial degree of the space discretization. For the tests using the ideal gas law \eqref{eq:ideal_gas}, the value $\gamma = 1.4$ is employed.

\subsection{Isentropic vortex}
\label{ssec:isentropic_vortex}

First, we consider the isentropic vortex benchmark studied in \cite{boscheri:2021a, yee:1999, zampa:2025} using the ideal gas law \eqref{eq:ideal_gas}, for which an analytical solution is available. Following, e.g., \cite{boscheri:2021a, zeifang:2019}, the steady solution of system \eqref{eq:euler_adim} as a function of the Mach number reads as follows:
\begin{subequations}
\begin{eqnarray}
    \mathbf{u}\left(\mathbf{x},t\right) &=& \mathbf{u}\left(\mathbf{x},0\right) = M\frac{\beta}{2\pi}\exp\left(\frac{1 - \tilde{r}^{2}}{2}\right)\begin{pmatrix}
        -(y - y_{0}) \\
        x - x_{0}
    \end{pmatrix} \\
    \rho\left(\mathbf{x},t\right) &=& \rho\left(\mathbf{x},0\right) = \left(1 + \delta T\right)^{\frac{1}{\gamma - 1}} \\
    p\left(\mathbf{x},t\right) &=& p\left(\mathbf{x},0\right) = M^{2}\left(1 + \delta T\right)^{\frac{\gamma}{\gamma - 1}} = M^{2}\rho^{\gamma},
\end{eqnarray}
\end{subequations}
where $\tilde{r}^{2} = \left(x - x_{0}\right)^{2} + \left(y - y_{0}\right)^{2}$, with $x_{0}$ and $y_{0}$ denoting the coordinates of the center of the vortex. Moreover, $\beta$ is the vortex strength and the temperature perturbation $\delta T$ is given by
\begin{equation}
    \delta T = -M^{2}\frac{\gamma - 1}{\gamma}\frac{\beta^{2}}{8\pi^{2}}\exp\left(1 - \tilde{r}^{2}\right).
\end{equation}
A travelling vortex configuration can be found instead in \cite{orlando:2025}, to which we refer also for the impact of different time discretization strategies. To avoid problems with the definition of the boundary conditions, we choose a sufficiently large domain $\Omega = \left(-10, 10\right)^{2}$, with $x_{0} = y_{0} = 0$, and periodic boundary conditions. Finally, we set $\beta = 5$ and $T_{f} = 10$. The purpose of this test is twofold. First, we assess the convergence properties of the method, employing the IMEX schemes in Table \ref{tab:rk3_butch} and Table \ref{tab:rk4_butch_ARS}. Next, we verify the asymptotic expansion in the small Mach number limit outlined in Section \ref{ssec:ap_single_length}. We set $M = 10^{-3}$ for the convergence analysis, which is performed at fixed acoustic Courant number $C \approx 3.5$. We report results for the case of polynomial degree $r = 2$ in combination with the third order time discretization scheme (Table \ref{tab:vortex_convergence_third_order}) and for the case of polynomial degree $r = 3$ in combination with the fourth order time discretization scheme (Table \ref{tab:vortex_convergence_fourth_order}). The expected convergence rate is achieved for the third order method, whereas an order reduction is experienced for the fourth order method as the resolution increases and the time step decreases. Since the solution is steady, this is likely related to an early manifestation of a low Mach inaccuracy (see the discussion below).

Next, we analyze the behaviour in the $M \to 0$ limit. We employ the third order time discretization scheme in Table \ref{tab:rk3_butch} with $r = 2$ and $N_{el} = 120$ elements along each coordinate direction. Until $M = 10^{-5}$, the density fluctuations scale as $\mathcal{O}(M^{2})$ and the divergence of the velocity field scales as $\mathcal{O}(M)$ (Table \ref{tab:vortex_Mach_third_order}), as expected \cite{jung:2024}. The convergence with respect to $M$ of the divergence of the velocity field deserves some comments. Since $\dive\bar{\mathbf{u}} = 0$, the first order term of the energy equation for a steady state solution reduces to $\bar{p}\dive\mathbf{u}^{'} = 0$ (see \eqref{eq:energy_limit_first_order_two_scale}). The initial velocity field is indeed solenoidal and a quadratic convergence with respect to $M$ could be therefore expected. However, the divergence-free property is not imposed pointwise and, since a term proportional to $M$ is present in the velocity field, we obtain a linear scaling with respect to $M$ for $\dive\mathbf{u}$ (Table \ref{tab:vortex_Mach_third_order}). A quadratic convergence was obtained recently in \cite{zampa:2025} for the Taylor-Green vortex. This is likely related to the use of compatible finite elements which allow the imposition of the divergence-free property for the initial datum, so as to observe the quadratic convergence rate predicted by the asymptotic expansion of the continuous model. Indeed, in our framework, the error associated to $\dive\mathbf{u}$ is basically constant in time and it is therefore related to the interpolation of the initial datum into the employed finite element space.

At $M = 10^{-6}$, we observe a small degradation for the density fluctuations. This is likely related to the well-known \textit{inf-sup} stability condition for DG discretizations of incompressible flows \cite{toselli:2002}. Indeed, we have verified that a slightly improved scaling of the density fluctuations is obtained employing a polynomial degree $r + 1$ for the velocity field, i.e. third order polynomials (Table \ref{tab:vortex_Mach_third_order}). Moreover, it is worth to remark that at $M = 10^{-6}$, the density and the pressure are basically constant and the $L^{2}$ error for the pressure is below the machine precision, so that round-off errors play a relevant role. Indeed, as remarked in \cite{boscheri:2021b}, the use of quadruple precision is crucial to maintain the theoretical scaling at  very small Mach numbers.

A similar behaviour is experienced employing the fourth order time discretization scheme in Table \ref{tab:rk4_butch_ARS} with $r = 3$ and $N_{el} = 80$. Here we notice that the density fluctuations start scaling as $\mathcal{O}(M)$ from $M = 10^{-4}$ (Table \ref{tab:vortex_Mach_fourth_order}). The loss of low Mach accuracy from $M = 10^{-4}$ is also related to well-known order reduction phenomenon experienced for very stiff problems using high-order time discretization methods \cite{kennedy:2019, wanner:1996} (see also \cite{orlando:2025}). Indeed, we have verified that, using the third order time discretization scheme with polynomial degree $r = 3$, the correct scaling is established up to $M = 10^{-5}$. The use of polynomial degree $r + 1 = 4$ for the velocity field allows us to recover the correct scaling up to $M = 10^{-5}$ (Table \ref{tab:vortex_Mach_fourth_order}). At $M = 10^{-6}$, for which round-off errors play a major role, a degradation is still experienced.

It has to be recalled that well known issues arise using quadrilateral elements for strongly subsonic flows. The seminal work of Guillard and Viozat \cite{guillard:1999} showed through an asymptotic analysis of the first order Roe scheme that a pressure term of order $\mathcal{O}(M)$ appears on Cartesian meshes as $M \to 0$. A number of fixes for numerical fluxes that preserve contact discontinuities (HLLC, Roe, etc...) have been proposed in the literature \cite{dellacherie:2016, dellacherie:2010b, galie:2024, rieper:2011}. They have been developed in the framework of the Finite Volume method, but they can in principle be extended to the DG method. The high-order accuracy of the DG method typically counterbalances the lack of low Mach accuracy for strongly subsonic flows. A loss of accuracy in this limit was already reported in \cite{bassi:2009}. In the recent work of Jung and Perrier \cite{jung:2024}, the authors show that the DG method employing numerical fluxes that preserve contact discontinuities is low Mach number accurate using triangular elements, while the same does not hold for quadrilateral elements. Moreover, as reported in \cite{jung:2022}, low Mach fixes are similar to schemes based on specific IMEX time discretizations. As an example, the fix proposed in \cite{dellacherie:2010b} imposes a zero velocity jump in the artificial viscosity term, so as to obtain a centering of the pressure gradient in the momentum equation. The method presented in this work uses a centered flux for the quantities defined implicitly, including the pressure gradient (see Section \ref{sec:space}). Since the focus of this work is to show the AP property of a general class of IMEX time discretizations, we do not investigate further these issues, that are related to the spatial discretization. It can be observed, however, that the high-order accuracy of the DG method allows to simulate correctly flows down to $M=10^{-4}-10^{-5}.$ This is lower than the typical values of the Mach number for fluids like water, that are modelled as incompressible in most realistic applications. Moreover, in \cite{thornber:2008} (see also \cite{hope:2023}), the authors show that inaccuracies of standard Godunov schemes at low Mach number are linked to spurious entropy generation. Entropy-stable DG schemes, as those developed, e.g., in \cite{gassner:2021, waruszewski:2022}, could therefore improve the low Mach accuracy. The use of entropy-stable DG methods and of exterior calculus and compatible finite element to further improve the low Mach accuracy of the spatial discretization will be the the focus of future work.

\begin{table}[pos=H]
    \begin{tabularx}{0.95\linewidth}{XXXXXXX}
	\toprule
        $N_{el}$ & $L^{2}$ rel. error $\mathbf{u}$ & $L^{2}$ rate $\mathbf{u}$ & $L^{2}$ rel. error $\rho$ & $L^{2}$ rate $\rho$ & $L^{2}$ rel. error $p$ & $L^{2}$ rate $p$ \\
        \midrule
        $15$ & $\num{4.70e-2}$ & & $\num{9.91e-8}$ & & $\num{1.31e-7}$ & \\
        \midrule
        $30$ & $\num{5.06e-3}$ & $3.2$ & $\num{1.77e-9}$ & $5.8$ & $\num{2.27e-9}$ & $5.9$ \\
        \midrule
        $60$ & $\num{6.42e-4}$ & $3.0$ & $\num{8.18e-11}$ & $4.4$ & $\num{3.51e-11}$ & $6.0$ \\
        \midrule
        $120$ & $\num{8.07e-5}$ & $3.0$ & $\num{1.07e-11}$ & $2.9$ & $\num{2.03e-12}$ & $4.1$ \\
        \midrule
        $240$ & $\num{1.02e-5}$ & $3.0$ & $\num{1.92e-12}$ & $2.5$ & $\num{2.46e-13}$ & $3.0$ \\
	\bottomrule
    \end{tabularx}
    \caption{Convergence analysis for the isentropic vortex test case using the time discretization scheme in Table \ref{tab:rk3_butch} together with polynomial degree $r = 2$. Here, $N_{el}$ denotes the number of elements along each direction.}
    \label{tab:vortex_convergence_third_order}
\end{table}

\begin{table}[pos=H]
    \begin{tabularx}{0.95\linewidth}{XXXXXXX}
        \toprule
        $N_{el}$ & $L^{2}$ rel. error $\mathbf{u}$ & $L^{2}$ rate $\mathbf{u}$ & $L^{2}$ rel. error $\rho$ & $L^{2}$ rate $\rho$ & $L^{2}$ rel. error $p$ & $L^{2}$ rate $p$ \\
        \midrule
        $10$ & $\num{3.27e-2}$ & & $\num{8.67e-7}$ & & $\num{1.21e-6}$ & \\
        \midrule
        $20$ & $\num{2.18e-3}$ & $3.9$ & $\num{2.97e-8}$ & $4.9$ & $\num{4.15e-8}$ & $4.9$ \\
        \midrule
        $40$ & $\num{1.47e-4}$ & $3.9$ & $\num{2.58e-9}$ & $3.5$ & $\num{3.62e-9}$ & $3.5$ \\
        \midrule
        $80$ & $\num{1.32e-5}$ & $3.5$ & $\num{4.76e-10}$ & $2.4$ & $\num{6.66e-10}$ & $2.4$ \\
        \midrule
        $160$ & $\num{1.82e-6}$ & $2.9$ & $\num{6.20e-11}$ & $2.9$ & $\num{8.68e-11}$ & $2.9$ \\
	\bottomrule
    \end{tabularx}
    \caption{Convergence analysis for the isentropic vortex test case using the time discretization scheme in Table \ref{tab:rk4_butch_ARS} together with polynomial degree $r = 3$. Here, $N_{el}$ denotes the number of elements along each direction.}
    \label{tab:vortex_convergence_fourth_order}
\end{table}

\begin{table}[pos=H]
    \centering  
    \begin{tabularx}{0.95\linewidth}{XXXXXXX}
	\toprule
        $M$ & $L^{2}$ norm $\dive\mathbf{u}$ & Rate $\dive\mathbf{u}$ & $L^{2}$ norm $\grad\rho$ & Rate $\grad\rho$ & $L^{2}$ norm $\grad\rho (Q_{3} - Q_{2})$ & Rate $\grad\rho (Q_{3} - Q_{2})$ \\
	\midrule
	$10^{-1}$ & $\num{3.52e-4}$ & & $\num{1.09e-2}$ & & &\\
	\midrule
        $10^{-2}$ & $\num{3.46e-5}$ & $1.0$ & $\num{1.09e-4}$ & $2.0$ & & \\
	\midrule
        $10^{-3}$ & $\num{3.44e-6}$ & $1.0$ & $\num{1.09e-6}$ & $2.0$ & & \\
	\midrule
        $10^{-4}$ & $\num{3.44e-7}$ & $1.0$ & $\num{1.10e-8}$ & $2.0$ & & \\
        \midrule
        $10^{-5}$ & $\num{3.44e-8}$ & $1.0$ & $\num{1.29e-10}$ & $1.9$ & & \\
        \midrule
        $10^{-6}$ & $\num{3.47e-9}$ & $1.0$ & $\num{9.66e-12}$ & $1.1$ & $\num{8.34e-12}$ & $1.2$ \\
	\bottomrule
    \end{tabularx}
    \caption{Mach number scaling of the density fluctuations and of the divergence of the velocity field for the isentropic vortex test case. The results are obtained using the third order time discretization scheme in Table \ref{tab:rk3_butch} together with polynomial degree $r = 2$ and $N_{el} = 120$. The last two columns report the results obtained using polynomial degree $r + 1 = 3$ for the velocity and polynomial degree $r = 2$ for the other variables. We recall that $Q_{r}$ denotes polynomial spaces  of degree $r$ for quadrilateral elements.}
    \label{tab:vortex_Mach_third_order}
\end{table}

\begin{table}[pos=H]
    \centering  
    \begin{tabularx}{0.95\linewidth}{XXXXXXX}
	\toprule
        $M$ & $L^{2}$ norm $\dive\mathbf{u}$ & Rate $\dive\mathbf{u}$ & $L^{2}$ norm $\grad\rho$ & Rate $\grad\rho$ & $L^{2}$ norm $\grad\rho (Q_{4} - Q_{3})$ & Rate $\grad\rho (Q_{4} - Q_{3})$ \\
	\midrule
	$10^{-1}$ & $\num{1.12e-4}$ & & $\num{1.09e-2}$ & & &\\
	\midrule
        $10^{-2}$ & $\num{1.21e-5}$ & $1.0$ & $\num{1.09e-4}$ & $2.0$ & & \\
	\midrule
        $10^{-3}$ & $\num{1.25e-6}$ & $1.0$ & $\num{1.19e-6}$ & $2.0$ & & \\
	\midrule
        $10^{-4}$ & $\num{1.26e-7}$ & $1.0$ & $\num{5.05e-8}$ & $1.4$ & $\num{1.09e-8}$ & $2.0$ \\
        \midrule
        $10^{-5}$ & $\num{1.26e-8}$ & $1.0$ & $\num{4.95e-9}$ & $1.0$ & $\num{1.10e-10}$ & $2.0$ \\
        \midrule
        $10^{-6}$ & $\num{1.27e-9}$ & $1.0$ & $\num{4.97e-10}$ & $1.0$ & $\num{1.53e-11}$ & $0.9$ \\
	\bottomrule
    \end{tabularx}
    \caption{Mach number scaling of the density fluctuations and of the divergence of the velocity field for the isentropic vortex test case. The results are obtained using the fourth order time discretization scheme in Table \ref{tab:rk4_butch_ARS} together with polynomial degree $r = 3$ and $N_{el} = 80$. The last two columns report the results obtained using polynomial degree $r + 1 = 4$ for the velocity and polynomial degree $r = 3$ for the other variables. We recall that $Q_{r}$ denotes polynomial spaces of degree $r$ for quadrilateral elements.}
    \label{tab:vortex_Mach_fourth_order}
\end{table}

\subsection{Colliding acoustic pulses}
\label{ssec:colliding_pulses}

\begin{figure}[h!]
   \begin{subfigure}{0.475\textwidth}
	\centering
        \includegraphics[width=0.95\textwidth]{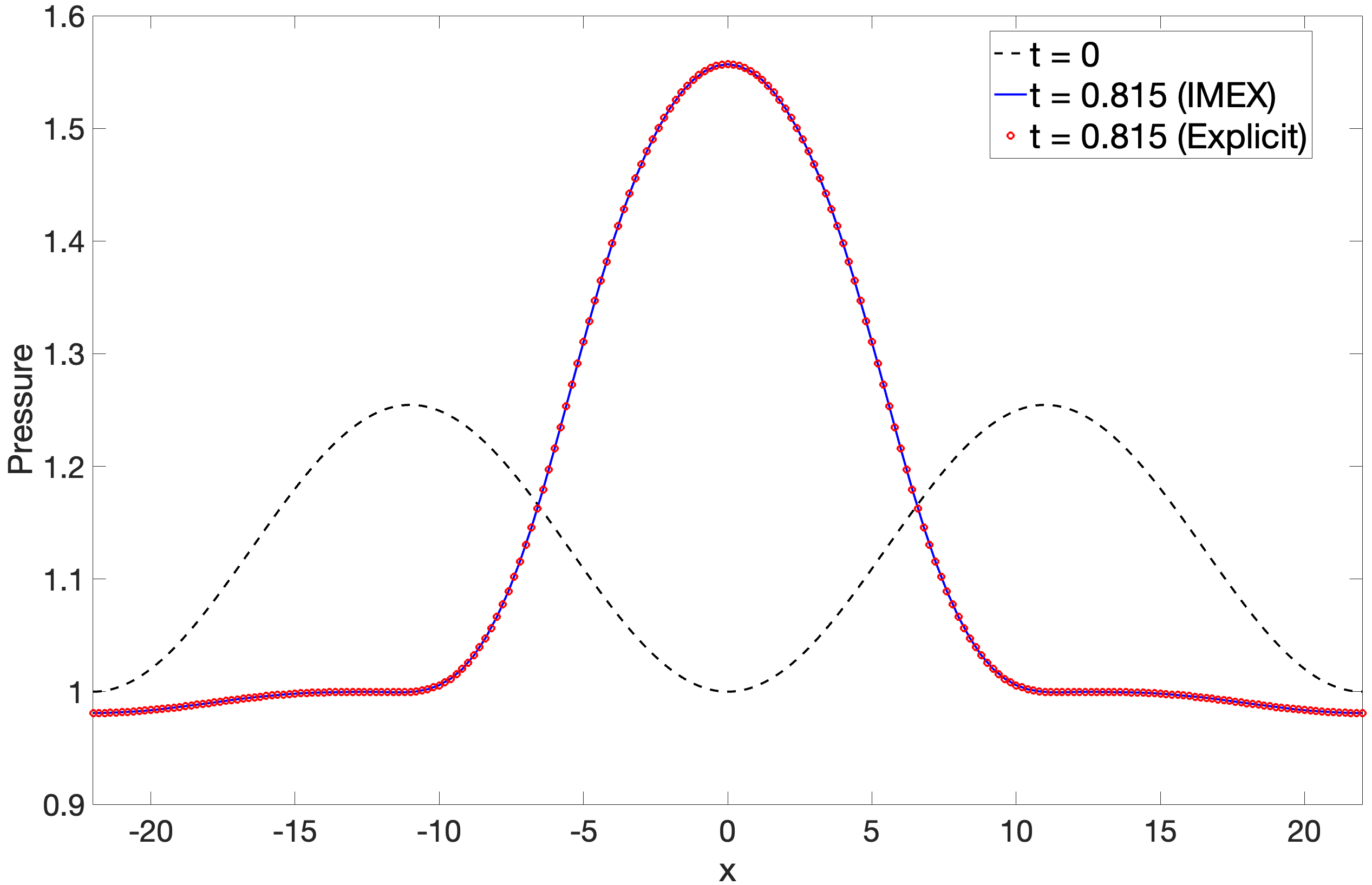}
    \end{subfigure}
    \begin{subfigure}{0.475\textwidth}
	\centering
        \includegraphics[width=0.95\textwidth]{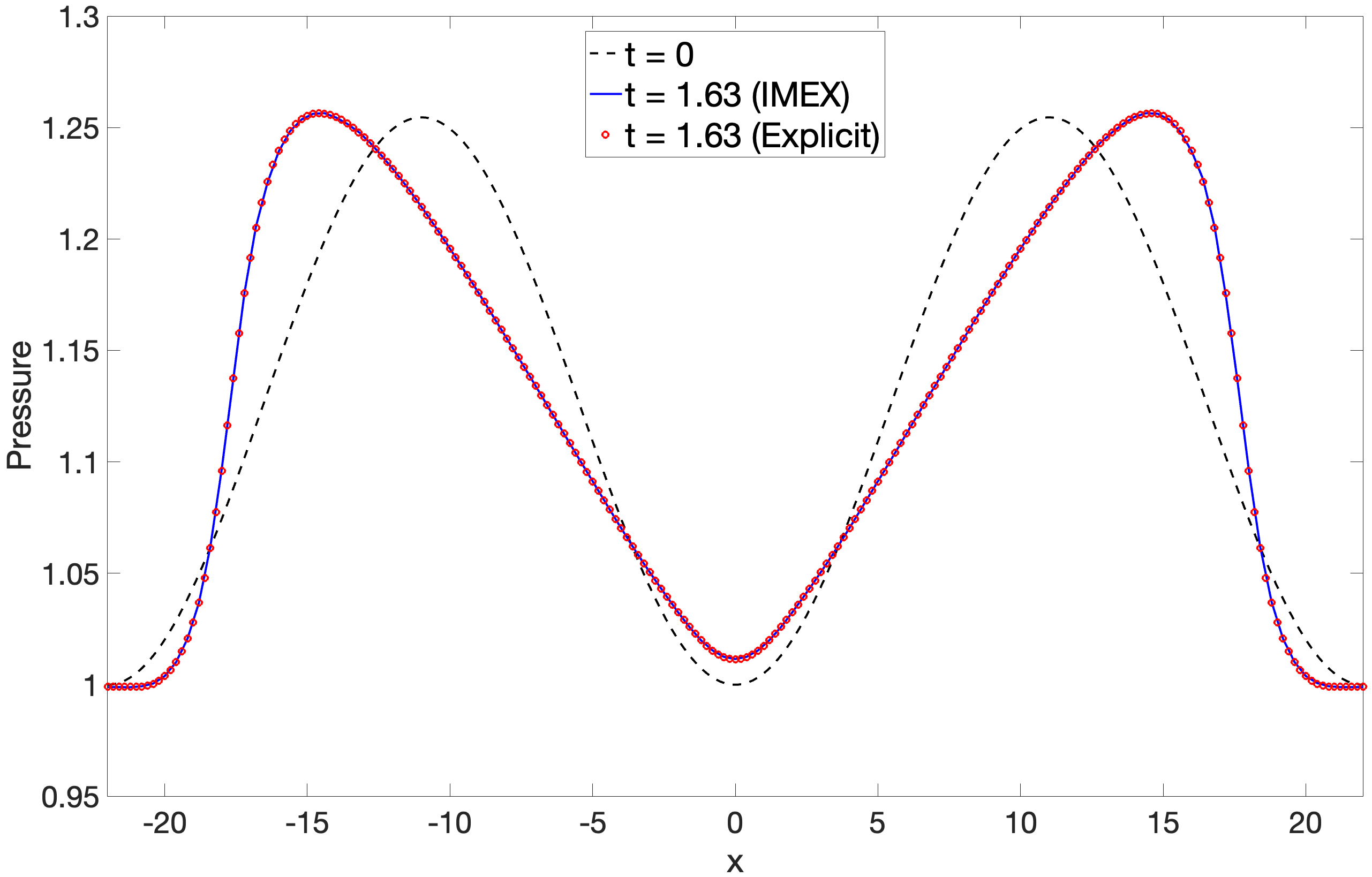}
   \end{subfigure}
    \caption{Colliding acoustic pulses test case, pressure profile. Left: $t = \frac{T_{f}}{2}$. Right: $t = T_{f}$. The initial profile is in dashed black line, the solid blue lines provide the results at the corresponding time obtained with the IMEX method at $C_{u} \approx 0.1$, whereas the red dots show the reference results obtained with the explicit method.}
    \label{fig:colliding_pulses_pressure}
\end{figure}

\begin{figure}[h!]
   \begin{subfigure}{0.475\textwidth}
	\centering
        \includegraphics[width=0.95\textwidth]{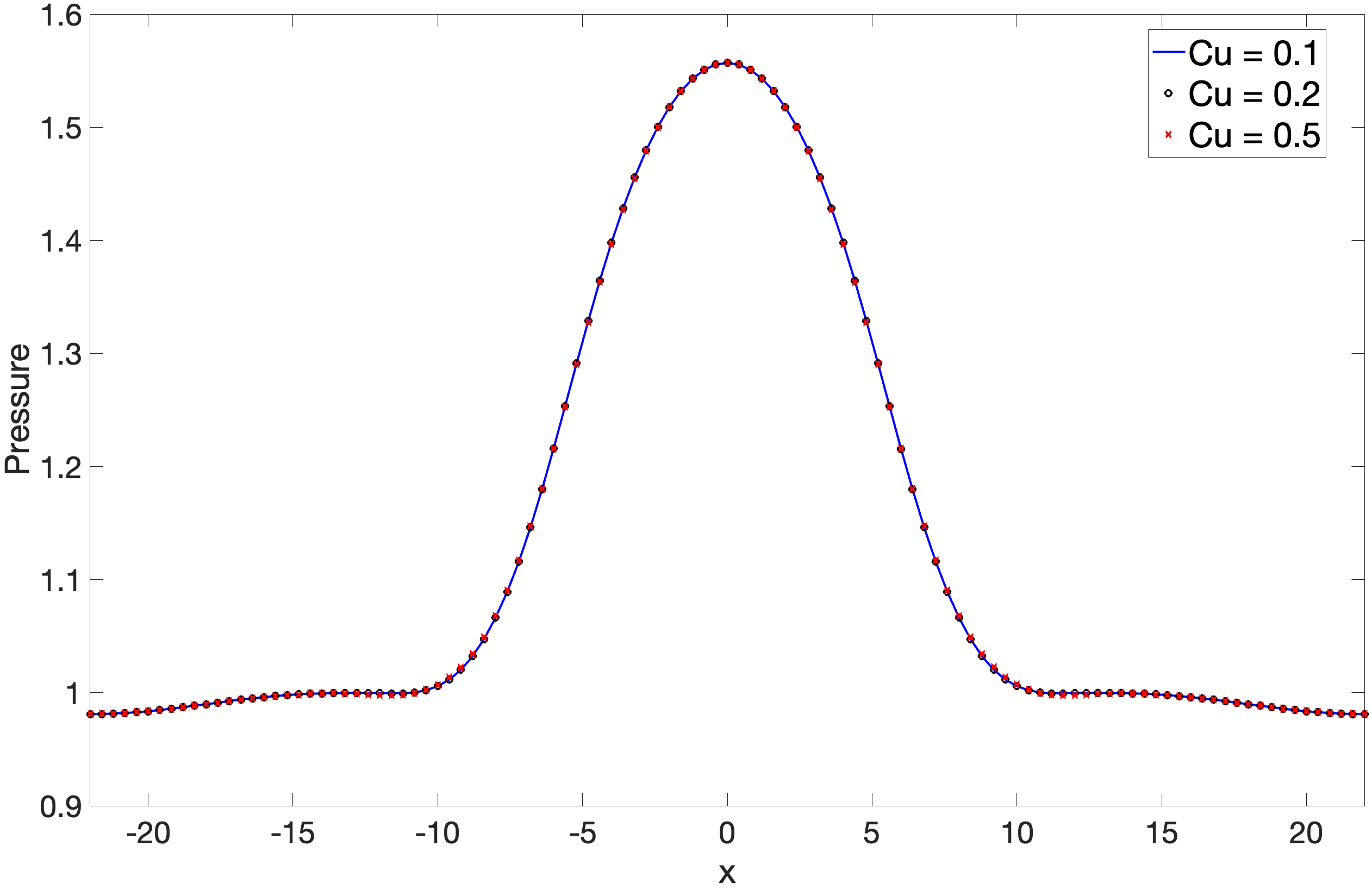}
    \end{subfigure}
    \begin{subfigure}{0.475\textwidth}
	\centering
        \includegraphics[width=0.95\textwidth]{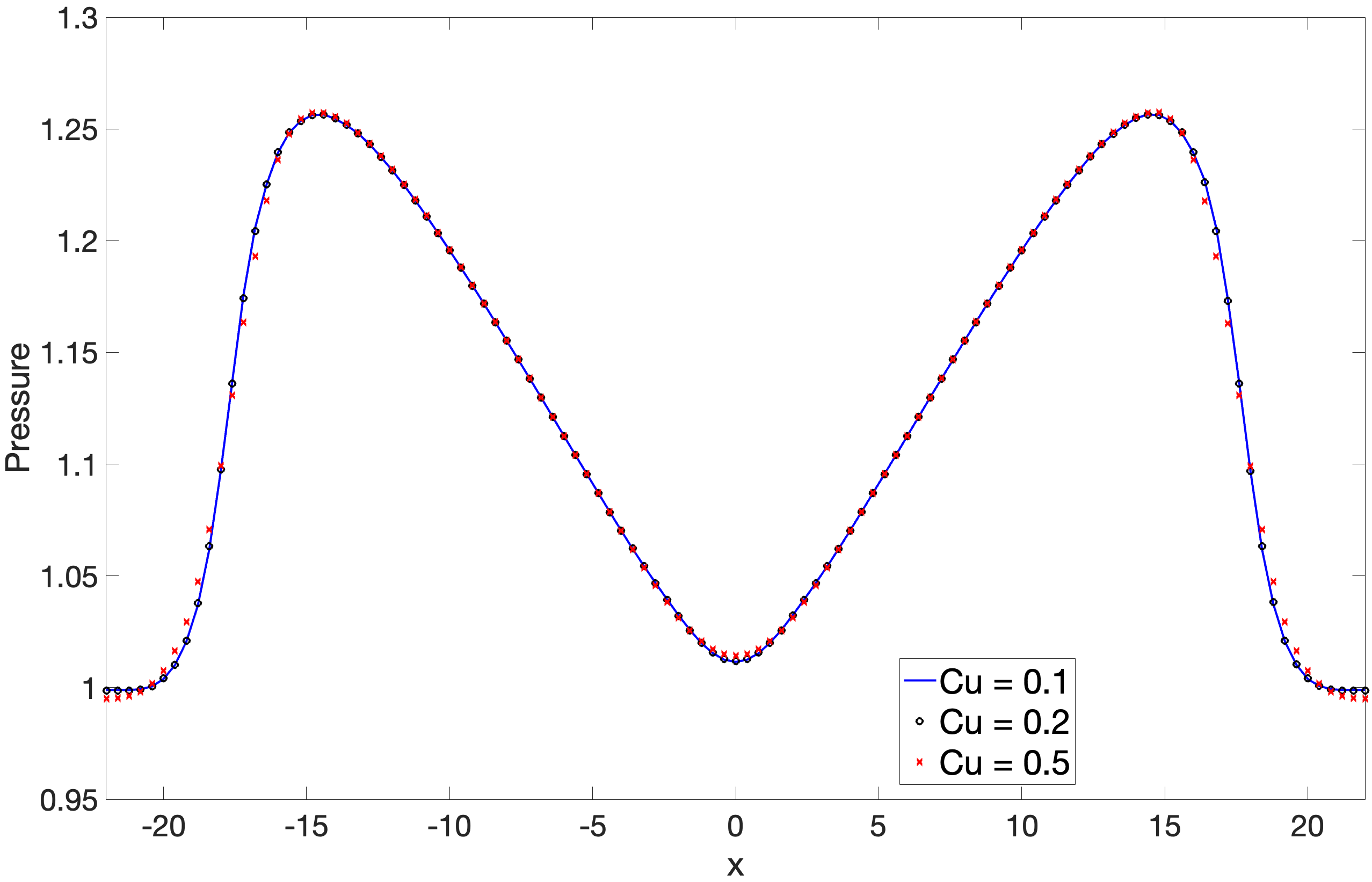}
   \end{subfigure}
    \caption{Colliding acoustic pulses test case, pressure profile. Comparison of the IMEX method employing different time step. Left: $t = \frac{T_{f}}{2}$. Right: $t = T_{f}$. The solid blue lines provide the results at the obtained at $C_{u} \approx 0.1$, the black dots show the results obtained at $C_{u} \approx 0.2$, whereas the red crosses report the results obtained at $C_{u} \approx 0.5$.}
    \label{fig:colliding_pulses_pressure_comparison_Cu}
\end{figure}

This benchmark, proposed in \cite{klein:1995}, consists of two colliding acoustic pulses in the domain $\Omega = \left(-L, L\right)$, namely, a right-running pulse initially located in $\left(-L, 0\right)$ and a left-running pulse initially located in $\left(0, L\right)$. Following \cite{klein:1995}, we set $M = \frac{1}{11}$ and we define the half-length of the domain $L = \frac{2}{M} = 22$. Periodic boundary conditions are prescribed. The initial conditions read as follows:
\begin{subequations}
\begin{eqnarray}
    \rho\left(x, 0\right) &=& \bar{\rho}_{0} + \frac{1}{2}M \rho_{0}^{'}\left(1 - \cos\left(\frac{2\pi x}{L}\right)\right) \qquad \bar{\rho}_{0} = 0.955 \quad \rho_{0}^{'} = 2 \\
    u\left(x, 0\right) &=& -\frac{1}{2}\sgn\left(x\right)\bar{u}_{0} \left(1 - \cos\left(\frac{2\pi x}{L}\right)\right) \quad \bar{u}_{0} = 2\sqrt{\gamma} \\
    p\left(x, 0\right) &=& \bar{p}_{0} + \frac{1}{2}M p_{0}^{'}\left(1 - \cos\left(\frac{2\pi x}{L}\right)\right) \qquad \bar{p}_{0} = 1 \quad p_{0}^{'} = 2\gamma
\end{eqnarray}
\end{subequations}
The final time is $T_{f} = 1.63$. We consider a number of elements $N_{el} = 55$ with $r = 2$, i.e. polynomial degree of order 2, whereas the time step is $\Delta t = 1.63 \times 10^{-2}$, leading to a maximum advective Courant number $C_{u} \approx 0.1$ and a maximum acoustic Courant number $C \approx 0.56$. A reference solution is computed using the explicit third order strong stability preserving (SSP) scheme described in \cite{gottlieb:2001}, to which we refer for all the details. We employ $N_{el} = 880$ elements with $\Delta t = 2.54687 \times 10^{-4}$, which corresponds to an acoustic Courant number $C \approx 0.14$. The pressure profiles at $t = \frac{T_{f}}{2} = 0.815$ and $t = T_{f}$ are in agreement with the reference results and with the results present in the literature \cite{cordier:2012, klein:1995, noelle:2014} (Figure \ref{fig:colliding_pulses_pressure}). One can easily notice that at $t = \frac{T_{f}}{2}$, the two pulses are superposed and the pressure reaches its maximum value. At $t = T_{f}$, the pulses are separated from each other and assume almost the initial configuration. However, as explained in \cite{cordier:2012, klein:1995}, weakly nonlinear acoustic effects start steepening the pulses and distort the final profile, since shocks are beginning to form around $x = \pm 18.5$. We also compare the accuracy of the IMEX scheme for increasing Courant numbers. More specifically, we consider $\Delta t = 3.26 \times 10^{-2}$ and $\Delta t = 8.15 \times 10^{-2}$, which lead to $C_{u} \approx 0.2, C \approx 1.16$ and $C_{u} \approx 0.5, C \approx 2.9$, respectively. For larger time step, stability restrictions imposed by the explicit component of the IMEX scheme arise \cite{cockburn:2001}. Moreover, as we will discuss later on, for very large values of the acoustic Courant number, the profile of the pulses is damped. One can easily notice that an excellent agreement is established and we can correctly capture the acoustic pulses also at acoustic Courant number which are moderately higher than $1$ (Figure \ref{fig:colliding_pulses_pressure_comparison_Cu}). Small differences arise at $t = T_{f}$, where the pulses start steepening and shocks are beginning to form.

\begin{figure}[h!]
   \begin{subfigure}{0.475\textwidth}
		\centering
        \includegraphics[width=0.95\textwidth]{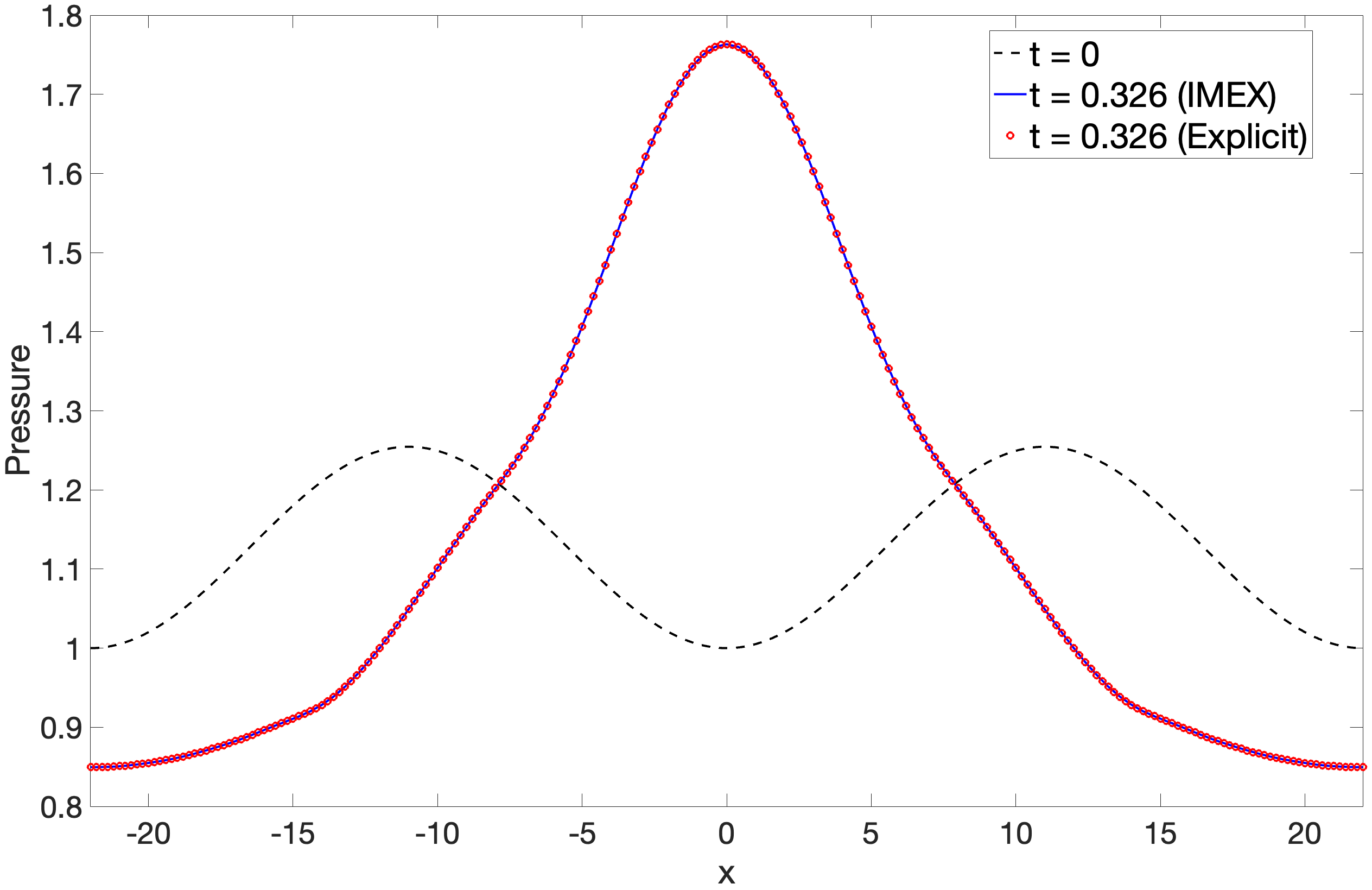}
    \end{subfigure}
    \begin{subfigure}{0.475\textwidth}
		\centering
        \includegraphics[width=0.95\textwidth]{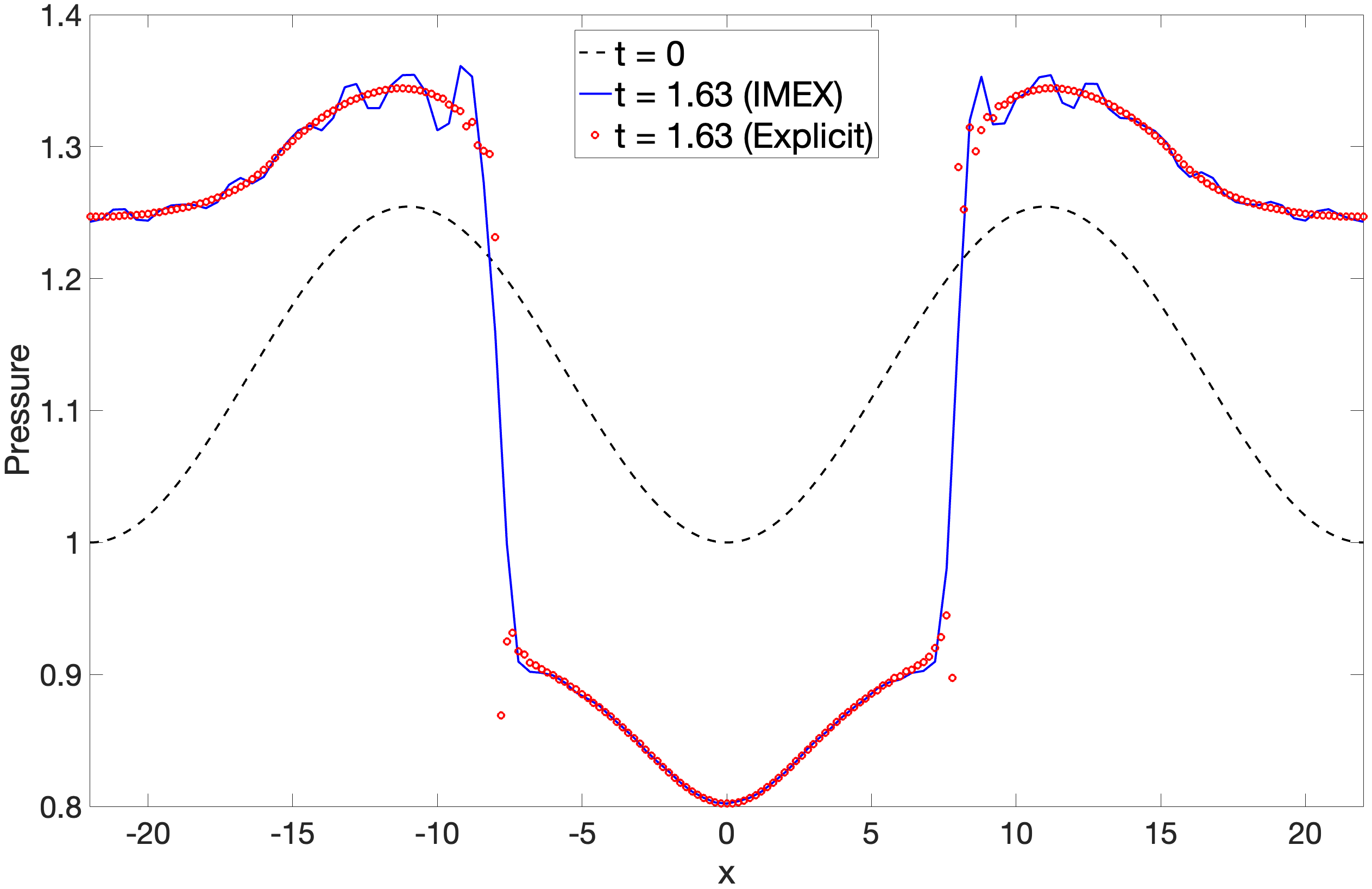}
   \end{subfigure}
    \caption{Colliding acoustic pulses test case employing the SG-EOS \eqref{eq:sg_eos} with $p_{\infty} = 6.8 \times 10^{-3}$, pressure profile. Left: $t = \frac{T_{f}}{5}$. Right: $t = T_{f}$. The initial profile is in dashed black line, the solid blue lines provide the results at the corresponding time obtained with the IMEX method, whereas the red dots show the reference results obtained with the explicit method.}
    \label{fig:colliding_pulses_pressure_SG}
\end{figure}

Finally, we employ the SG-EOS \eqref{eq:sg_eos}. We take $\gamma = 4.4, q_{\infty} = 0$, and we consider two different values of $\pi_{\infty}$, namely $\pi_{\infty} = 6.8 \times 10^{-3}$ and $\pi_{\infty} = 6.8 \times 10^{3}$. Notice that we do not modify the initial conditions, namely we keep $\tilde{u}_{0} = 2\sqrt{1.4}$ and $p_{0}^{'} = 2 \cdot 1.4 = 2.8$. First, we analyze the results with $\pi_{\infty} = 6.8 \times 10^{-3}$. A reference solution is computed using the third order SSP scheme with $\Delta t = 1.27344 \times 10^{-4}$ and $N_{el} = 880$ elements, leading to a maximum acoustic Courant number $C \approx 0.13$ and a maximum advective Courant number $C \approx 0.01$. The time step employed for the IMEX scheme is not modified, yielding a maximum acoustic Courant number $C \approx 1.09$ and a maximum advective Courant number $C \approx 0.09$. The pulses collide at $t = \frac{T_{f}}{5} = 0.326$ and an excellent agreement between the results obtained using the IMEX method and the reference ones is established (Figure \ref{fig:colliding_pulses_pressure_SG}). At $t = T_{f}$, shocks form around $x = \pm 7$ and spurious oscillations due to the high-order discretization methods arise.

\begin{figure}[h!]
   \begin{subfigure}{0.475\textwidth}   
        \centering
        \includegraphics[width=0.95\textwidth]{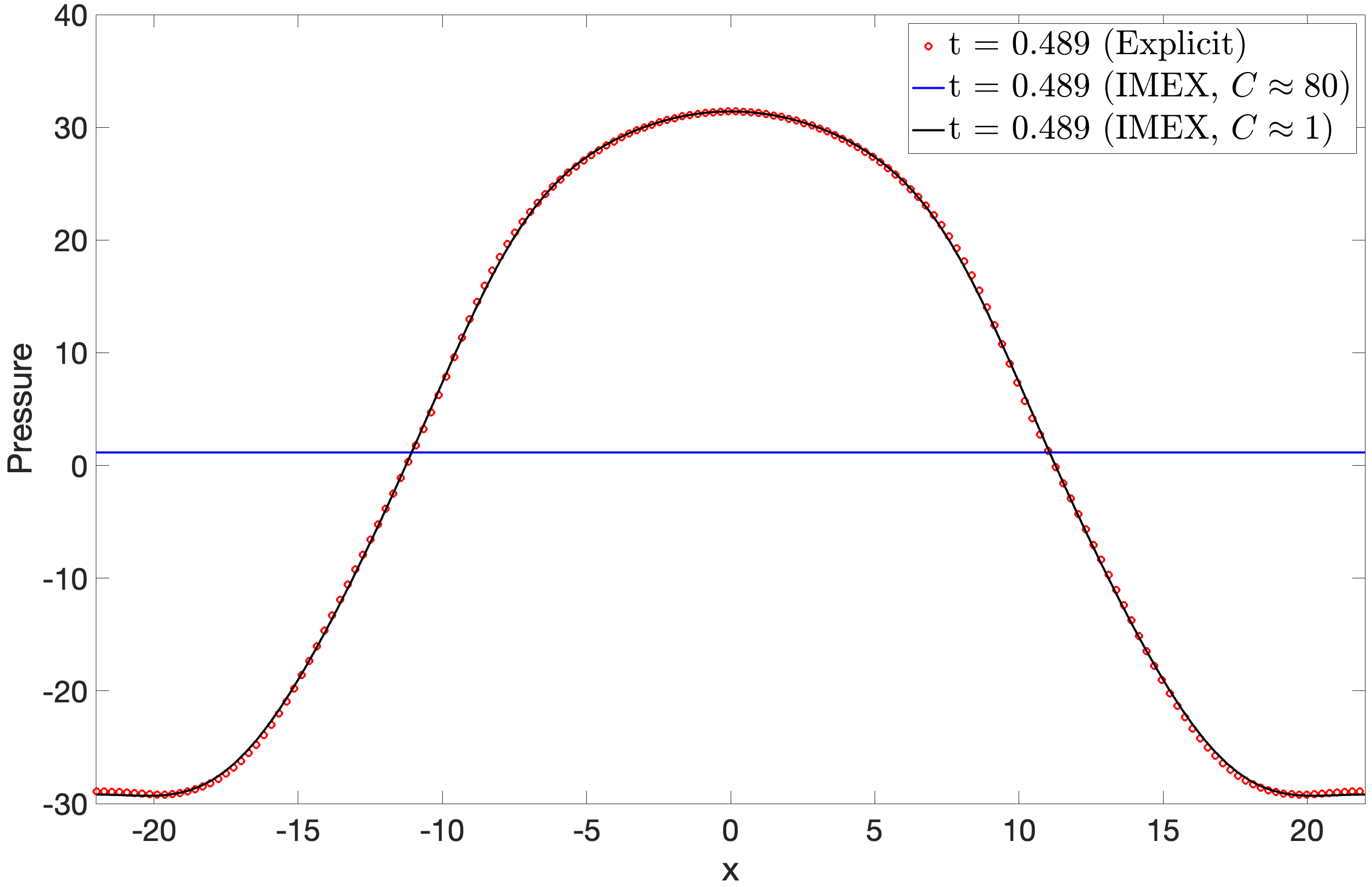}
    \end{subfigure}
    \begin{subfigure}{0.475\textwidth}
        \centering
        \includegraphics[width=0.95\textwidth]{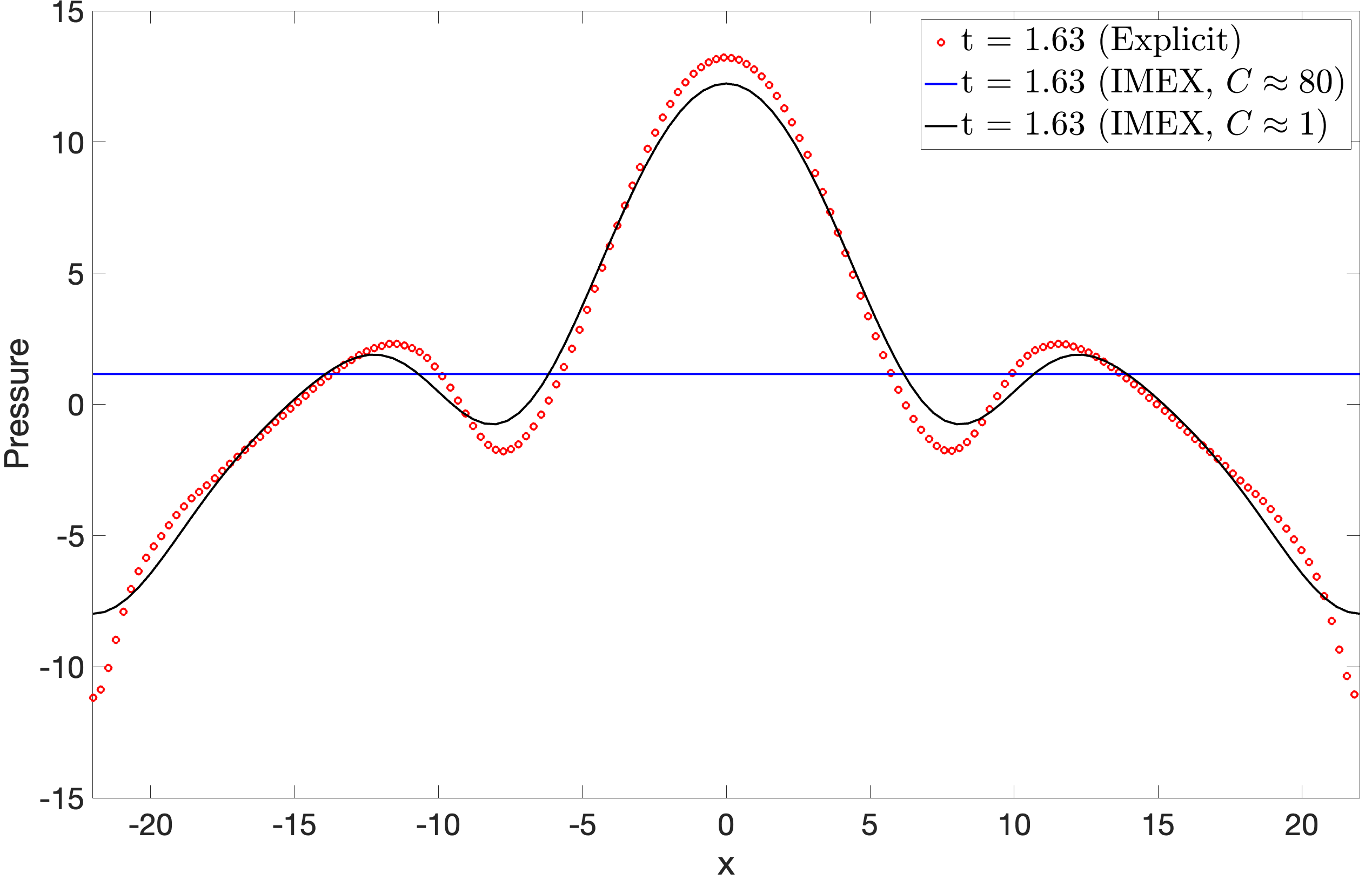}
   \end{subfigure}
    \caption{Colliding acoustic pulses test case employing the SG-EOS \eqref{eq:sg_eos} with $p_{\infty} = 6.8 \times 10^{3}$, pressure profile. Left: $t = \frac{3}{10}T_{f}$. Right: $t = T_{f}$. The solid blue lines provide the results at the corresponding time obtained with the IMEX method at acoustic Courant number $C \approx 80$, the solid black lines report the results obtained with the IMEX method at $C \approx 1$, whereas the red dots show the reference results obtained with the explicit method.}
    \label{fig:colliding_pulses_pressure_SG_low_Mach_acoustics}
\end{figure}

Next, we consider the configuration with $\pi_{\infty} = 6.8 \times 10^{3}$. A reference solution is computed using the third order SSP scheme with $\Delta t = 1.59179 \times 10^{-6}$ and $N_{el} = 880$ elements, yielding a maximum acoustic Courant number $C \approx 0.13$ and a maximum advective Courant number $C \approx 0.01$. For what concerns the IMEX method, a stable solution can be obtained without modifying the time step, but the pressure profiles are completely damped (Figure \ref{fig:colliding_pulses_pressure_SG_low_Mach_acoustics}). This is related to the fact the maximum acoustic Courant number is $C \approx 80$ and therefore we can no longer capture the acoustic wave. In order to correctly resolve the acoustic pulse, we have to decrease the value of the acoustic Courant number. We consider therefore $\Delta t = 2.0375 \times 10^{-4}$, namely a time step 80 times smaller than the previous one, so as to obtain $C \approx 1$. One can easily notice that a good agreement is established with the reference solution and the pressure profile is no longer damped, with the pulses colliding around $t = \frac{3}{10}T_{f}$ (Figure \ref{fig:colliding_pulses_pressure_SG_low_Mach_acoustics}). Moreover, no spurious oscillations arise. While the primary goal in the use of IMEX schemes is to obtain a method capable to resolve the material waves filtering out the acoustic waves, on the other hand, if a sufficiently small time step is employed, the method seems  to be naturally able to deal with low Mach acoustics. Notice that this is not valid in general for other low Mach fixes. In \cite{bruel:2019}, for example, the authors show that the low Mach fixes proposed in \cite{dellacherie:2010a} and \cite{rieper:2011} can suffer of spurious oscillations and of order reduction when applied to low Mach acoustics. A correction able to deal with low Mach acoustic was developed in \cite{bruel:2019} and, more recently, in \cite{galie:2024}. A more detailed analysis of the low Mach acoustic behaviour will be the focus of future work.

\subsection{Density layering}
\label{ssec:density_layering}

We consider now the test case II proposed in \cite{klein:1995} and described also, e.g., in \cite{noelle:2014}. The domain is $\Omega = \left(-L, L\right)$, with $L = \frac{1}{0.02} = 50$. The initial conditions read as follows: \\
\begin{subequations}
\begin{eqnarray}
    \rho\left(x, 0\right) &=& \bar{\rho}_{0} + \Phi\left(x\right)\tilde{\rho}_{0}\sin\left(\frac{40\pi x}{L}\right) + \frac{1}{2}M\rho_{1}\left(1 + \cos\left(\frac{\pi x}{L}\right)\right) \label{eq:density_layering_rho_init} \\
    u\left(x, 0\right) &=& \frac{1}{2}\tilde{u}_{0}\left(1 + \cos\left(\frac{\pi x}{L}\right)\right) \label{eq:density_layering_u_init} \\
    p\left(x,0\right) &=& \bar{p}_{0} + \frac{1}{2}Mp_{1}\left(1 + \cos\left(\frac{\pi x}{L}\right)\right), \label{eq:density_layering_p_init}
\end{eqnarray}
\end{subequations}
with $\bar{\rho}_{0} = 1, \tilde{\rho}_{0} = \frac{1}{2}, \rho_{1} = 2, \tilde{u}_{0} = 2\sqrt{\gamma} = 2\sqrt{1.4}, \bar{p}_{0} = 1$, and $p_{1} = 2\gamma = 2.8$. Finally, the function $\Phi\left(x\right)$ is given by
\begin{equation}
    \Phi\left(x\right) =
    \begin{cases}
        \frac{1}{2}\left(1 - \cos\left(\frac{5\pi x}{L}\right)\right) \qquad &\text{if } 0 \le x \le \frac{2}{5}L \\
        0 \qquad &\text{otherwise}.
    \end{cases}
\end{equation}
The initial data describe a density layering of large amplitude and small wavelengths, which is driven by the motion of a right-moving periodic acoustic wave with long wavelength. Periodic boundary conditions are prescribed. The final time is $T_{f} = 5.071$. We consider a computational grid composed by $N_{el} = 250$ elements with $r = 2$, whereas the time step is $\Delta t = 1.6903 \times 10^{-2}$. Following \cite{noelle:2014}, we start considering $M = \frac{1}{50}$. Hence, the advective Courant $C_{u}$ is around $0.2$, while the acoustic Courant number $C$ is around $7$. A comparison between the initial and the final time for both the density and the pressure displays a good agreement with the results presented in \cite{klein:1995, noelle:2014} (Figure \ref{fig:density_layering_IG_M_0,02}). One can easily notice that the acoustic wave transports the density layer and the shape of the layer is undistorted. As in the previous test case, due to weakly non-linear effects, the pulse starts steepening, leading to shock formation. A reference solution has been computed using the explicit third order SSP scheme. The time step employed for the explicit scheme is $\Delta t = 5.071 \times 10^{-4}$, namely around 33 times smaller than that used with the IMEX scheme. An excellent agreement is established between the two solutions. 

\begin{figure}[h!]
    \centering
    \begin{subfigure}{0.475\textwidth}
        \centering
        \includegraphics[width = 0.95\textwidth]{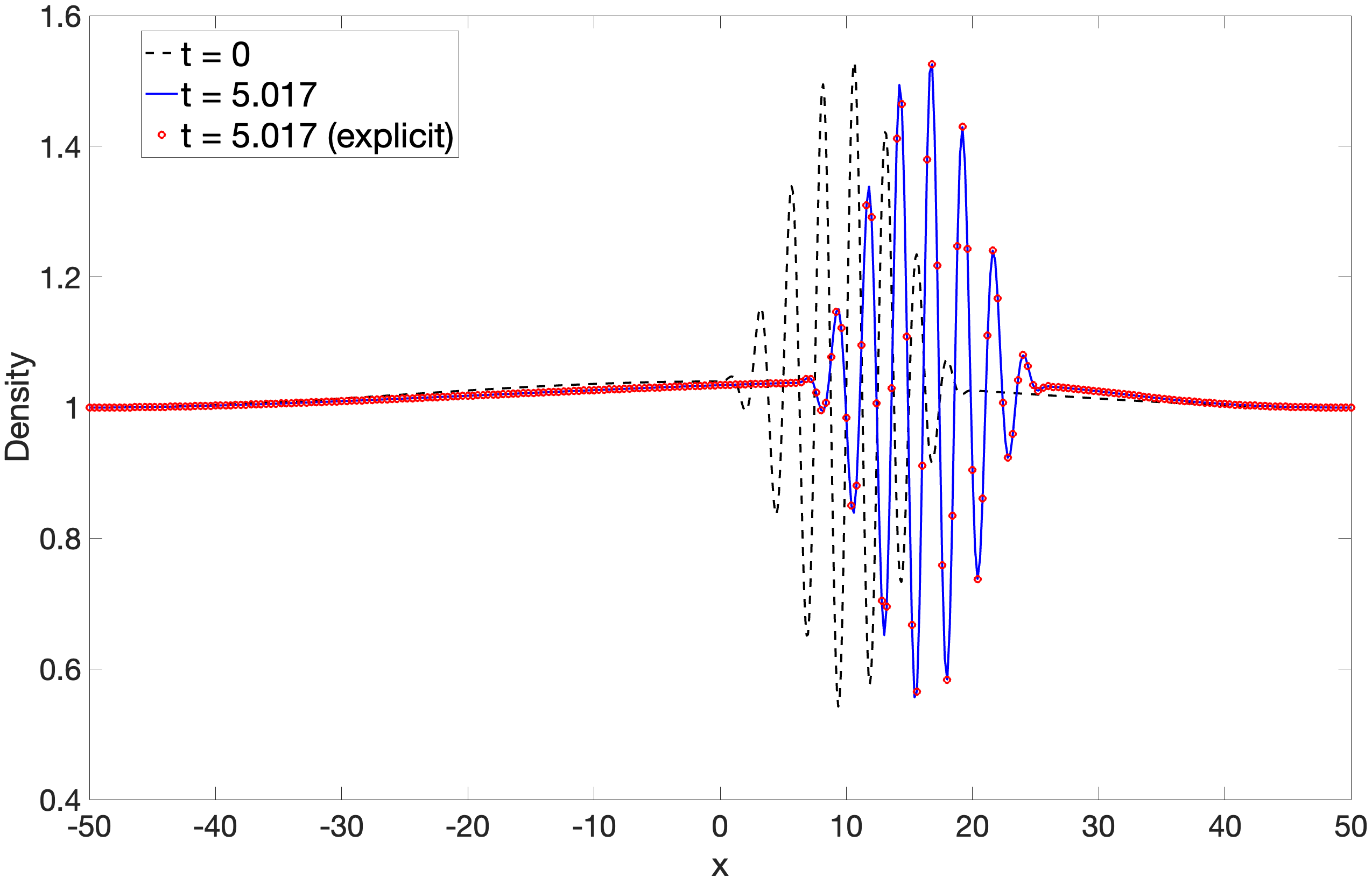}
    \end{subfigure}
    \begin{subfigure}{0.475\textwidth}
	\centering
        \includegraphics[width = 0.95\textwidth]{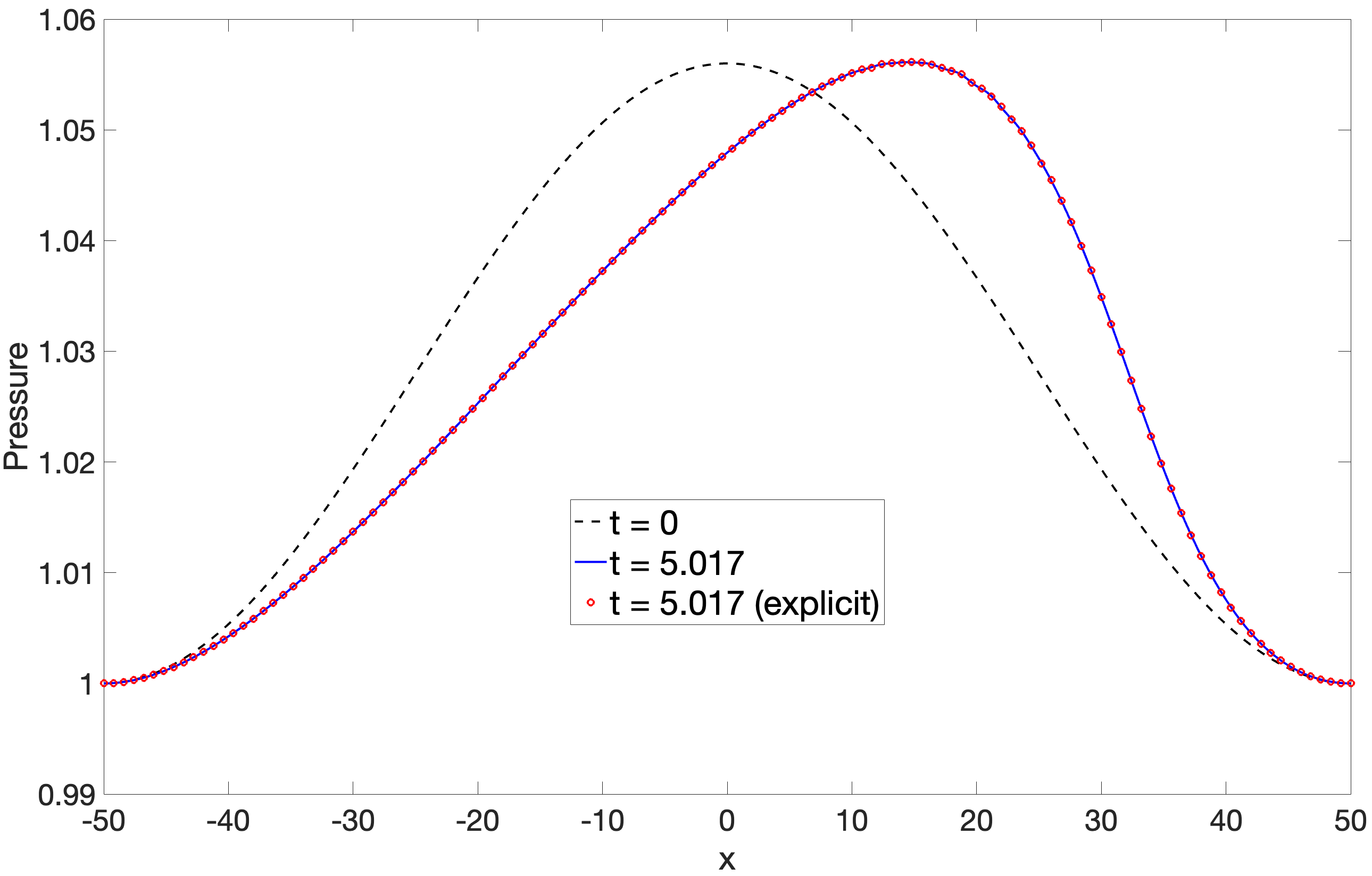}
    \end{subfigure}
    \caption{Density layering test case at $M = 0.02$ with the ideal gas law \eqref{eq:ideal_gas}. Left: density. Right: pressure. The dashed black lines represent the initial condition, the continuous blue lines show the solution at the final time, whereas the red dots report the solution obtained with the third order optimal explicit SSP scheme.}
    \label{fig:density_layering_IG_M_0,02}
\end{figure}

For the sake of completeness, we also consider a case even closer to the incompressible regime, taking $M = 10^{-4}$, which results in an acoustic Courant number $C \approx 1400$. The analytical solution of the leading order term of the limit model \eqref{eq:euler_adim_ap_two_scale} with initial conditions \eqref{eq:density_layering_rho_init}-\eqref{eq:density_layering_p_init} is
\begin{equation}
    \bar{\rho} = \bar{\rho}_{0} + \Phi\left(x - \bar{u}\left(t\right)t\right)\tilde{\rho}_{0}\sin\left(\frac{40\pi\left(x - \bar{u}\left(t\right)t\right)}{L}\right) \qquad \bar{u} = \bar{u}\left(t\right) \qquad \bar{p} = \bar{p}_{0},
\end{equation}
$\bar{u}\left(t\right)$ being a function only of time. Since we are considering periodic boundary conditions, the integral over the domain of $\rho u$ is constant and therefore the steady value of $\bar{u}$ is
\begin{equation}
    \bar{u} = \frac{\int_{\Omega}\left(\bar{\rho}_{0} + \Phi\left(x\right)\tilde{\rho}_{0}\sin\left(\frac{40\pi\left(x\right)}{L}\right)\right)\left(\frac{1}{2}\tilde{u}_{0}\left(1 + \cos\left(\frac{\pi x}{L}\right)\right)\right)d\Omega}{\int_{\Omega}\left(\frac{1}{2}\tilde{u}_{0}\left(1 + \cos\left(\frac{\pi x}{L}\right)\right)\right)d\Omega} = \frac{\tilde{u}_{0}}{2} - \frac{164375 - 32875\sqrt{5}}{1320441408\pi}\frac{\tilde{\rho}_{0}\tilde{u}_{0}}{\bar{\rho}_{0}} \approx \frac{\tilde{u}_{0}}{2} = \sqrt{\gamma}.
\end{equation}
A comparison at $t = T_{f}$ between the analytical solution as $M \to 0$ and the numerical results shows an excellent agreement for both the density and the pressure profile (Figure \ref{fig:density_layering_IG_M_0,0001}). Notice that the initial velocity field is not divergence-free, namely it is not well-prepared. However, the numerical method leads to the incompressible limit, as already discussed in Section \ref{ssec:ap_single_length} (Figure \ref{fig:density_layering_IG_M_0,0001}). For further reference, we include the solution obtained employing the explicit scheme with $\Delta t = 2.5355\times 10^{-6}$, i.e. a time step around $6666$ times smaller. While on the one hand, the use of high-order discretization schemes reduces the numerical diffusion and allows for preserving the shape of the layer also employing the explicit method, on the other hand, the incompressible limit is not achieved (Figure \ref{fig:density_layering_IG_M_0,0001}). This result confirms the necessity to employ asymptotic-preserving methods as $M \to 0$ to obtain reliable results as well as to be much more efficient.

\begin{figure}[h!]
    \centering
    \begin{subfigure}{0.475\textwidth}
        \centering
        \includegraphics[width = 0.95\textwidth]{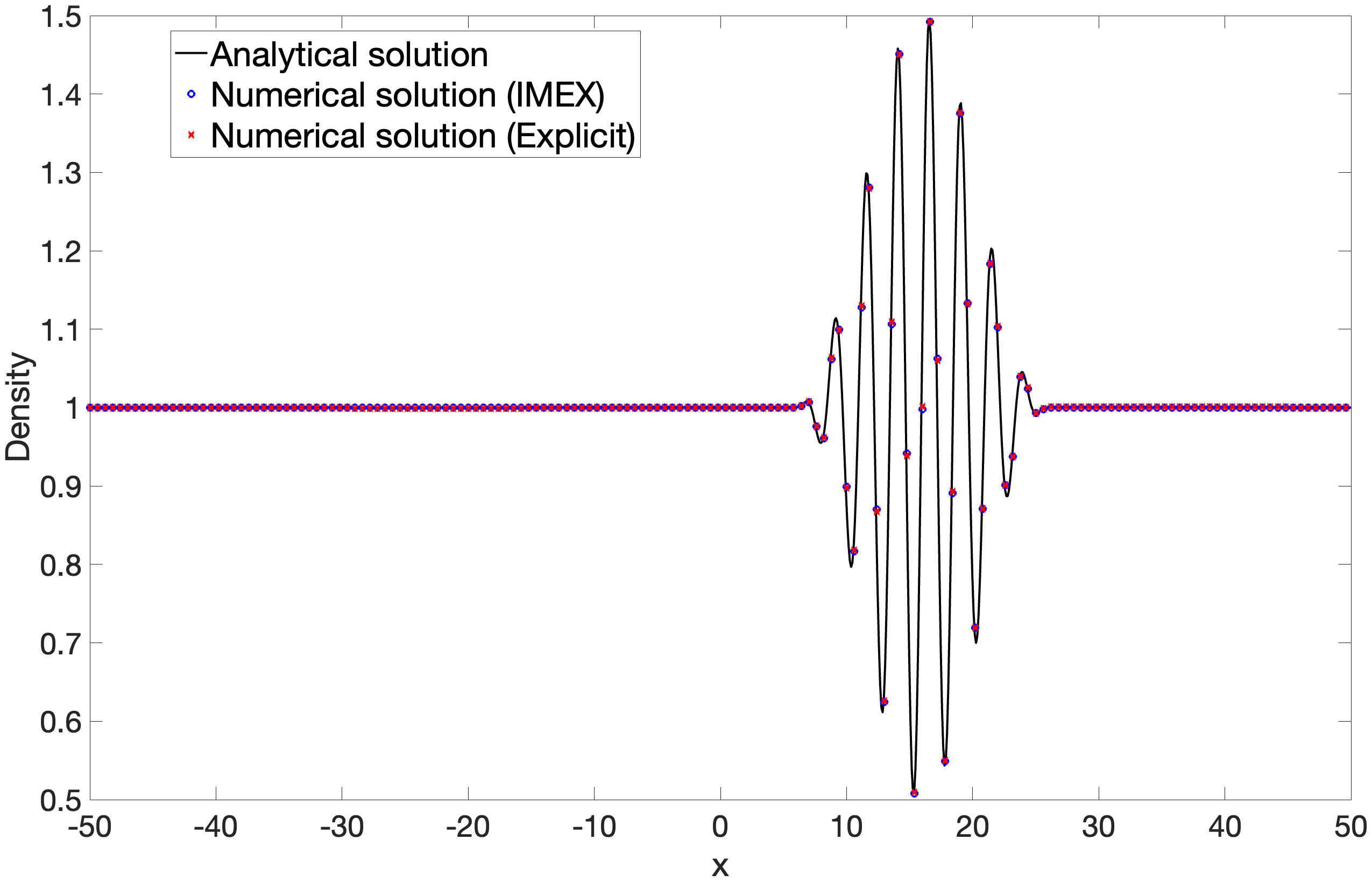} a)
    \end{subfigure}
    \begin{subfigure}{0.475\textwidth}
        \centering
        \includegraphics[width = 0.95\textwidth]{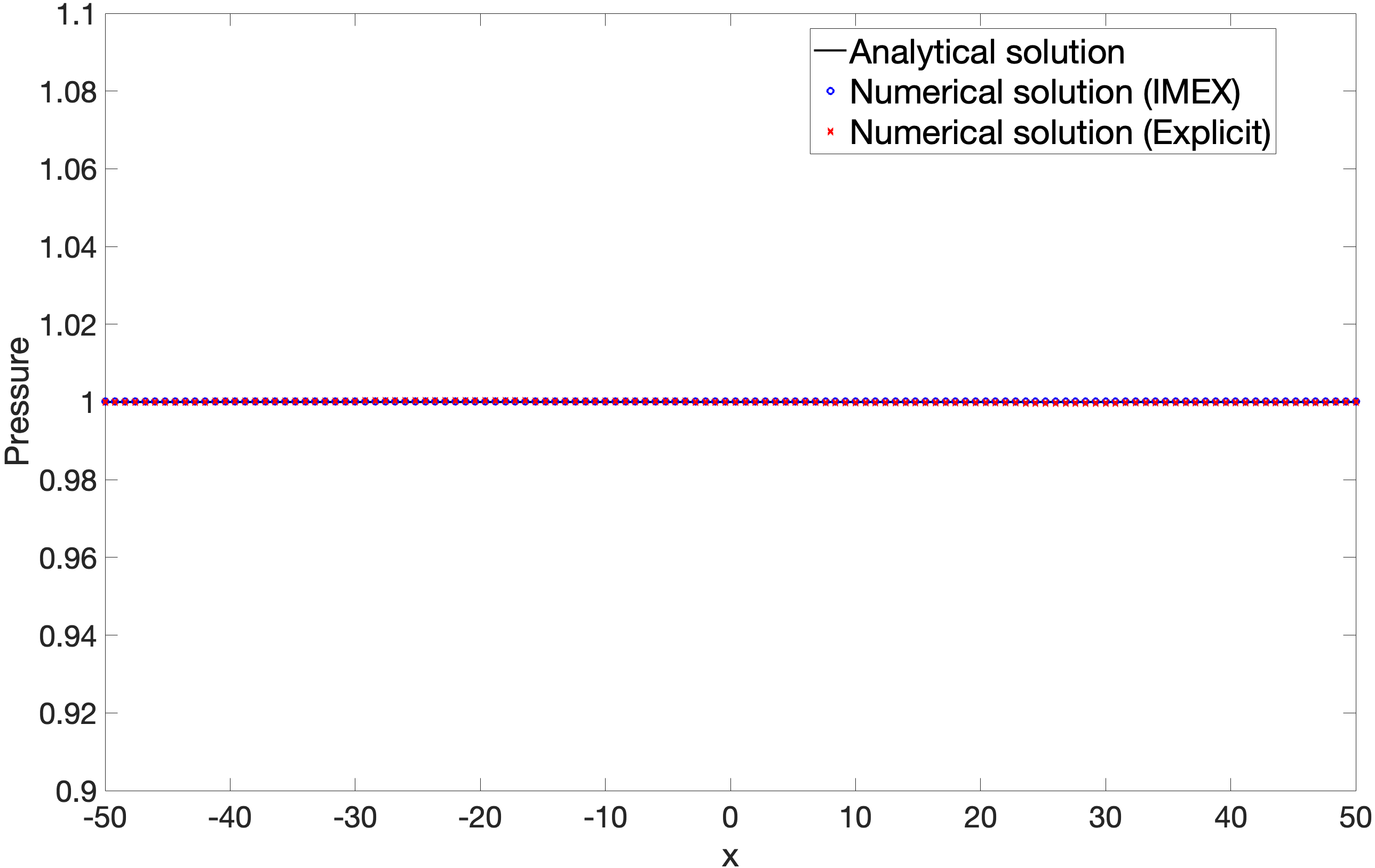} b)
    \end{subfigure}
    \begin{subfigure}{0.475\textwidth}
        \centering
        \includegraphics[width = 0.95\textwidth]{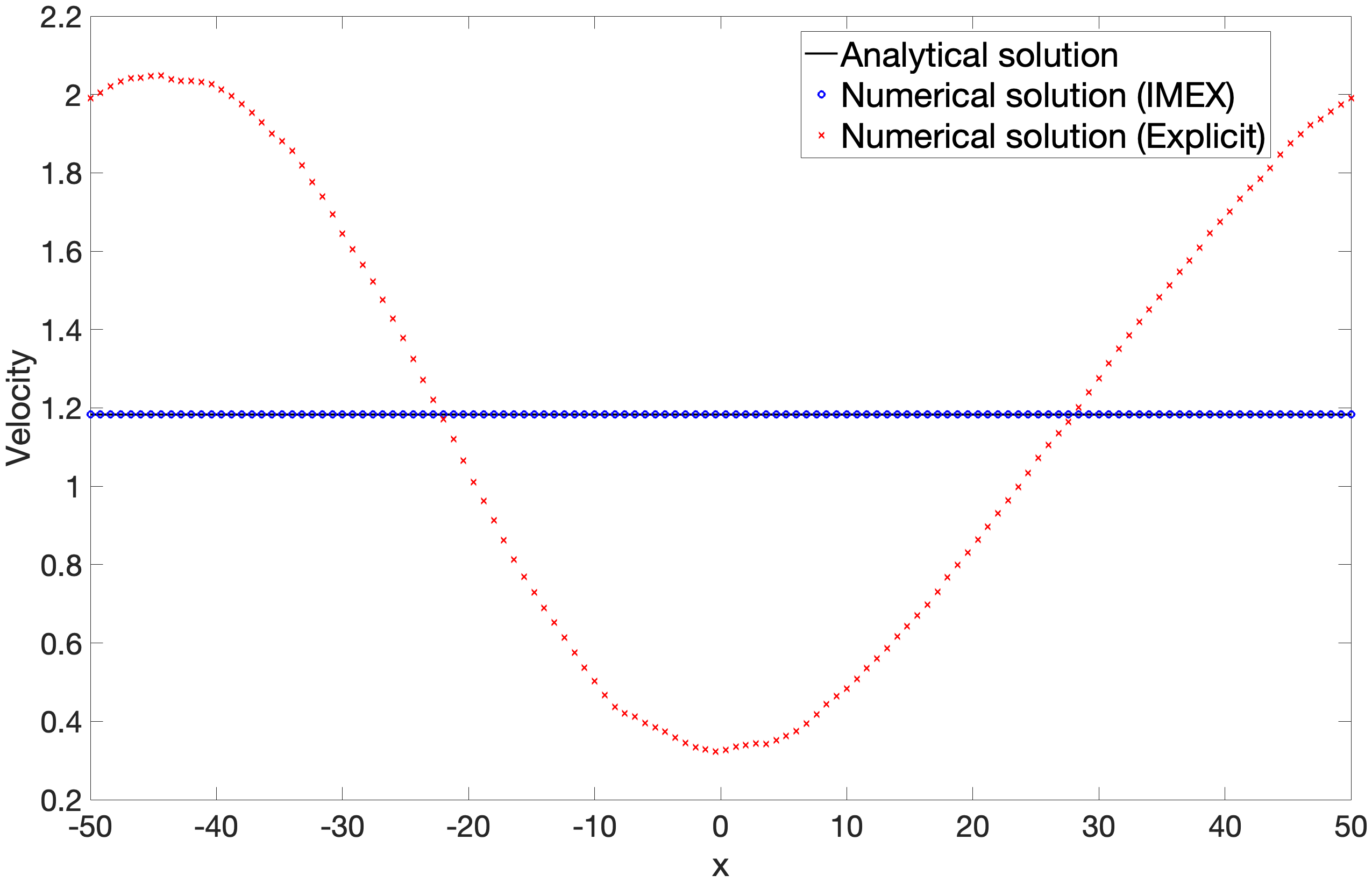} c)
    \end{subfigure}
    \caption{Density layering test case at $M = 10^{-4}$ with the ideal gas law \eqref{eq:ideal_gas}. a) density, b) pressure, c) velocity. The continuous black lines represent the analytical solution of the limit model \eqref{eq:euler_adim_ap_two_scale}, the blue dots report the numerical results obtained with the IMEX method, whereas the red crosses show the results achieved with the fully explicit scheme.}
    \label{fig:density_layering_IG_M_0,0001}
\end{figure}

We now consider a configuration of this test case for the SG-EOS \eqref{eq:sg_eos}. We take $\gamma = 4.4, \pi_{\infty} = 6.8 \times 10^{-3}$, and $q_{\infty} = 0$ in \eqref{eq:sg_eos}. Notice that, we do not modify the initial conditions \eqref{eq:density_layering_rho_init}-\eqref{eq:density_layering_p_init}, namely we keep $\tilde{u}_{0} = 2\sqrt{1.4}$ and $p_{1} = 2.8$. We start considering $M = \frac{1}{50} = 0.02$, which yields an acoustic Courant number $C \approx 12.4$ and a maximum advective Courant number $C_{u} \approx 0.2$. Figure \ref{fig:density_layering_SG_M_0,02} shows a comparison between the initial and the final time for both the density and the pressure. A reference solution has been computed using the third order explicit SSP scheme, with a time step $\Delta t = 2.5355 \times 10^{-4}$, namely a time step around $66$ times smaller than that employed with the IMEX scheme. One can easily notice that the density layer is transported without too much damping. Moreover, an excellent agreement with the explicit solution is established. Finally, for what concerns the incompressible limit at $M = 10^{-4}$, since $\frac{\partial\bar{p}_{0}}{\partial t} = 0$, all the equations of state lead to the same limit (see \eqref{eq:dive_u}). This is further confirmed by the density and pressure profiles reported in Figure \ref{fig:density_layering_SG_M_0,0001}.

\begin{figure}[h!]
    \centering
    \begin{subfigure}{0.475\textwidth}
	\centering
        \includegraphics[width = 0.95\textwidth]{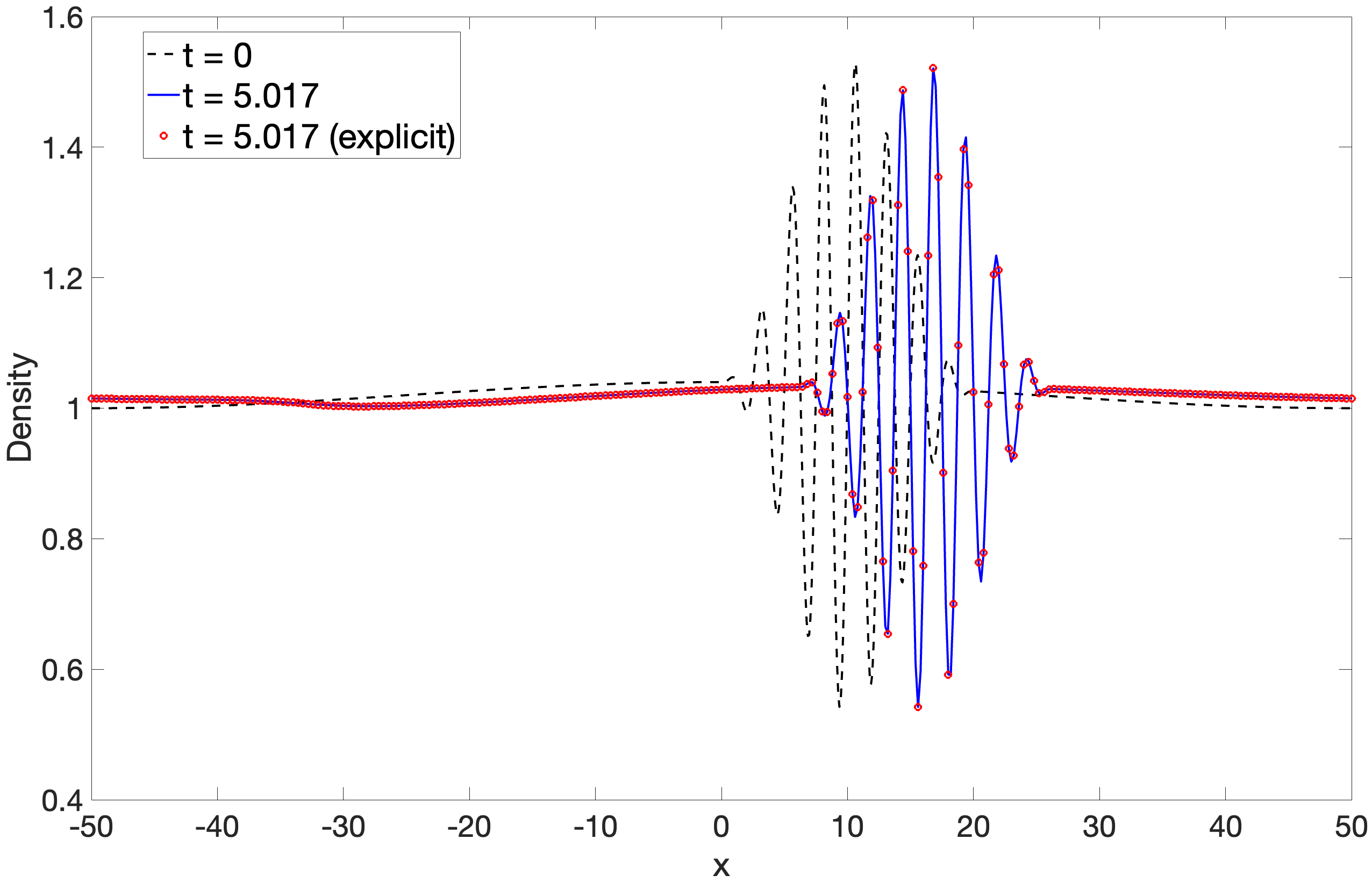}
    \end{subfigure}
    \begin{subfigure}{0.475\textwidth}
	\centering
        \includegraphics[width = 0.95\textwidth]{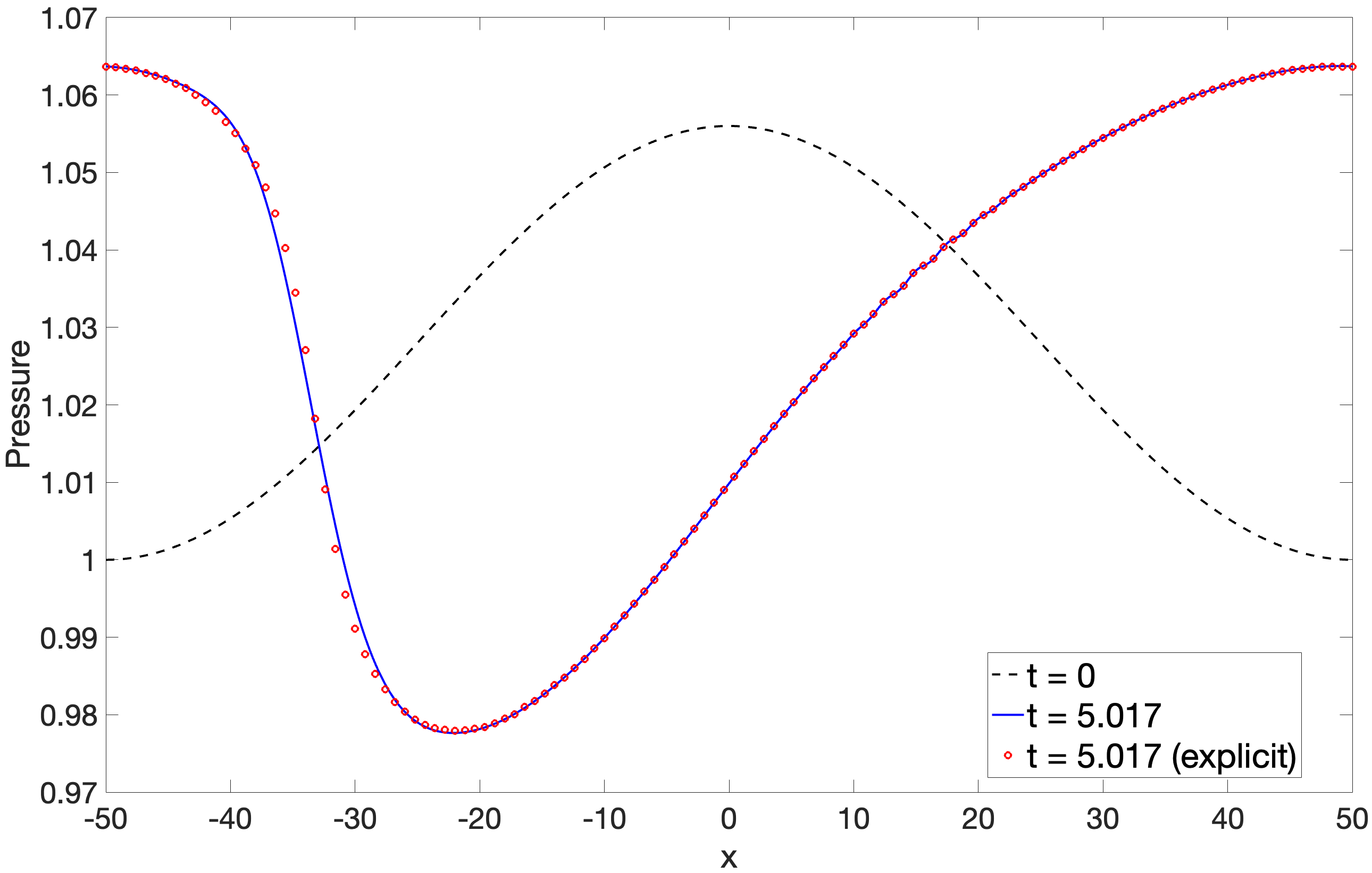}
    \end{subfigure}
    \caption{Density layering test case at $M = 0.02$ with the SG-EOS \eqref{eq:sg_eos}. Left: density. Right: pressure. The dashed black lines represent the initial condition, the continuous blue lines show the solution at the final time, whereas the red dots report the solution obtained with the third order optimal explicit SSP scheme.}
    \label{fig:density_layering_SG_M_0,02}
\end{figure}

\begin{figure}[h!]
    \centering
    \begin{subfigure}{0.475\textwidth}
	\centering
        \includegraphics[width = 0.9\textwidth]{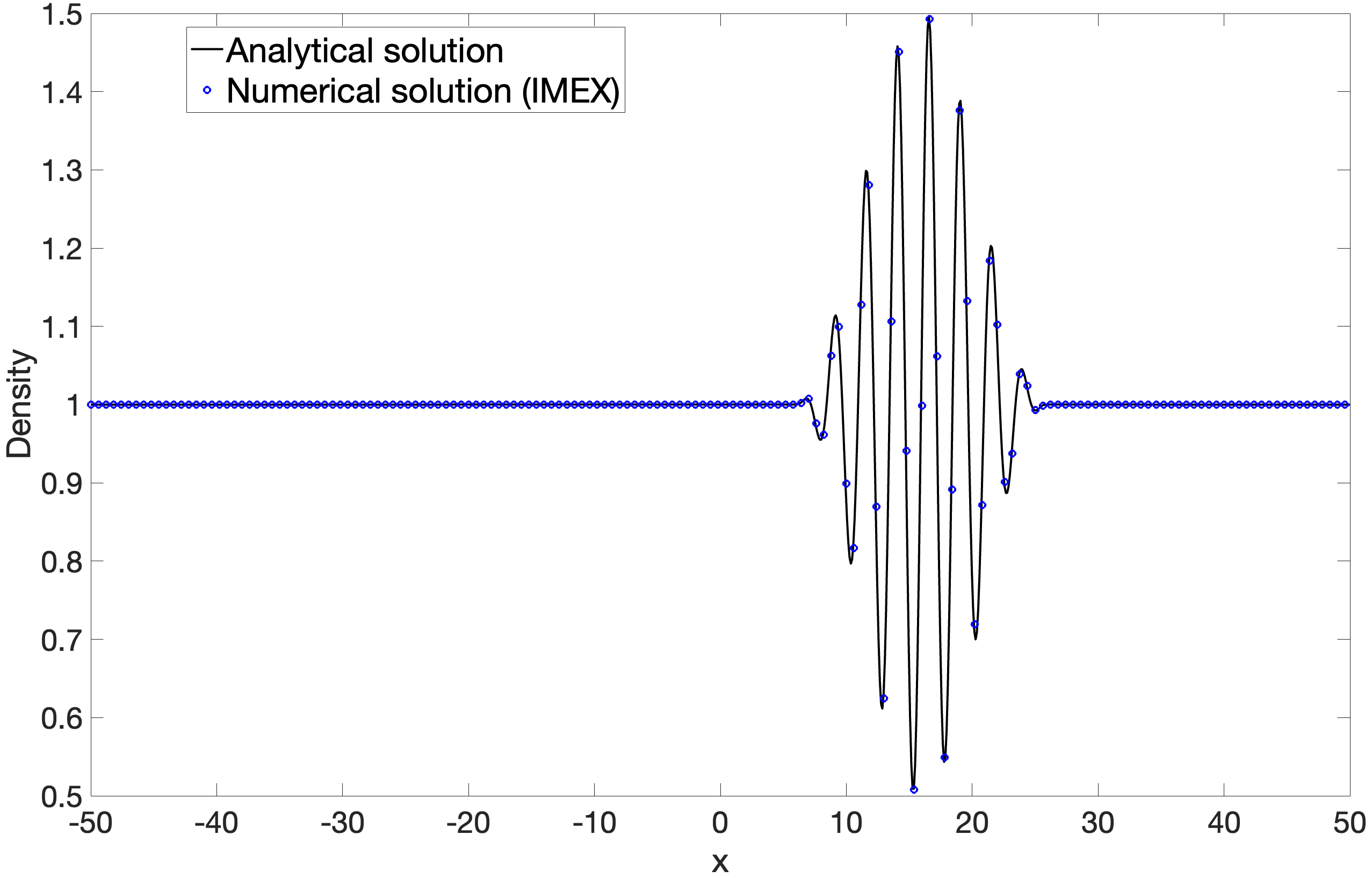}
    \end{subfigure}
    \begin{subfigure}{0.475\textwidth}
	\centering
        \includegraphics[width = 0.9\textwidth]{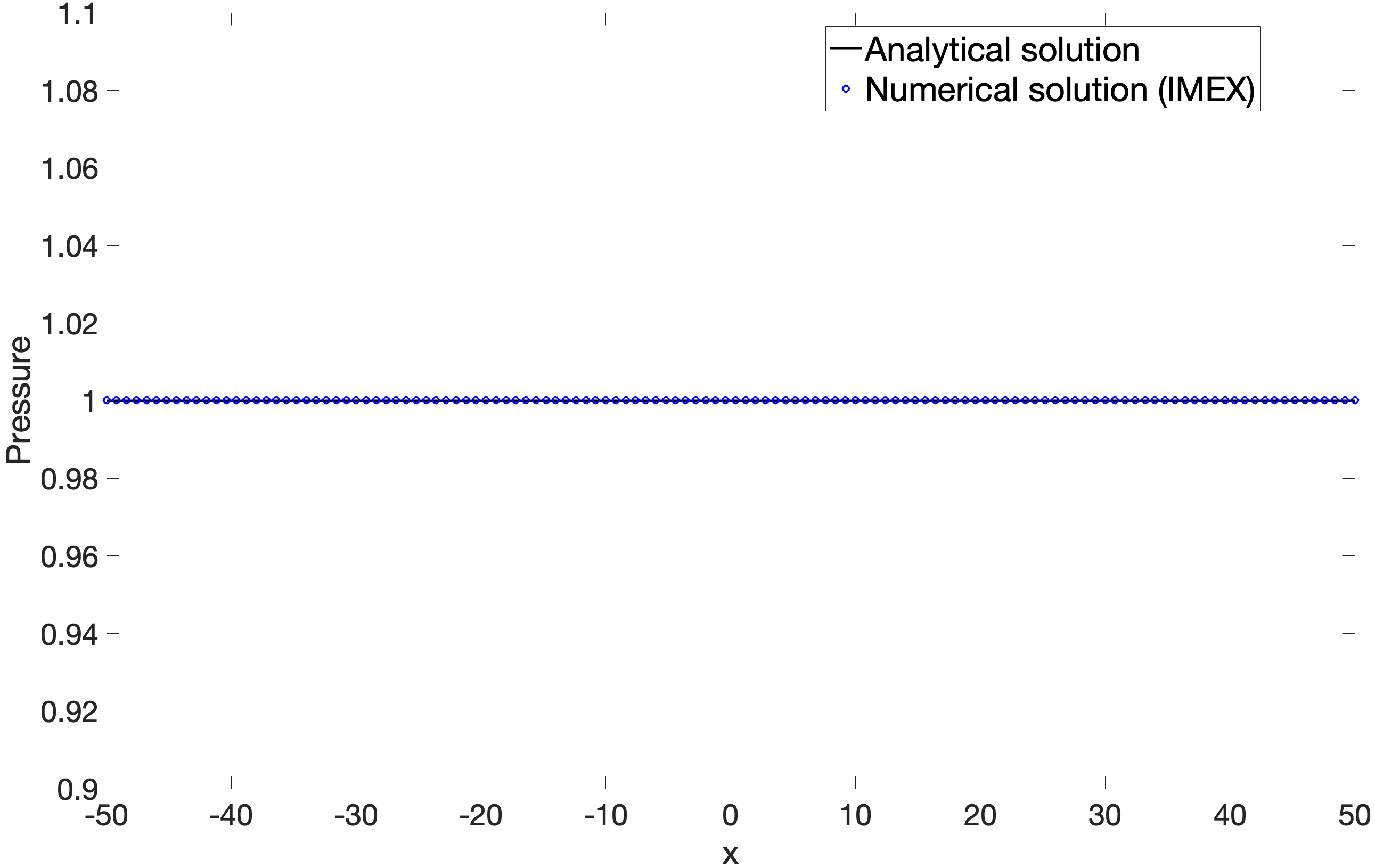}
    \end{subfigure}
    \caption{Density layering test case at $M = 10^{-4}$ with the SG-EOS \eqref{eq:sg_eos}. Left: density. Right: pressure. The continuous black lines represent the analytical solution of the limit model \eqref{eq:euler_adim_ap_two_scale}, while the blue dots report the numerical results.}
    \label{fig:density_layering_SG_M_0,0001}
\end{figure}

\subsection{Flow in an open tube}
\label{ssec:open_tube}

\begin{figure}[h!]
    \centering
    \begin{subfigure}{0.475\textwidth}
	\centering
        \includegraphics[width = 0.9\textwidth]{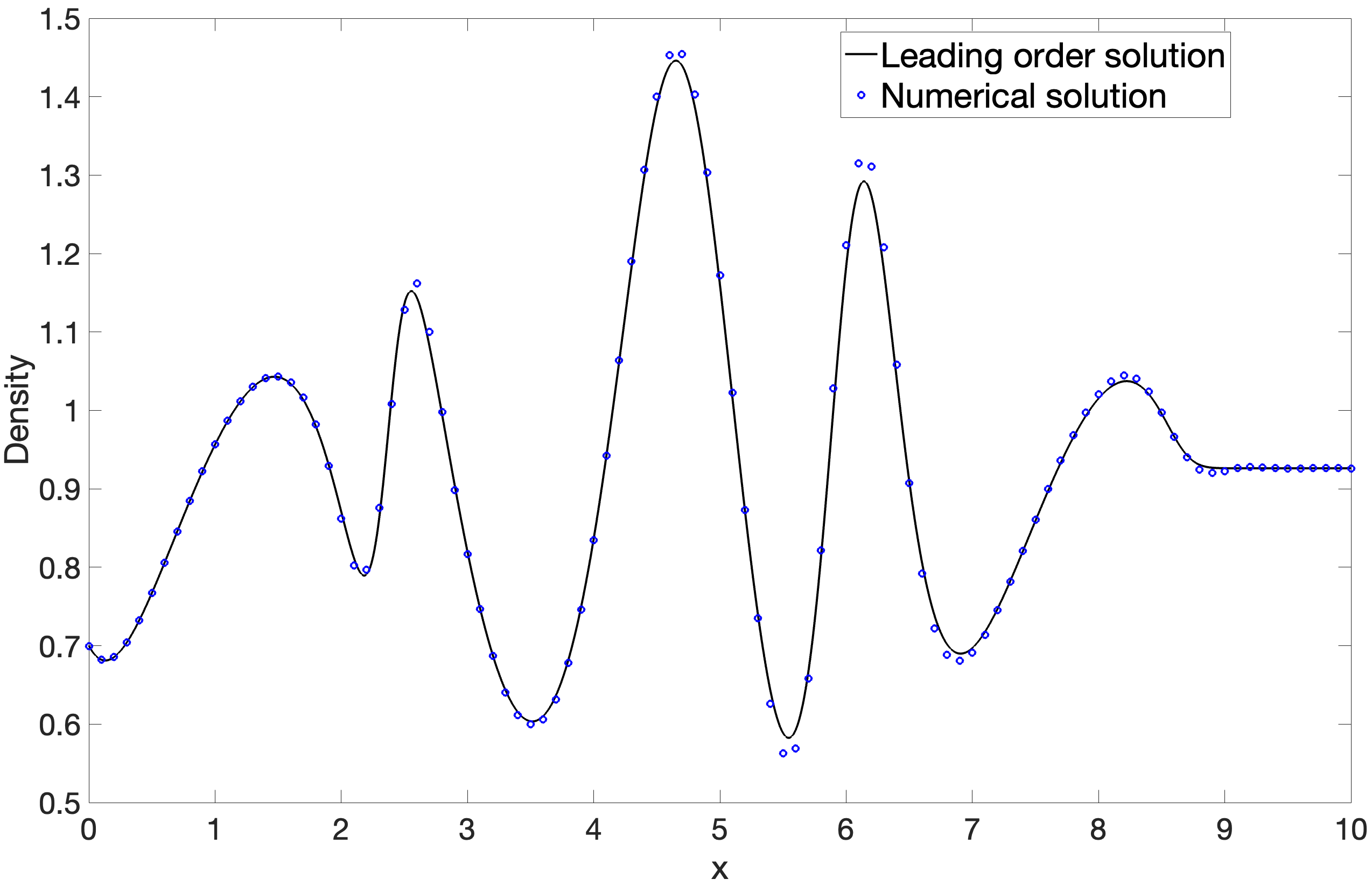}
    \end{subfigure}
    \begin{subfigure}{0.475\textwidth}
	\centering
        \includegraphics[width = 0.9\textwidth]{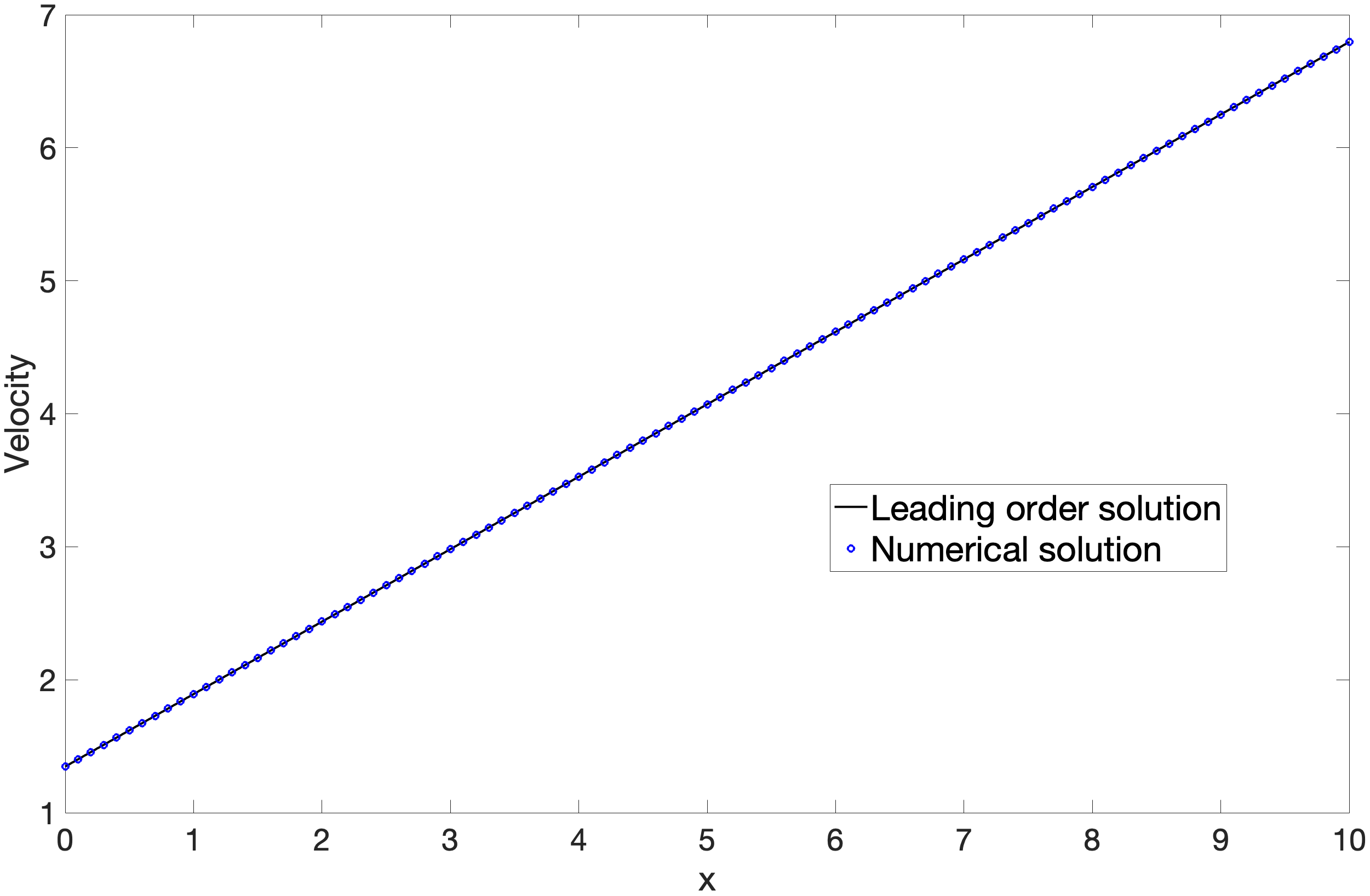}
    \end{subfigure}
    \caption{Open tube test case with the ideal gas law \eqref{eq:ideal_gas}, results at $t = T_{f} = 7.47$. Left: density. Right: velocity. The continuous black line shows the leading order solution as $M \to 0$, whereas the blue dots report the numerical results.}
    \label{fig:open_tube_IG} 
\end{figure}

We consider now the test case III proposed in \cite{klein:1995} for an ideal gas, which we recall here for the convenience of the reader. A flow in an open tube represented by the domain $\Omega = \left(0, 10\right)$ is analyzed; at the left-end a time dependent density and velocity are prescribed, whereas at the right-end a time dependent outflow pressure with large amplitude variation is imposed. This kind of boundary conditions is employed e.g. in the case of subsonic inflow and subsonic outflow \cite{carlson:2011}. More specifically, the initial conditions read as follows:
\begin{equation}
    \left(\rho, u, p\right)\left(x, 0\right) = \left(1, 1, 1\right),
\end{equation}
while the boundary conditions are
\begin{equation}
    \rho\left(0, t\right) = 1 + \frac{3}{10}\sin\left(4t\right) \qquad
    u\left(0, t\right) = 1 + \frac{1}{2}\sin\left(2t\right) \qquad
    p\left(L, t\right) = 1 + \frac{1}{4}\sin\left(3t\right), 
\end{equation}
with $L = 10$. The final time is $T_{f} = 7.47$. The Mach number is set to $M = 10^{-4}$. We consider a number of elements $N_{el} = 50$ with $r = 2$, whereas the time step is $\Delta t = 9.3375 \times 10^{-4}$, leading to a maximum advective Courant number $C_{u} \approx 0.07$ and a maximum acoustic Courant number $C \approx 155$. The results at $t = T_{f}$ are those expected by the asymptotic analysis for both the density and velocity profile (Figure \ref{fig:open_tube_IG}). The limit solution as $M \to 0$ has been included for comparison in Figure \ref{fig:open_tube_IG}. For an ideal gas, \eqref{eq:energy_limit_incomp} reduces to
\begin{equation}
    \dive\bar{\mathbf{u}} = -\frac{1}{\gamma\bar{p}}\frac{\partial\bar{p}}{\partial t} = -\frac{1}{\gamma\bar{p}}\frac{d\bar{p}}{dt},
\end{equation}
since $\grad\bar{p} = \mathbf{0}$. Hence, for $M \to 0$, in one space dimension, $\frac{\partial\bar{\mathbf{u}}}{\partial x}$ is a function only of time and, therefore, the velocity is a linear function of space with a given time dependent slope and boundary value at $x = 0$. For what concerns the density, we rewrite \eqref{eq:continuity_limit} as follows:
\begin{equation}
    \frac{\partial\bar{\rho}}{\partial t} + \bar{\mathbf{u}} \cdot \grad\bar{\rho} + \bar{\rho}\dive\bar{\mathbf{u}} = \frac{\partial\bar{\rho}}{\partial t} + \bar{\mathbf{u}} \cdot \grad\bar{\rho} - \frac{1}{\gamma}\frac{\bar{\rho}}{\bar{p}}\frac{d\bar{p}}{dt} = 0,
\end{equation} 
or, equivalently,
\begin{equation}
    \frac{D\log\bar{\rho}}{Dt} = \frac{1}{\gamma}\frac{d\log\bar{p}}{dt},
\end{equation}
with $\frac{D}{Dt} = \frac{\partial}{\partial t} + \mathbf{u} \cdot \grad$ denoting the Lagrangian derivative. Hence, as discussed in \cite{klein:1995}, the material elements undergo a quasi-static adiabatic compression and expansion following the particle paths described by $\bar{u}$. One can easily notice from the density profile in Figure \ref{fig:open_tube_IG} that mass elements, after entering the domain at the left-end, are correctly compressed and expanded.

\begin{figure}[h!]
    \centering
    \begin{subfigure}{0.475\textwidth}
	\centering
	\includegraphics[width = 0.9\textwidth]{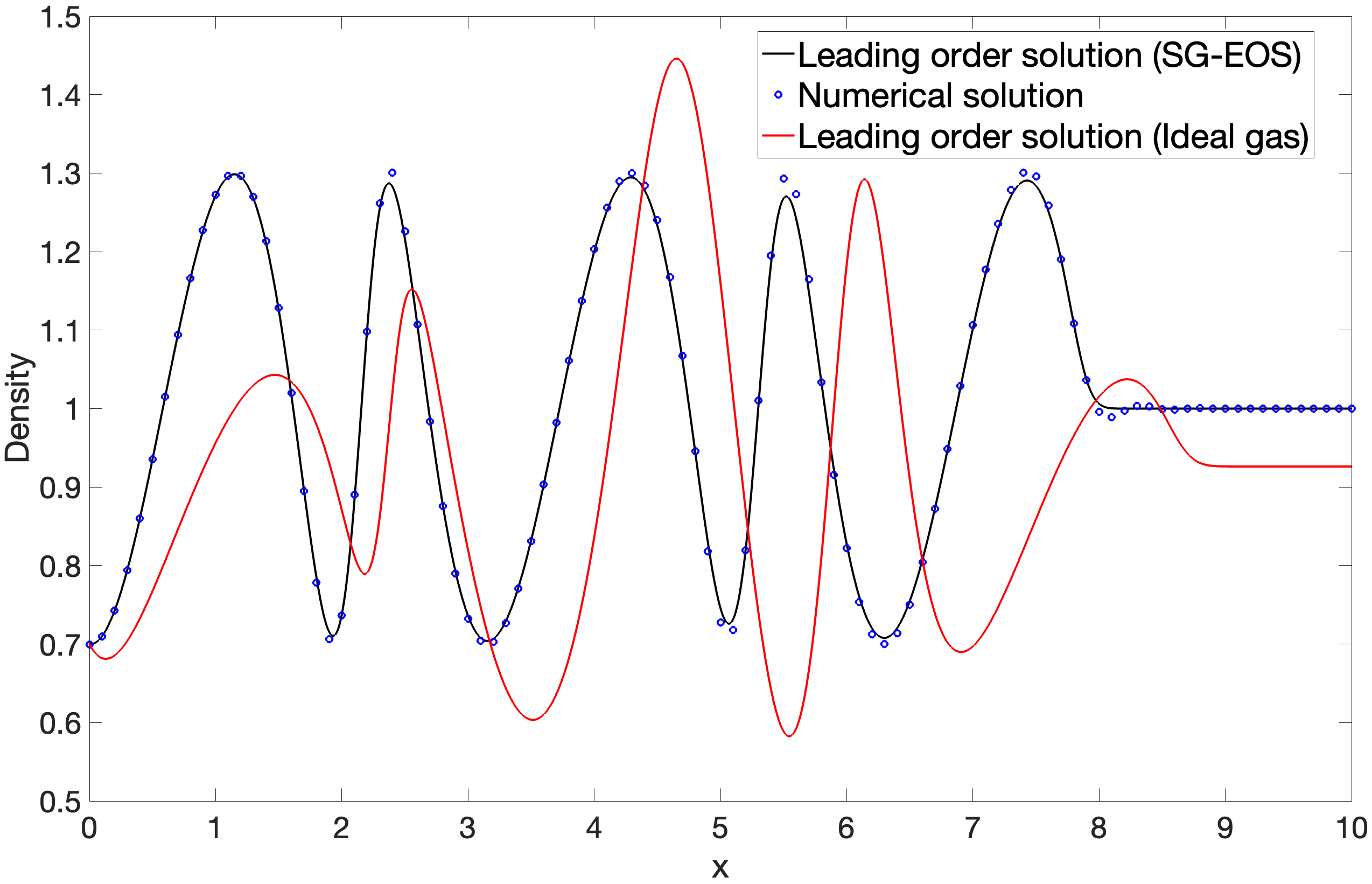}
    \end{subfigure}
    \begin{subfigure}{0.475\textwidth}
	\centering
	\includegraphics[width = 0.9\textwidth]{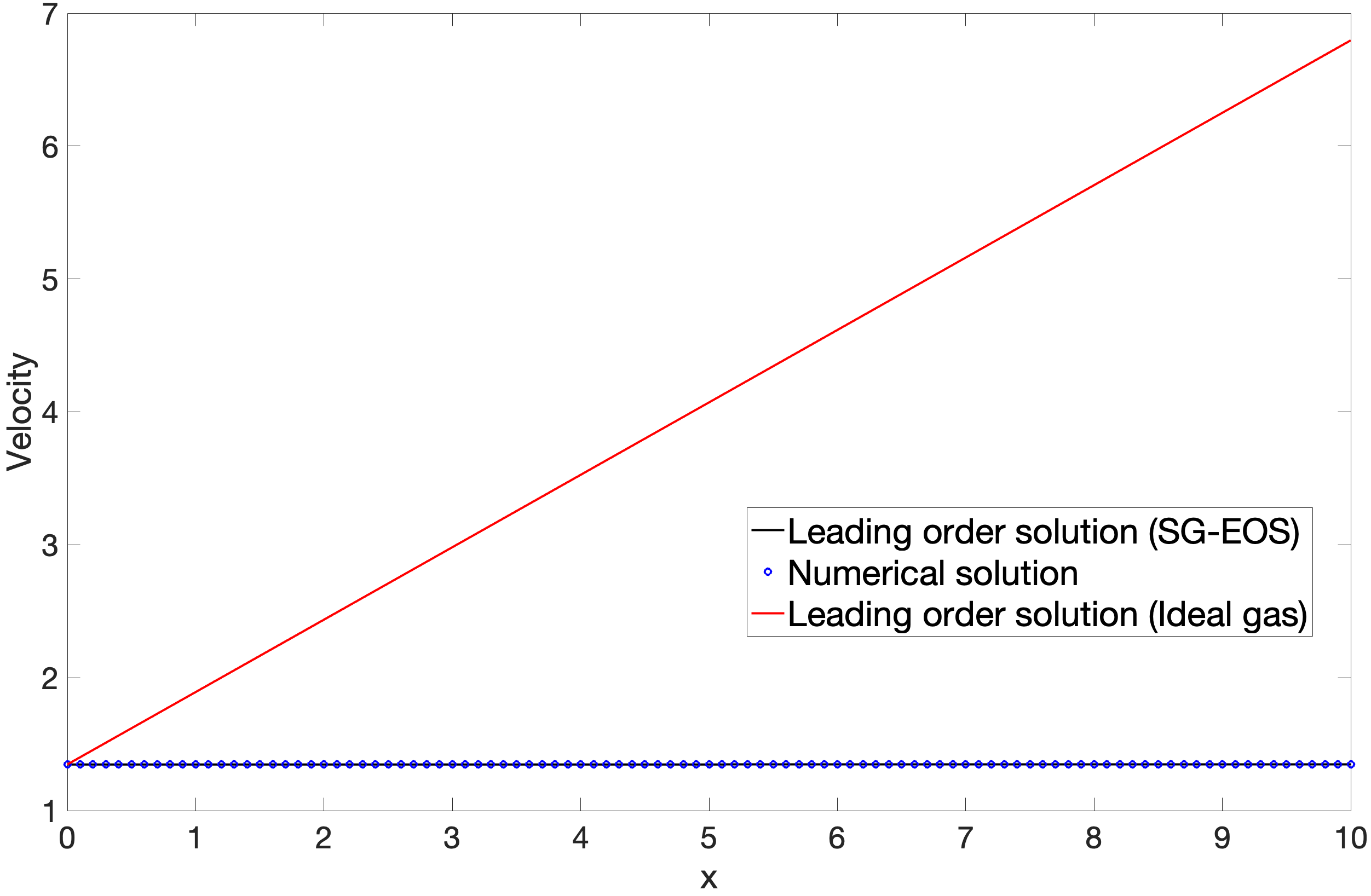}
    \end{subfigure}
    \caption{Open tube test case with the SG-EOS \eqref{eq:sg_eos}, results at $t = T_{f} = 7.47$. Left: density. Right: velocity. The continuous black line shows the leading order solution as $M \to 0$, the blue dots report the numerical results, while the red line shows the leading order solution as $M \to 0$ with the ideal gas employing $\gamma$ = 1.4.}
    \label{fig:open_tube_SG} 
\end{figure}

We consider now an extension of this test case for the SG-EOS \eqref{eq:sg_eos}. Equation \eqref{eq:energy_limit_incomp} reduces to
\begin{equation}
    \dive\bar{\mathbf{u}} = -\frac{1}{\gamma\left(\bar{p} + \pi_{\infty}\right)}\frac{d\bar{p}}{dt}.
\end{equation}
Hence, the velocity is still a linear function of space with a different time dependent slope with respect to that of the ideal gas. Analogous considerations hold for the continuity equation, which reduces to
\begin{equation}
    \frac{D\log\bar{\rho}}{Dt} = \frac{1}{\gamma\left(\bar{p} + \pi_{\infty}\right)}\frac{d\bar{p}}{dt} = \frac{1}{\gamma}\frac{d\log\left(\bar{p} + \pi_{\infty}\right)}{dt}.
\end{equation} 
We take $\gamma = 4.4, \pi_{\infty} = 6.8 \times 10^{3}$, and $q_{\infty} = 0$ in \eqref{eq:sg_eos}. The time step is not modified. Hence, the maximum advective Courant number is $C_{u} \approx 0.015$, while the maximum acoustic Courant number is $C \approx 19300$. Notice that, with this configuration, an explicit scheme would require a time step around 65000 times smaller to achieve a stable solution, yielding therefore a computational cost orders of magnitude larger. A comparison at the final time between the numerical results and the leading order solution for both the density and the velocity displays a good agreement for both profiles (Figure \ref{fig:open_tube_SG}). The leading order term solution as $M \to 0$ for the ideal gas with $\gamma = 1.4$, i.e. the previous configuration, has been included in Figure \ref{fig:open_tube_SG}. One can easily notice a visible difference in the behaviour of both density and velocity. In particular, considering the large value of $\pi_{\infty}$, the velocity field is almost constant (Figure \ref{fig:open_tube_SG}). Hence, if large amplitude pressure variations are considered, the limit regime depends on the equation of state and on its parameters and does not necessarily correspond to the incompressible Euler equations.  

\begin{figure}[h!]
    \centering
    \begin{subfigure}{0.475\textwidth}
	\centering
	\includegraphics[width = 0.9\textwidth]{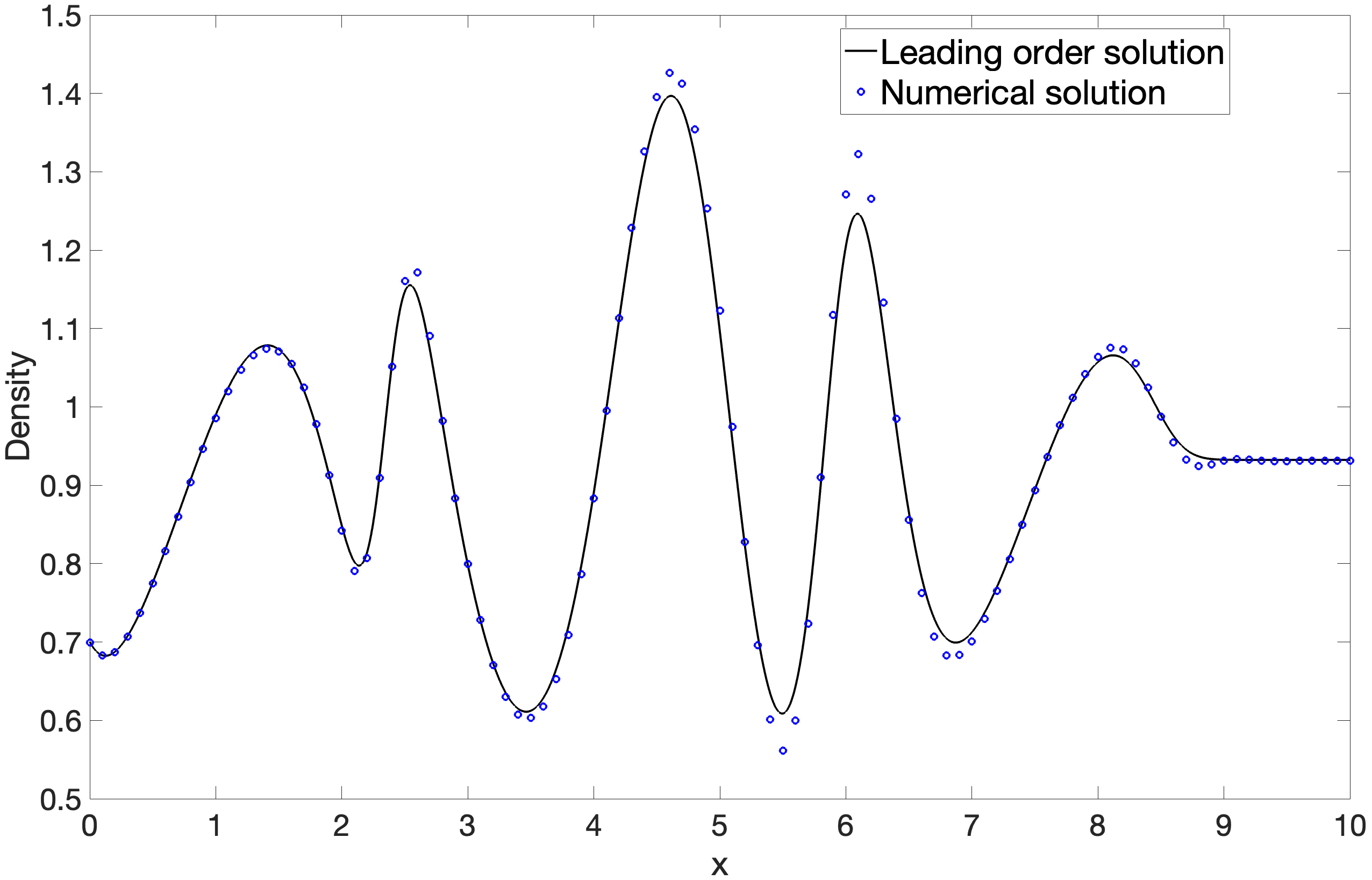}
    \end{subfigure}
    \begin{subfigure}{0.475\textwidth}
	\centering
	\includegraphics[width = 0.9\textwidth]{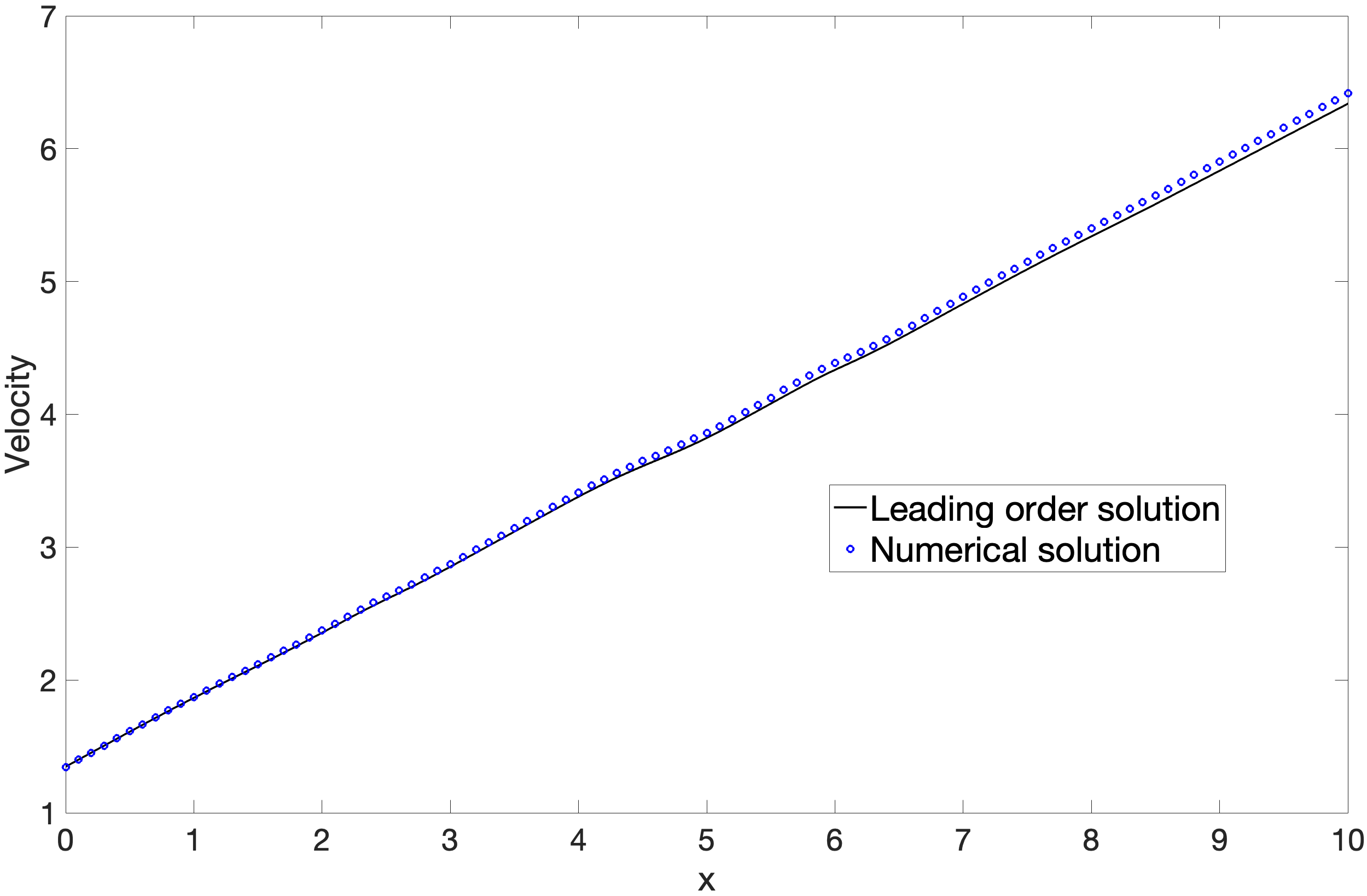}
    \end{subfigure}
    \caption{Open tube test case with the Peng-Robinson EOS \eqref{eq:general_cubic_eos}, results at $t = T_{f} = 7.47$. Left: density. Right: velocity. The continuous black line shows the leading order solution as $M \to 0$, whereas the blue dots report the numerical results.}
    \label{fig:open_tube_PR} 
\end{figure}

Finally, we consider the Peng-Robinson EOS \eqref{eq:general_cubic_eos}. The asymptotic analysis becomes much more involved. First of all, notice that
\begin{equation}
    \frac{\partial\bar{\rho}\bar{e}}{\partial\bar{p}} = \frac{1 - \bar{\rho}b}{\gamma - 1} \qquad \frac{\partial\bar{h}}{\partial\bar{\rho}} = - \frac{\gamma}{\gamma - 1}\frac{\bar{p}}{\bar{\rho}^{2}} + g\left(\bar{\rho}\right),
\end{equation}
where we have set
\begin{equation}
    g\left(\bar{\rho}\right) = \frac{a\left(1 - 2\bar{\rho}b\right)}{\left(\gamma - 1\right)\left(1 - \bar{\rho}br_{1}\right)\left(1 - \bar{\rho}br_{2}\right)} + \frac{ab\bar{\rho}\left(1 - \bar{\rho}b\right)\left(r_{1}\left(1 - \bar{\rho}br_{2}\right) + r_{2}\left(1 - \bar{\rho}br_{1}\right)\right)}{\left(\gamma - 1\right)\left(1 - \bar{\rho}br_{1}\right)^{2}\left(1 - \bar{\rho}br_{2}\right)^{2}} + \frac{a}{b}\frac{\partial U}{\partial\bar{\rho}}.
\end{equation}
Hence, \eqref{eq:energy_limit_incomp} reduces to
\begin{equation}\label{eq:energy_limit_PR}
    \left(-\frac{\gamma}{\gamma - 1}\bar{p} + \bar{\rho}^{2}g\left(\bar{\rho}\right)\right)\dive\bar{\mathbf{u}} = \frac{1 - \bar{\rho}b}{\gamma - 1}\frac{d\bar{p}}{dt},
\end{equation}
or, equivalently, to
\begin{equation}
    \dive\bar{\mathbf{u}} = -\frac{1 - \bar{\rho}b}{\gamma\bar{p} - \left(\gamma - 1\right)\bar{\rho}^{2}g\left(\bar{\rho}\right)}\frac{d\bar{p}}{dt}.
\end{equation}
Notice that, $\dive\bar{\mathbf{u}}$ is now a function of both space and time. Hence, in one space dimension, the velocity is no longer a linear profile. The continuity equation \eqref{eq:continuity_limit} reads as follows:
\begin{equation}
    \frac{D\log\bar{\rho}}{Dt} = -\dive\bar{\mathbf{u}} = \frac{1 - \bar{\rho}b}{\gamma\bar{p} - \left(\gamma - 1\right)\bar{\rho}^{2}g\left(\bar{\rho}\right)}\frac{d\bar{p}}{dt}.
\end{equation}
We take $\gamma = 1.4, a = 1$, and $b = 0.15$. The time step is $\Delta t = 9.3375 \times 10^{-4}$, yielding a maximum advective Courant number $C_{u} \approx 0.07$ and a maximum acoustic Courant number $C \approx 150$. A comparison at the final time between the numerical results and the leading order solution for both the density and the velocity shows a good agreement for both profiles (Figure \ref{fig:open_tube_PR}). The results are similar to those obtained with the ideal gas. Weakly non-ideal gas effects are present in particular between $x = 4$ and $x = 6$, namely in correspondence of the peak density.

\subsection{Gresho vortex}
\label{ssec:gresho}

\begin{figure}[h!]
    \centering
    \begin{subfigure}{0.475\textwidth}
	\centering
        \includegraphics[width = 0.9\textwidth]{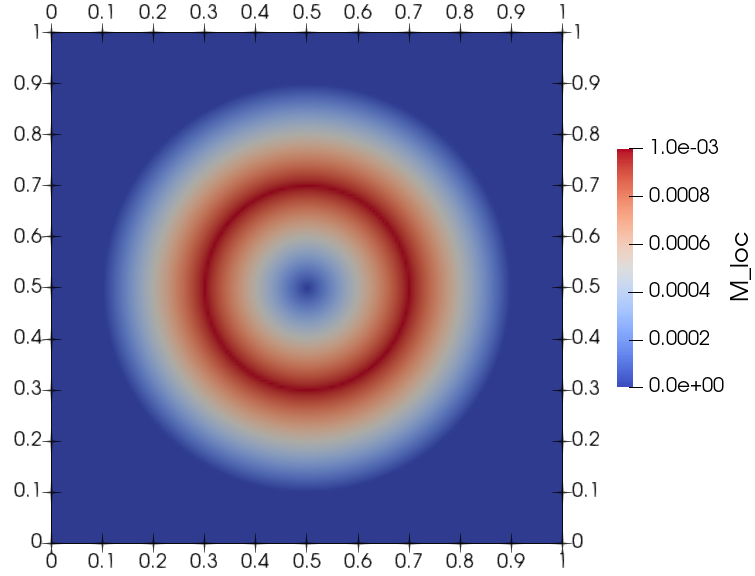}
    \end{subfigure}
    \begin{subfigure}{0.475\textwidth}
	\centering
        \includegraphics[width = 0.9\textwidth]{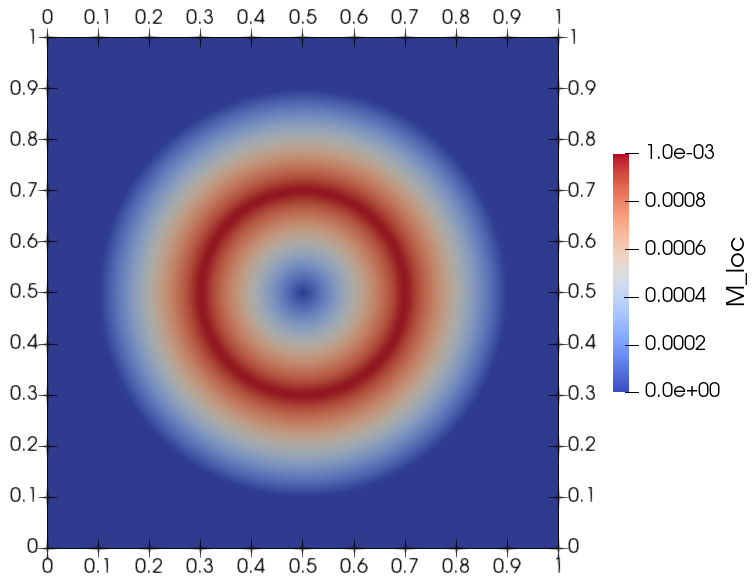}
    \end{subfigure}
    \begin{subfigure}{0.475\textwidth}
	\centering
        \includegraphics[width = 0.9\textwidth]{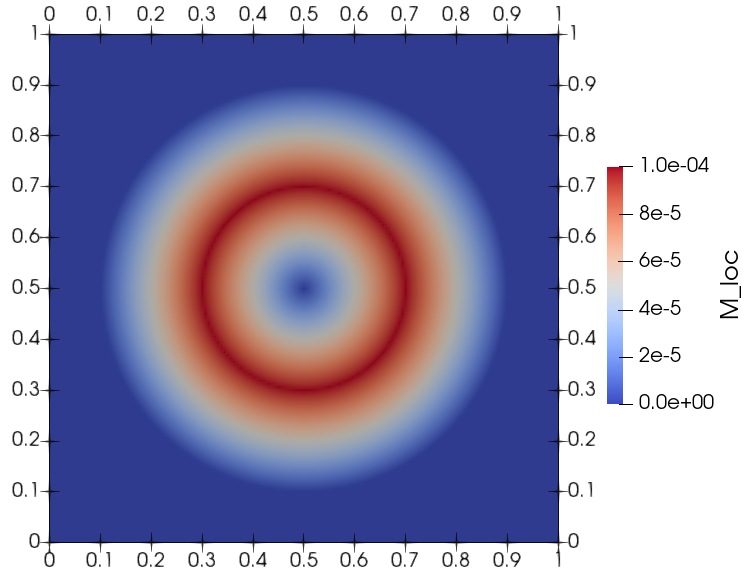}
    \end{subfigure}
    \begin{subfigure}{0.475\textwidth}
	\centering
        \includegraphics[width = 0.9\textwidth]{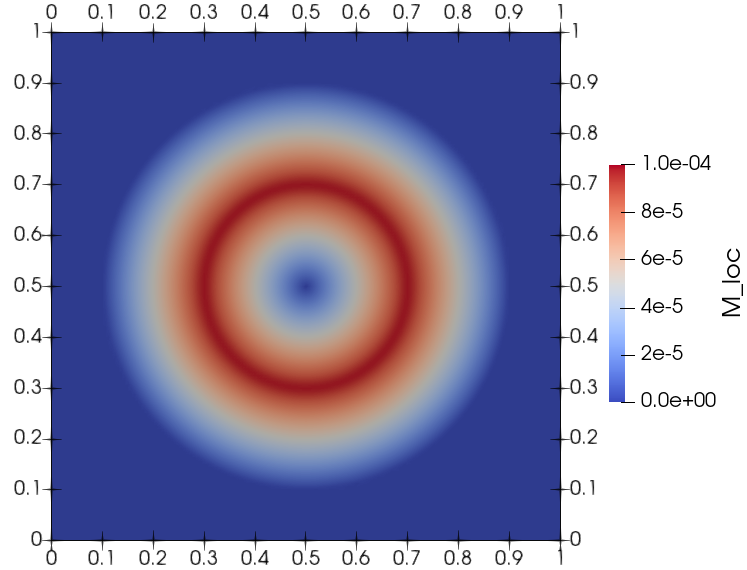}
    \end{subfigure}
    \caption{Gresho vortex test case with the ideal gas law \eqref{eq:ideal_gas}, comparison of local Mach number $M_{loc} = \frac{M\left|\mathbf{u}\right|}{c}$. From bottom to top: results at $M = 10^{-4}$ and $M = 10^{-3}$. From left to right: initial condition and results at $t = T_{f} = 3$, after three full rotations.}
    \label{fig:Gresho_IG}
\end{figure}

In this Section, we perform simulations of the so-called Gresho vortex \cite{gresho:1990, liska:2003}, which is a stationary solution of the incompressible Euler equations. The centrifugal force, indeed, is balanced by the gradient of the pressure. A rotating vortex is positioned at the center $\left(0.5, 0.5\right)$ of the computational domain $\Omega = \left(0,1\right)^{2}$. The initial conditions for dimensional variables read as follows:
\begin{subequations}
\begin{eqnarray}
    \rho\left(x,0\right) &=& 1 \qquad
    u\left(x,0\right) = -u_{\varphi}\sin\left(\varphi\right) \qquad
    v\left(x,0\right) = u_{\varphi}\cos\left(\varphi\right) 
    \label{eq:rho_u_init_Gresho} \\
    p\left(x,0\right) &=&
    \begin{cases}
        p_{0} + \frac{25}{2}\rho r^{2} \qquad &\text{if }0 \le \tilde{r} < 0.2 \\
        p_{0} + \frac{25}{2}\rho\tilde{r}^{2} + 4\rho\left(1 - 5\tilde{r} - \log(0.2) + \log(\tilde{r})\right) \qquad &\text{if }0.2 \le \tilde{r} < 0.4 \\
        p_{0} - \rho\left(2 - 4\log(2)\right) \qquad &\text{if } \tilde{r} \ge 0.4.
    \end{cases}\label{eq:pres_init_Gresho}
\end{eqnarray}
\end{subequations}
Here, $\tilde{r} = \sqrt{\left(x - 0.5\right)^{2} + \left(y - 0.5\right)^{2}}, \varphi = \arctan\left(\frac{y - 0.5}{x - 0.5}\right), p_{0} = \frac{\rho_{0}u_{\varphi,max}^{2}}{\gamma M^{2}}$, with $\rho_{0} = \SI{1}{\kilogram\per\meter\cubed}$ and $u_{\varphi,max} = \SI{1}{\meter\per\second}$ for $\tilde{r} = 0.2$. Finally, $u_{\varphi}$ is
\begin{equation}
    u_{\varphi} =
    \begin{cases}
	5\tilde{r} \qquad &\text{if }0 \le \tilde{r} < 0.2 \\
	2 - 5\tilde{r} \qquad &\text{if }0.2 \le \tilde{r} < 0.4 \\
	0 \qquad &\text{if } \tilde{r} \ge 0.4.
    \end{cases}
\end{equation}
Notice that, as discussed in \cite{happenhofer:2011}, the pressure $p_{0}$ is chosen in such a way that the maximum value of $\frac{\left|\mathbf{u}\right|}{c}$ matches $M$, so as to consider low Mach effects. We transform the initial conditions in non-dimensional quantities by using $\mathcal{R} = \SI{1}{\kilogram\per\meter\cubed}, \mathcal{L} = \SI{1}{\meter}$, and $\mathcal{U} = \SI[parse-numbers=false]{1}{\meter\per\second}$. Periodic boundary conditions are imposed for all the boundaries. We simulate the flow until $T_{f} = 3$, when three full rotations are completed. The computational grid is composed by $80 \times 80$ elements with polynomial degree $r = 2$, whereas the time step is $\Delta t = 2 \times 10^{-3}$, leading an advective Courant number $C_{u} \approx 0.32$. We consider $M = 10^{-3}$ and $M = 10^{-4}$. Hence, the acoustic Courant number is $C \approx 320$ for $M = 10^{-3}$ and $C \approx 3200$ for $M = 10^{-4}$. A comparison of the local Mach number $M_{loc} = \frac{M\left|\mathbf{u}\right|}{c}$ at initial time and the final time for the two tests shows that the numerical method accurately preserves the shape of the vortex (Figure \ref{fig:Gresho_IG}). We also monitor the behaviour over time of the kinetic energy, which should be conserved. Table \ref{tab:relative_ke_Gresho_IG} reports the total kinetic energy relative to the initial one after each rotation. The kinetic energy is conserved and these results compare very well with those presented in \cite{abbate:2019}, \cite{thomann:2019} where a loss of about 1.5 percent of the initial kinetic energy occurs after one rotation of the vortex. Analogous results are achieved for $M = 10^{-4}$. Hence, the preservation of the kinetic energy holds independently of the Mach number.

\begin{table}[pos=H]
    \centering
    \small
    \begin{tabularx}{0.45\columnwidth}{cccc}
	\toprule
	$M$ & $t = 1$ & $t = 2$ & $t = 3$ \\
	\midrule
	$10^{-3}$ & $0.999981$ & $0.999977$ & $0.999968$ \\
	\midrule
	$10^{-4}$ & $0.999981$ & $0.999977$ & $0.999968$ \\
	\bottomrule		
    \end{tabularx}
    \caption{Total kinetic energy relative to its initial value for different Mach numbers after each full rotation of the Gresho vortex with the ideal gas law \eqref{eq:ideal_gas}.}
    \label{tab:relative_ke_Gresho_IG}
\end{table}

\begin{figure}[h!]
    \centering
    \begin{subfigure}{0.475\textwidth}
	\centering
        \includegraphics[width = 0.9\textwidth]{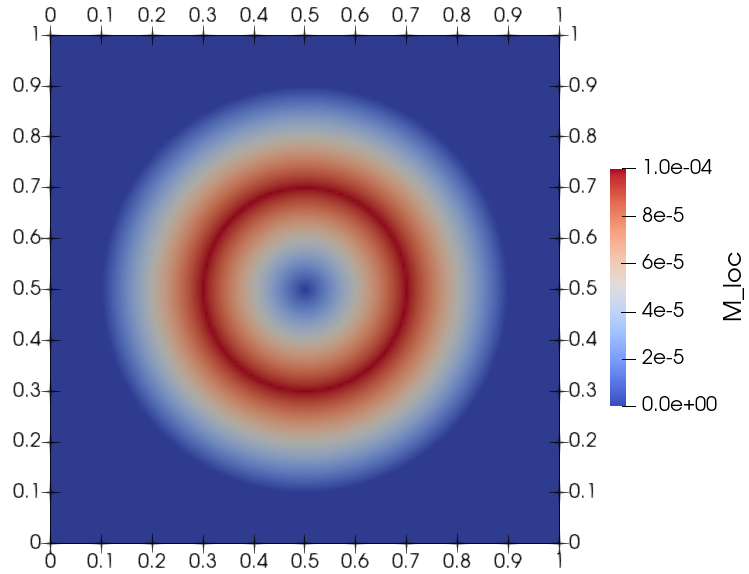}
    \end{subfigure}
    \begin{subfigure}{0.475\textwidth}
	\centering
        \includegraphics[width = 0.9\textwidth]{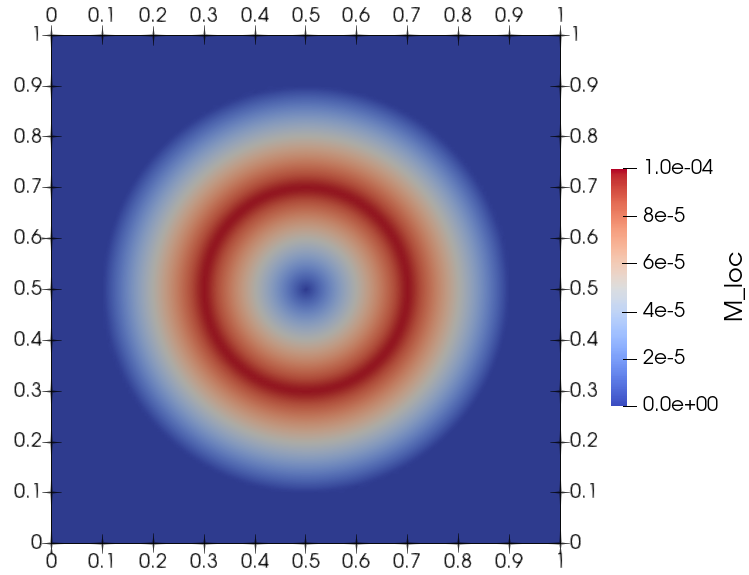}
    \end{subfigure}
    \caption{Gresho vortex test case with SG-EOS \eqref{eq:sg_eos}, comparison of local Mach number $M_{loc} = \frac{M\left|\mathbf{u}\right|}{c}$. Left: initial condition. Right: results at $t = T_{f} = 3$, after three full rotations.}
    \label{fig:Gresho_SG_EOS}
\end{figure}

We now adapt the standard Gresho vortex test case to a water flow. As discussed in \cite{abbate:2019}, it suffices to modify $p_{0}$ for the SG-EOS \eqref{eq:sg_eos} as follows:
\begin{equation}
    p_{0} = \frac{\rho_{0}u_{\varphi,max}^{2}}{\gamma M^{2}} - \pi_{\infty},
\end{equation}
with $M = 10^{-4}, \rho_{0} = \SI{1000}{\kilogram\per\meter\cubed}, \gamma = 4.4,$ and $\pi_{\infty} = \SI[parse-numbers=false]{6.8 \times 10^{8}}{\pascal}$. We also take $q_{\infty} = 0$ in \eqref{eq:sg_eos}. The initial density is now $\rho\left(x,0\right) = \rho_{0}$ and we employ $\mathcal{R} = \SI{1000}{\kilogram\per\meter\cubed}$ to compute the non-dimensional counter part of initial conditions \eqref{eq:rho_u_init_Gresho}-\eqref{eq:pres_init_Gresho}. A comparison of $M_{loc}$ between the initial and the final time shows that the shape of the vortex is accurately preserved also for a fluid with parameters corresponding to those of water (Figure \ref{fig:Gresho_SG_EOS}). Table \ref{tab:relative_ke_Gresho_SG} reports the total kinetic energy relative to the initial one after each rotation, from which we notice that the kinetic energy is conserved.

\begin{table}[h!]
    \centering
    \small
    \begin{tabularx}{0.45\columnwidth}{cccc}
	\toprule
	$M$ & $t = 1$ & $t = 2$ & $t = 3$ \\
	\midrule
	$10^{-4}$ & $0.999984$ & $0.999981$ & $0.999977$ \\
	\bottomrule		
    \end{tabularx}
    \caption{Total kinetic energy relative to its initial value after each full rotation of the Gresho vortex with the SG-EOS \eqref{eq:sg_eos}.}
    \label{tab:relative_ke_Gresho_SG}
\end{table}

\begin{figure}[h!]
    \centering
    \begin{subfigure}{0.475\textwidth}
	\centering
        \includegraphics[width = 0.9\textwidth]{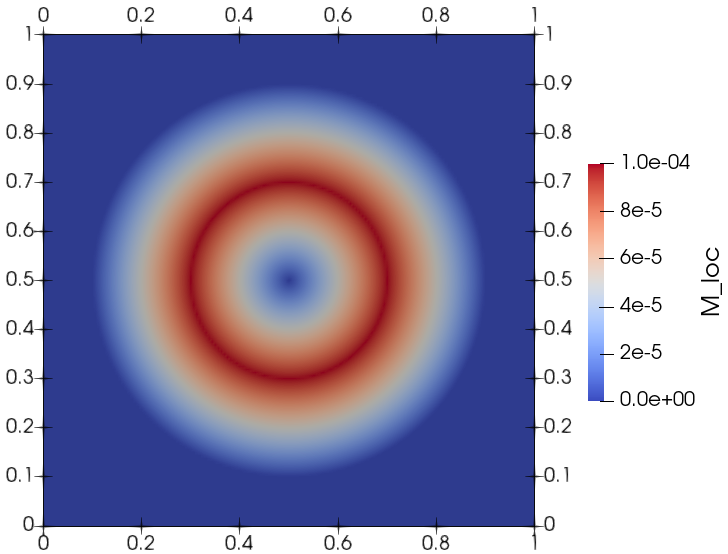}
    \end{subfigure}
    \begin{subfigure}{0.475\textwidth}
	\centering
        \includegraphics[width = 0.9\textwidth]{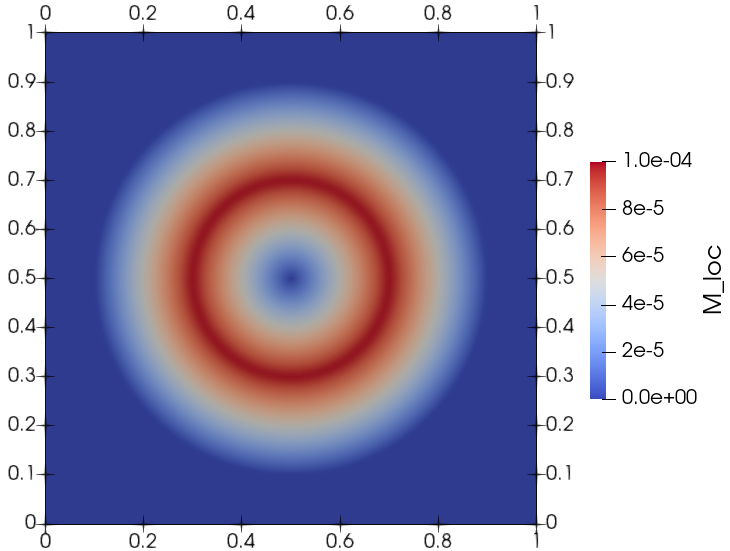}
    \end{subfigure}
    \caption{Gresho vortex test case with the Peng-Robinson EOS \eqref{eq:general_cubic_eos}, comparison of local Mach number $M_{loc} = \frac{M\left|\mathbf{u}\right|}{c}$. Left: initial condition. Top: results at $t = T_{f} = 3$, after three full rotations.}
    \label{fig:Gresho_PR_EOS}
\end{figure}

Finally, we consider a configuration of the Gresho vortex for the Peng-Robinson EOS \eqref{eq:general_cubic_eos}. The new expression of the background pressure $p_{0}$ reads as follows:
\begin{equation}
    p_{0} = \left[\frac{u_{\varphi,max}^{2}}{M^{2}} + f\left(\rho_{0}\right)\right]\frac{\rho_{0}\left(1 - \rho_{0}b\right)}{\gamma},
\end{equation}
with
\begin{equation}
    f\left(\rho_{0}\right) = \frac{a\rho_{0}}{1 - \rho_{0}b}\left(\frac{\frac{\partial U}{\partial\rho_{0}}}{b}\left(\gamma - 1\right) + \frac{1 - 2\rho_{0}b}{\left(1 - \rho_{0}br_{1}\right)\left(1 - \rho_{0}br_{2}\right)}\right) + ab\rho_{0}^{2}\frac{r_{1}\left(1 - \rho_{0}br_{2}\right) + r_{2}\left(1 - \rho_{0}br_{1}\right)}{\left(1 - \rho_{0}br_{1}\right)^{2}\left(1 - \rho_{0}br_{2}\right)^{2}}.
\end{equation}
We take $\gamma = 1.4, \rho_{0} = \SI{1}{\kilogram\per\meter\cubed}, a = \SI{500}{\meter\tothe{5}\per\second\squared\per\kilogram}$, and $b = \SI[parse-numbers=false]{10^{-3}}{\meter\cubed\per\kilogram}$. Finally, we consider $\rho\left(x,0\right) = \SI{1}{\kilogram\per\meter\cubed}$ and this value is also employed to compute non-dimensional quantities. Figure \ref{fig:Gresho_PR_EOS} shows a comparison of $M_{loc}$ between the initial and the final time, while Table \ref{tab:relative_ke_Gresho_PR} reports the total kinetic energy relative to the initial one after each rotation, from which we notice that the kinetic energy is conserved. The same considerations done for the ideal gas law \eqref{eq:ideal_gas} and for the SG-EOS \eqref{eq:sg_eos} are therefore valid also for this particularly challenging and complex equation of state.

\begin{table}[h!]
    \centering
    \small
    \begin{tabularx}{0.45\columnwidth}{cccc}
	\toprule
	$M$ & $t = 1$ & $t = 2$ & $t = 3$ \\
	\midrule
	$10^{-4}$ & $0.999981$ & $0.999977$ & $0.999968$ \\
	\bottomrule		
    \end{tabularx}
    \caption{Total kinetic energy relative to its initial value after each full rotation of the Gresho vortex with the Peng-Robinson EOS \eqref{eq:general_cubic_eos}.}
    \label{tab:relative_ke_Gresho_PR}
\end{table}

\subsection{Baroclinic vorticity generation problem}
\label{ssec:baroclinic}

We now consider a test case proposed in \cite{geratz:1998} and discussed also in \cite{klein:2002, noelle:2014}, which consists of a right-going acoustic wave crossing a density fluctuation in the vertical direction. This test case illustrates the nontrivial interaction between large-scale acoustic waves and small-scale density fluctuations. Following the discussion in \cite{klein:1995, klein:2002} and in Section \ref{ssec:ap_num_two_scale_length}, we notice that the mass and momentum balance for \eqref{eq:euler_adim_ap_two_scale} read as follows:
\begin{eqnarray}
    \frac{\partial\bar{\rho}}{\partial t} &=& 0 \\
    \frac{\partial\bar{\rho}\bar{\mathbf{u}}}{\partial t} + \gradxi p^{'} &=& 0.
\end{eqnarray}
Suppose now that two neighbouring mass elements characterized by densities $\rho_{1}$ and $\rho_{2}$, with $\rho_{1} \neq \rho_{2}$ as in the present test case, are accelerated by a common large-scale acoustic pressure gradient. Since the time derivative of the momentum is the same for both mass elements, their velocities must differ by a factor of $\frac{\rho_{2}}{\rho_{1}}$. As a consequence of different accelerations, vorticity is generated. This phenomenon is also known as \textit{baroclinic effect} and it is the result of mutual interaction between the quasi-incompressible small-scale and the large-scale acoustic flow. Indeed, baroclinic instabilities are well known to play a major role in large scale atmospheric dynamics \cite{charney:1947, jablonowski:2006}, as well as in other areas of compressible fluid dynamics, see, e.g., \cite{brouillette:2002}.

\begin{figure}[h!]
    \centering
    \begin{subfigure}{\textwidth}
	\centering
        \includegraphics[width=0.9\textwidth]{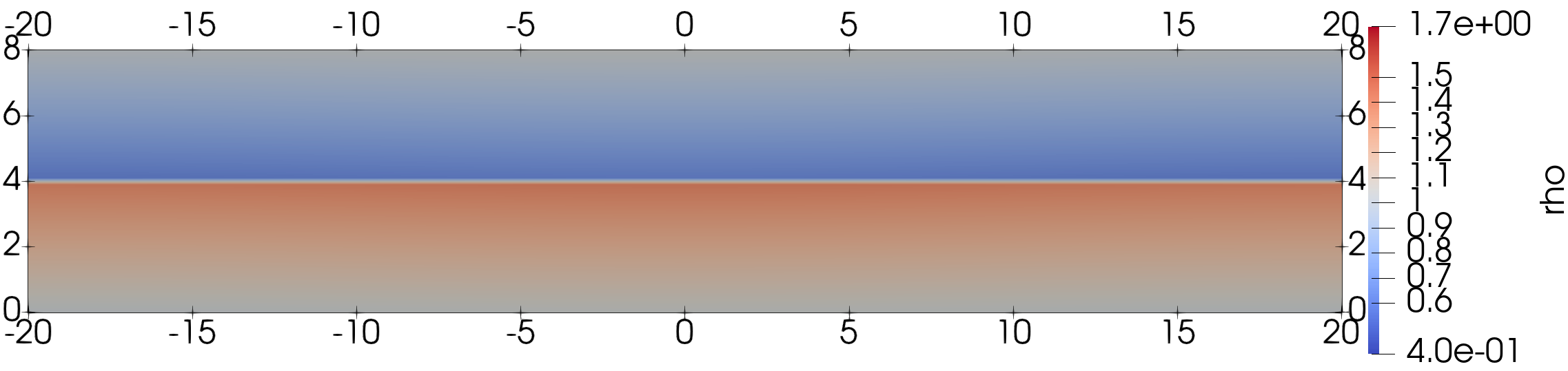}
    \end{subfigure}
    \begin{subfigure}{\textwidth}
	\centering
        \includegraphics[width=0.9\textwidth]{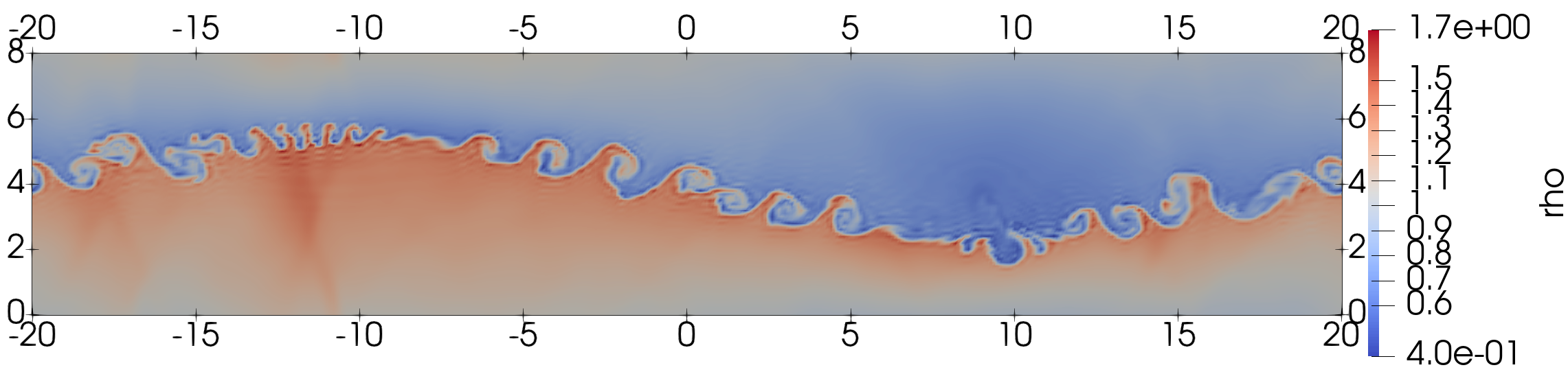}
    \end{subfigure}
    \caption{Baroclinic vorticity generation with the ideal gas law \eqref{eq:ideal_gas}, contour plot of the density. Top: $t = 0$. Bottom: $t = T_{f} = 16$.}
    \label{fig:baroclinic_density_IG}
\end{figure}

We first consider this test assuming that the ideal gas law \eqref{eq:ideal_gas} holds. The computational domain is $\Omega = \left(-L, L\right) \times \left(0, \frac{2}{5}L\right)$. Following \cite{klein:2002, noelle:2014}, we set $M = 5 \times 10^{-2}$ and we take $L = \frac{1}{M}$. The initial conditions read as follows:
\begin{subequations}
\begin{eqnarray}
    \bar{\rho}\left(\mathbf{x},0\right) &=& \bar{\rho}_{0} + M\rho_{0}^{'}\left(1 + \cos\left(\frac{\pi x}{L}\right)\right) + \Phi\left(y\right) \\
    u\left(\mathbf{x},0\right) &=& \bar{u}_{0}\left(1 + \cos\left(\frac{\pi x}{L}\right)\right) \\
    v\left(\mathbf{x},0\right) &=& 0 \\
    p\left(\mathbf{x},0\right) &=& \bar{p}_{0} + Mp_{0}^{'}\left(1 + \cos\left(\frac{\pi x}{L}\right)\right),
\end{eqnarray}
\end{subequations}
with $\bar{\rho}_{0} = 1, \rho_{0}^{'} = 0.2, \bar{p}_{0} = 1. p_{0}^{'} = \gamma$, and $\bar{u}_{0} = \sqrt{\gamma}$. The function $\Phi\left(y\right)$ is defined by
\begin{equation}
    \Phi\left(y\right) =
    \begin{cases}
        \rho_{2}\frac{y}{\frac{2}{5}L} \qquad &\text{if } y \le \frac{1}{5}L - \varepsilon \\
        \rho_{2}\left(\frac{y}{\frac{2}{5}L} - \frac{1}{2}\right) - 0.4 \qquad &\text{if } y \ge \frac{1}{5}L + \varepsilon \\
        \rho_{2}\frac{\frac{1}{5}L - \varepsilon}{\frac{2}{5}L} + \frac{1}{2\varepsilon}\left(\rho_{2}\left(\frac{\frac{1}{5}L + \varepsilon}{\frac{2}{5}L} - \frac{1}{2}\right) - 0.4 - \rho_{2}\frac{\frac{1}{5}L - \varepsilon}{\frac{2}{5}L}\right)\left(y - \frac{1}{5}L + \varepsilon\right) \qquad &\text{otherwise},
    \end{cases}
\end{equation}
where $\rho_{2} = 0.8$ and $\varepsilon = 10^{-2}$. Notice that, unlike in \cite{klein:2002, noelle:2014}, the function $\Phi$ is regularized to obtain a continuous profile. Periodic boundary conditions are prescribed, whereas the final time is $T_{f} = 16$. The computational grid is composed by $200 \times 40$ elements with $r = 2$. The time step is $\Delta t = 4 \times 10^{-3}$, yielding a maximum advective Courant number $C_{u} \approx 0.27$ and a maximum acoustic Courant number $C \approx 2.5$. Figure \ref{fig:baroclinic_density_IG} shows a comparison of the density between the initial and the final time. The initial density profile consists of two layers with different acceleration. Hence, a rotational motion is induced along the separating layer and a Kelvin-Helmholtz instability develops. 

\begin{figure}[h!]
    \centering
    \includegraphics[width = 0.9\textwidth]{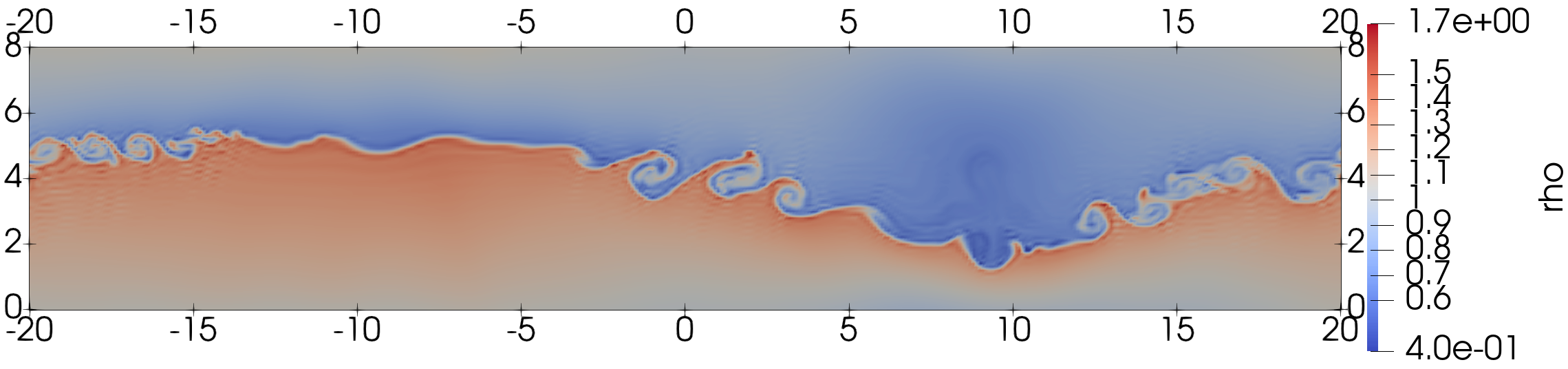}
    \caption{Baroclinic vorticity generation with the SG-EOS \eqref{eq:sg_eos}, contour plot of the density at $t = T_{f} = 16$.}
    \label{fig:baroclinic_density_SG}
\end{figure}

We then replicate the test considering the SG-EOS \eqref{eq:sg_eos}. We take $\gamma = 4.4, \pi_{\infty} = 6.8 \times 10^{-3}$, and $q_{\infty} = 0$. The same initial conditions of the configuration with the ideal gas law are employed. The maximum acoustic Courant number is $C \approx 3.5$, whereas the maximum advection Courant number $C_{u} \approx 0.23$. One can easily notice that the development of the Kelvin-Helmholtz instability depends on the EOS and on the fluid parameters (Figure \ref{fig:baroclinic_density_SG}).

\section{Conclusions}
\label{sec:conclu}

We have presented the asymptotic-preserving (AP) analysis of a general class of IMEX-RK schemes for the time discretization of the compressible Euler equations. Based on the results of \cite{casulli:1984, dumbser:2016b}, these approaches consider an implicit coupling between the momentum and the energy balance, while treating the density explicitly. Third order and fourth order time discretization schemes, in combination with a Discontinuous Galerkin (DG) for the space discretization, have been employed for numerical simulations. The AP property of the proposed method is valid for a general equation of state as well as for two length scales models. A number of classical benchmarks for ideal gases and their non trivial extension for equations of state of real gases, in particular for the general cubic equation of state, validate the proposed method in the low Mach number regime and in the limit of incompressible flows.
In particular, in spite of the use of quadrilateral meshes, the proposed method yields correct results for Mach number values that are typical of fluids, such as water, usually modelled as incompressible. Notice that no operator splitting, flux splitting or relaxation techniques have been employed, differently from the approaches proposed, e.g., in \cite{abbate:2019, chalons:2013, chalons:2016, klein:1995, noelle:2014}. In future work, we aim to consider gravity effects, so as to perform an asymptotic analysis in the limit of low Froude numbers, and to consider an extension to two-phase flows. Moreover, as already mentioned at the end of Section \ref{ssec:isentropic_vortex}, we aim to analyze more in detail the spatial discretization and, in particular, the use of compatible finite elements, as recently done in \cite{zampa:2025}, and entropy-stable schemes.

\section*{Acknowledgements}

We thank the three anonymous reviewers and the Associate Editor who handled the paper for their very useful and constructive comments and remarks, which have greatly helped in improving the quality of the presentation of our results. G.O. is part of the INdAM-GNCS National Research Group. The simulations have been partly run at CINECA thanks to the computational resources made available through the ISCRA-C project FEM-GPU - HP10CQYKJ1. We acknowledge the CINECA award for the availability of high-performance computing resources and support.

\appendix

\section{Eigenvalues of the implicit and explicit part}
\label{app:eigenvalues}

In this Appendix, we analyze the eigenvalues for the Euler equations \eqref{eq:euler_adim}. More specifically, we compute the eigenvalues for the two subsystems obtained considering the IMEX approach described in Section \ref{sec:ap_num}. For the sake of simplicity, we focus on 1D case, so that the equations can be written as follows:
\begin{eqnarray}
    \frac{\partial\rho}{\partial t} + \frac{\partial q}{\partial x} &=& 0 \nonumber \\
    \frac{\partial q}{\partial t} + \frac{\partial}{\partial x}\left(\frac{q^{2}}{\rho}\right) + \frac{1}{M^{2}}\frac{\partial p}{\partial x} &=& 0 \\
    \frac{\partial\hat{E}}{\partial t} + \frac{\partial hq}{\partial x} + \frac{1}{2}M^{2}\frac{\partial}{\partial x}\left(\frac{q^{3}}{\rho^{2}}\right) &=& 0, \nonumber
\end{eqnarray}
with $q = \rho u$ and $\hat{E} = \rho E$. Hence, considering the time discretization reported in Section \ref{sec:ap_num}, the system can be written in the following quasi-linear form:
\begin{equation}\label{eq:euler_adim_quasi_linear_orig}
    \frac{\partial\tilde{\mathbf{W}}}{\partial t} + \tilde{\mathbf{A}}_{I}\frac{\partial\tilde{\mathbf{W}}}{\partial x} + \tilde{\mathbf{A}}_{E}\frac{\partial\tilde{\mathbf{W}}}{\partial x} = \mathbf{0},
\end{equation}
with
$$\tilde{\mathbf{W}} =
    \begin{bmatrix}
	\rho \\
	q \\
	\hat{E}
    \end{bmatrix}
    \qquad
    \tilde{\mathbf{A}}_{I} = 
    \begin{bmatrix}
	0 & 0 & 0 \\
        \frac{1}{M^{2}}\frac{\partial p}{\partial\rho} & \frac{1}{M^{2}}\frac{\partial p}{\partial q} & \frac{1}{M^{2}}\frac{\partial p}{\partial\hat{E}} \\    
        q\frac{\partial h}{\partial\rho} & q\frac{\partial h}{\partial q} + h & q\frac{\partial h}{\partial\hat{E}}
    \end{bmatrix} 
    \qquad
    \tilde{\mathbf{A}}_{E} = 
    \begin{bmatrix}
	0 & 1 & 0 \\
	-u^{2} & 2u & 0 \\
	-M^{2}u^{3} & \frac{3}{2}M^{2}u^{2} & 0
    \end{bmatrix}.$$
Here, $\tilde{\mathbf{A}}_{I}$ and $\tilde{\mathbf{A}}_{E}$ denote matrices related to the fluxes discretized implicitly and explicitly, respectively. After some manipulations (see \cite{orlando:2022b}), we can rewrite \eqref{eq:euler_adim_quasi_linear_orig} as follows:
\begin{equation}\label{eq:euler_adim_quasi_linear}
    \frac{\partial\mathbf{W}}{\partial t} + \mathbf{A}_{I}\frac{\partial\mathbf{W}}{\partial x} + \mathbf{A}_{E}\frac{\partial\mathbf{W}}{\partial x} = \mathbf{0},
\end{equation}
with
$$
\mathbf{W} =
\begin{bmatrix}
    \rho \\ 
    u \\
    p
\end{bmatrix} \qquad
\mathbf{A}_{I} = 
\begin{bmatrix}
    0 & 0 & 0 \\
    0 & 0 & \frac{1}{\rho M^{2}} \\
    0 & \frac{\frac{p}{\rho} - \rho\frac{\partial e}{\partial\rho}}{\frac{\partial e}{\partial p}} & u
\end{bmatrix} \ \ \ \
\mathbf{A}_{E} = 
\begin{bmatrix}
    u & \rho & 0 \\
    0 & u & 0 \\
    0 & 0 & 0
\end{bmatrix}.
$$
The eigenvalues of $\mathbf{A}_{I}$ are
$$\frac{u}{2} - \sqrt{\frac{c^{2}}{M^{2}} + \frac{u^{2}}{4}} \qquad 0 \qquad \frac{u}{2} + \sqrt{\frac{c^{2}}{M^{2}} + \frac{u^{2}}{4}},$$
where the expression of the speed of sound $c$ is reported in \eqref{eq:speed_sound}, while the eigenvalues of $\mathbf{A}_{E}$ are
$$0 \qquad u \qquad u.$$
The eigenvalues of $\mathbf{A}_{E}$ are always real and the subsystem discretized explicitly does not take into account any acoustic effect. However, the subsystem is only weakly hyperbolic, since $\mathbf{A}_{E}$ is not diagonalizable. This is related to the fact that the terms treated explicitly in the continuity equation and in the momentum balance form the well-known pressureless gas dynamics system \cite{bouchut:2003, leveque:2002}. Since this system is weakly hyperbolic, delta-shocks can develop and the vacuum state can occur, yielding an expansion which propagates at infinite velocity. Nevertheless, the vacuum state cannot form spontaneously and we need to start from the vacuum to obtain the infinite velocity expansion \cite{bouchut:2003}. Moreover, in the case of regular solutions, as we are mainly interested in this work, the momentum equation decouples from the continuity equation and reduces to the Burgers' equation \cite{leveque:2002}. Hence, the velocity field can be computed solving the Burgers' equation and the continuity equation reduces to
\begin{equation}
    \frac{\partial\rho}{\partial t} + u\frac{\partial\rho}{\partial x} = -\rho\frac{\partial u}{\partial x},
\end{equation}
which is an evolution equation for $\rho$ along the characteristics for which the advecting field $u$ and the source term contribution $\frac{\partial u}{\partial x}$ are known.

\bibliographystyle{cas-model2-names}
\bibliography{AP_Euler.bib}
	
\end{document}